\input amstex \documentstyle{amsppt} 
\loadbold
\let\bk\boldkey
\let\bs\boldsymbol
\loadeusm \let\scr\eusm
\loadeurm 
\font\Rm=cmr12

\font\Rrm=cmr17

\font\scap=eusb10 
\define\cind{\text{\it c-\/\rm Ind}} 
\define\ind{\text{\it ind}} 
\define\Ind{\text{\rm Ind}}
\define\Aut#1#2{\text{\rm Aut}_{#1}(#2)}
\define\End#1#2{\text{\rm End}_{#1}(#2)}
\define\Hom#1#2#3{\text{\rm Hom}_{#1}({#2},{#3})} 
\define\GL#1#2{\roman{GL}_{#1}(#2)}
\define\M#1#2{\roman M_{#1}(#2)}
\define\Gal#1#2{\text{\rm Gal\hskip.5pt}(#1/#2)}
\define\N#1#2{\roman N_{#1/#2}}
\define\nrd#1{\roman{Nrd}_{#1}}
\define\Scr#1{\text{\scap #1}} 
\define\G#1#2{\scr G_{#1}(#2)}
\define\A#1#2{\scr A_{#1}(#2)} 
\define\Ao#1#2{\scr A^0_{#1}(#2)} 
\define\T#1#2{\scr T_{#1}(#2)}
\define\Tt#1#2{\scr T(#1;#2)} 
\define\LR#1{{#1}^\roman{ell}_\roman{reg}}
\define\SR#1{{#1}^\roman{ss}_\roman{reg}} 
\define\ags#1{\langle{#1}\rangle} 
\define\Ags#1{\big\langle{#1}\big\rangle} 
\define\bsi{\bs\iota} 
\define\Asq#1#2{\scr A^{\sssize\square}_{#1}(#2)} 
\define\Gsq#1#2{\scr G^{\sssize\square}_{#1}(#2)} 
\define\Go#1#2{\scr G^0_{#1}(#2)} 
\define\upr#1#2{{}^{#1\!}{#2}} 
\define\pre#1#2{{}_{#1}{#2}} 
\font\Symb=msam10 scaled \magstep2 
\define\Boxplus{{\Symb \char'001}}
\define\sqsum#1{\underset{#1}\to{\text{\Boxplus}}} 
\define\wt#1{\widetilde{#1}} 
\define\wG#1{\widehat{\text{\rm GL}}_{#1}} 
\let\ge\geqslant
\let\le\leqslant
\let\eps\epsilon
\let\ups\upsilon
\let\vD\varDelta 
\let\ve\varepsilon 
\let\vf\varphi 
\let\vF\varPhi
\let\vG\varGamma
\let\vk\varkappa
\let\vL\varLambda 
\let\vO\varOmega 
\let\vU\varUpsilon 
\let\vp\varpi

\let\vS\varSigma
\let\vt\vartheta
\let\vT\varTheta
 
\define\wP#1{\widehat{\scr P}_{#1}} 
\define\wI#1{\widehat{\scr I}_{#1}} 
\define\wW#1{\widehat{\scr W}_{#1}} 
\define\Gg#1#2#3{\scr G^0_{#1}(#2;#3)} 
\define\Aa#1#2#3{\scr A^0_{#1}(#2;#3)} 
\hfuzz=1.25pt
\topmatter 
\title\nofrills 
\hfill\hbox{\rm To appear in Memoirs Amer\. Math\. Soc\.}\\
\Rrm \vphantom{X} \\ 
To an effective local Langlands Correspondence 
\endtitle 
\rightheadtext{Effective Langlands correspondence} 
\author 
Colin J. Bushnell and Guy Henniart 
\endauthor 
\leftheadtext{C.J. Bushnell and G. Henniart}
\affil 
King's College London and Universit\'e de Paris-Sud 
\endaffil 
\address 
King's College London, Department of Mathematics, Strand, London WC2R 2LS, UK. 
\endaddress
\email 
colin.bushnell\@kcl.ac.uk 
\endemail
\address 
Universit\'e de Paris-Sud, Laboratoire de Math\'ematiques d'Orsay,
Orsay Cedex, F-91405; CNRS, Orsay cedex, F-91405. 
\endaddress 
\email 
Guy.Henniart\@math.u-psud.fr 
\endemail 
\date March 2011, corrected August 2012 \enddate 
\thanks 
Much of the work in this programme was carried out while the first-named author was visiting, and partly supported by, l'Universit\'e de Paris-Sud. 
\endthanks 
\abstract 
Let $F$ be a non-Archimedean local field. Let $\Cal W_F$ be the Weil group of $F$ and $\Cal P_F$ the wild inertia subgroup of $\Cal W_F$. Let $\widehat{\Cal W}_F$ be the set of equivalence classes of irreducible smooth representations of $\Cal W_F$. Let $\Cal A^0_n(F)$ denote the set of equivalence classes of irreducible cuspidal representations of $\roman{GL}_n(F)$ and set $\widehat{\roman{GL}}_F = \bigcup_{n\ge1} \Cal A^0_n(F)$. If $\sigma\in \widehat{\Cal W}_F$, let $\upr L\sigma \in \widehat{\roman{GL}}_F$ be the cuspidal representation matched with $\sigma$ by the Langlands Correspondence. If $\sigma$ is totally wildly ramified, in that its restriction to $\Cal P_F$ is irreducible, we treat $\upr L\sigma$ as known. From that starting point, we construct an explicit bijection $\Bbb N:\widehat{\Cal W}_F \to \widehat{\roman{GL}}_F$, sending $\sigma$ to $\upr N\sigma$. We compare this ``na\"\i ve correspondence'' with the Langlands correspondence and so achieve an effective description of the latter, modulo the totally wildly ramified case. A key tool is a novel operation of ``internal twisting'' of a suitable representation $\pi$ (of $\Cal W_F$ or $\roman{GL}_n(F)$) by tame characters of a tamely ramified field extension of $F$, canonically associated to $\pi$. We show this operation is preserved by the Langlands correspondence. 
\endabstract 
\endtopmatter 
\document \nopagenumbers 
We consider the local Langlands correspondence for the general linear group $\GL nF$ over a non-Archimedean local field $F$ of residual characteristic $p$. 
\subhead 
1 
\endsubhead 
Let $\scr W_F$ be the Weil group of $F$ relative to a chosen separable algebraic closure $\bar F/F$. For each integer $n\ge1$, let $\Go nF$ be the set of equivalence classes of smooth, complex representations of $\scr W_F$ which are irreducible of dimension $n$. On the other side, let $\Ao nF$ be the set of equivalence classes of smooth complex representations of $\GL nF$ which are irreducible and cuspidal. The Langlands correspondence provides a canonical bijection $\Go nF \to \Ao nF$, which we denote $\Bbb L:\sigma\mapsto \upr L\sigma$. 
\par 
The irreducible representations of $\scr W_F$ can be difficult to describe individually in any detail. On the other hand, the irreducible cuspidal representations of $\GL nF$ are the subject of an explicit classification theory \cite{18}. One cannot avoid asking how this wealth of detail on the one side gets translated, via the Langlands correspondence, into concrete information on the other. Our work \cite{10}, \cite{11}, \cite{13} on {\it essentially tame\/} representations revealed the connection in that case to be clear and intuitively obvious (up to a complicated, but functionally minor, technical correction). The aim of this paper is to bring a corresponding degree of clarity to the general case: while the objects under consideration are inherently more complicated, the transparency of the essentially tame case persists to a surprising degree.  
\subhead 
2 
\endsubhead  
The set $\Go1F$ (resp\. $\Ao1F$) is the group of smooth characters $\scr W_F \to \Bbb C^\times$ (resp\. $F^\times \to \Bbb C^\times$), and the Langlands correspondence $\sigma\mapsto \upr L\sigma$ is local class field theory. In the general case $n\ge2$, the correspondence is specified \cite{25} in terms of its behaviour relative to $L$-functions and local constants of pairs of representations, as in \cite{31}, \cite{36}. Its existence is established, in these terms, in \cite{34}, \cite{22}, \cite{26}, \cite{27}. 
\par 
Local constants of pairs of representations remain resistant to explicit computation, so the standard characterization of the correspondence is not helpful as a tool for elucidating it. On the other hand, it does imply that the correspondence is compatible, in a straightforward sense \cite{17}, with base change and automorphic induction for cyclic base field extensions \cite{1}, \cite{28}, \cite{30}. When the cyclic base field extension is also tamely ramified, base change and automorphic induction are, to a useful degree, expressible in terms of the classification of cuspidal representations \cite{3}, \cite{4}, \cite{8}, \cite{9}. This paper is a systematic exploitation of that connection. 
\subhead 
3  
\endsubhead 
There is a basic case, outside the scope of this paper. An irreducible smooth representation $\sigma$ of $\scr W_F$ is called {\it totally wildly ramified\/} if its restriction $\sigma|_{\scr P_F}$ to the wild inertia subgroup $\scr P_F$ of $\scr W_F$ is irreducible. In particular, the dimension of such a representation is a power $p^r$ of $p$, for some integer $r\ge0$. Little is known of the detailed structure of such representations or their behaviour relative to the Langlands correspondence. The case $r=1$ has been worked through in detail in \cite{32} when $p=2$ (but see also \cite{12}), and in \cite{23} when $p=3$. For general $p$, see \cite{33} and \cite{35}, plus \cite{4} for a final detail. Some first steps towards the general case $r\ge2$ are made in \cite{4}, \cite{5}, \cite{6}. We do not pursue that matter any further here. We proceed knowing that $\upr L\sigma$ is defined, but otherwise treat it as effectively unknown. This degree of ignorance has little practical effect, for simple reasons discussed in 7.8 below. 
\par 
From that initial position, we use explicit, elementary, techniques to construct a ``na\"\i ve correspondence'' $\Bbb N:\Go nF \overset\approx\to\longrightarrow \Ao nF$, sending $\sigma$ to $\upr N\sigma$. In particular, $\upr N\sigma = \upr L\sigma$ when $\sigma$ is totally wildly ramified. We estimate, closely and uniformly, the difference between the na\"\i ve correspondence and the Langlands correspondence. We thereby reveal a family of transparent and intuitively obvious connections. We simultaneously uncover a novel (and useful) arithmetic property of the Langlands correspondence. 
\subhead 
4 
\endsubhead 
We give an overview of the ideas leading to our main results. We start on the Galois side, with a description of the irreducible smooth representations of $\scr W_F$. Since it determines the flavour of all that follows, we give a detailed summary now but refer to \S1 below for proofs. 
\par 
Let $\wP F$ denote the set of equivalence classes of irreducible smooth representations of $\scr P_F$. The group $\scr W_F$ acts on $\wP F$ by conjugation. Taking $\alpha\in \wP F$, the $\scr W_F$-isotropy group of $\alpha$ is of the form $\scr W_E$, for a finite, tamely ramified field extension $E/F$. We use the notation $E = Z_F(\alpha)$, and call $E$ the {\it $F$-centralizer field of $\alpha$.} We write $\scr O_F(\alpha)$ for the $\scr W_F$-orbit of $\alpha$. 
\par
Let $\wW F = \bigcup_{n\ge1} \Go nF$. Starting with $\sigma \in \wW F$, the restriction $\sigma|_{\scr P_F}$ is a direct sum of irreducible representations $\alpha$ of $\scr P_F$, all lying in the same $\scr W_F$-orbit. Thus $\sigma$ determines an orbit $r^1_F(\sigma) = \scr O_F(\alpha) \in \scr W_F\backslash \wP F$. If $r^1_F(\sigma)$ consists of the trivial representation of $\scr P_F$, we say that $\sigma$ is {\it tamely ramified.} If it consists of characters of $\scr P_F$, we say that $\sigma$ is {\it essentially tame.} 
\par 
In the opposite direction, take $\alpha\in \wP F$ and set $E = Z_F(\alpha)$. There then exists $\rho\in \wW E$ such that the restriction of $\rho$ to $\scr P_E = \scr P_F$ is equivalent to $\alpha$. This condition determines $\rho$ up to tensoring with a tamely ramified character of $\scr W_E$. Let $\sigma \in \wW F$ contain $\alpha$ with multiplicity $m\ge1$, that is, let 
$$ 
\dim\Hom{\scr P_F}\alpha\sigma = m \ge 1. 
$$ 
There is then a tamely ramified representation $\tau\in \wW E$, of dimension $m$, such that 
$$
\sigma \cong \Ind_{\scr W_E}^{\scr W_F}\,\rho\otimes \tau. 
\tag 1 
$$
Here, $\Ind_{\scr W_E}^{\scr W_F}$ denotes the functor of {\it smooth\/} induction, for which we tend to use the briefer notation $\Ind_{E/F}$. This presentation of $\sigma$ has strong uniqueness properties. 
\par 
The expression (1) has an alternative version, often more useful in proofs. Let $E_m/E$ be unramified of degree $m$. There is then an $E_m/E$-regular, tamely ramified character $\xi$ of $\scr W_{E_m}$ such that $\tau\cong \Ind_{E_m/E}\,\xi$. Setting $\rho_m = \rho|_{\scr W_{E_m}}$, we get 
$$
\sigma \cong \Ind_{E_m/F}\,\rho_m\otimes\xi. 
\tag 2
$$ 
\indent 
The presentation (1) carries with it another structure. Let $X_1(E)$ (resp\. $X_0(E)$) denote the group of tamely ramified (resp\. unramified) characters of $E^\times$ or, via class field theory, of $\scr W_E$. If $\sigma\in \wW F$ contains $\alpha \in \wP F$, we write it in the form (1) and, for $\phi\in X_1(E)$, we define 
$$
\phi\odot_\alpha\sigma = \Ind_{E/F}\,\phi\otimes\rho\otimes\tau. 
\tag 3
$$
This gives an action of $X_1(E)$ on the set $\Gg mF{\scr O_F(\alpha)}$ of $\sigma\in \wW F$ which contain $\alpha$ with multiplicity $m$. As the notation indicates, the action does depend on $\alpha$ rather than just $\scr O_F(\alpha)$. Any other choice of $\alpha$ is of the form $\alpha^\gamma$, for some $F$-automorphism $\gamma$ of $E$, and we have the relation 
$$
\phi^\gamma\odot_{\alpha^\gamma}\sigma = \phi\odot_\alpha\sigma. 
$$ 
The case $m=1$ is of particular interest, since then $\Gg1F{\scr O_F(\alpha)}$ is a {\it principal homogeneous space\/} over $X_1(E)$. 
\subhead 
5 
\endsubhead 
Using the standard classification theory \cite{18}, analogous features are not hard to find in the irreducible cuspidal representations of groups $\GL n F$ although it takes rather longer to describe them. We need to use the theory of endo-equivalence classes (or ``endo-classes'') of simple characters and their tame lifting, as in \cite{3} and the summary in \cite{8, \S1}. We recall some of the main points. 
\par 
Let $\theta$ be a simple character in $G = \GL nF$, attached to a simple stratum $[\frak a,\beta]$ in $A = \M nF$. Let $\vT = \text{\it cl\/}(\theta)$ denote the endo-class of $\theta$. This class does not determine the field $F[\beta]$, but it does determine both the degree $[F[\beta]{:}F]$ and the ramification index $e(F[\beta]|F)$ of the field extension $F[\beta]/F$. We accordingly write $\deg\vT = [F[\beta]{:}F]$ and $e(\vT) = e(F[\beta]|F)$. 
\par 
We say that $\theta$ is {\it m-simple\/} if $\frak a$ is maximal among hereditary orders stable under conjugation by $F[\beta]^\times$. Let $\theta$ be an m-simple character in $G$, attached to a simple stratum $[\frak a,\beta]$. By definition, $\theta$ is a character of the group $H^1_\theta = H^1(\beta,\frak a)$. The $G$-normalizer $\bk J_\theta$ of $\theta$ is an open, compact modulo centre subgroup of $G$, having a unique maximal compact subgroup $J^0_\theta = J^0(\beta,\frak a)$. An {\it extended maximal simple type over $\theta$} is an irreducible representation $\vL$ of $\bk J_\theta$, containing $\theta$, and such that $\vL|_{J^0_\theta}$ is a maximal simple type in $G$, in the sense of \cite{18}. Let $\scr T(\theta)$ denote the set of equivalence classes of extended maximal simple types over $\theta$. 
\par
Taking a different standpoint, let $\pi$ be an irreducible cuspidal representation of $G$: in our earlier notation, $\pi\in \Ao nF$. The representation $\pi$ contains an m-simple character $\theta_\pi$. Any two choices of $\theta_\pi$ are $G$-conjugate, hence endo-equivalent. In other words, the endo-class $\vt(\pi) = \text{\it cl\/}(\theta_\pi)$ depends only on $\pi$. Certainly, $\deg\vt(\pi)$ divides $n$. 
\par 
Let $\vT$ be an endo-class of degree $d$ dividing $n$, say $n=md$. Let $\Aa mF\vT$ be the set of $\pi\in \Ao nF$ such that $\vt(\pi) = \vT$. There is a unique $G$-conjugacy class of m-simple characters $\theta$ in $G$ satisfying $\text{\it cl\/}(\theta) = \vT$. The principal results of \cite{18} assert that, for any such $\theta$, we have a bijection 
$$
\aligned 
\scr T(\theta) &\longrightarrow \Aa mF\vT, \\ 
\vL &\longmapsto \cind_{\bk J_\theta}^G\,\vL. 
\endaligned 
\tag 4 
$$ 
\subhead 
6 
\endsubhead  
We introduce a new element of structure. We start with an m-simple character $\theta$ in $G = \GL nF$, attached to a simple stratum $[\frak a,\beta]$. We let $T/F$ be the maximal tamely ramified sub-extension of the field extension $F[\beta]/F$. The character $\theta$ then determines $T$ uniquely, up to {\it unique\/} $F$-isomorphism. We refer to $T/F$ as a {\it tame parameter field for $\theta$.} 
\par 
We define an action of $X_1(T)$ on the set $\scr T(\theta)$. Let $G_T$ denote the $G$-centralizer of $T^\times$. Thus $G_T\cong \GL{n_T}T$, where $n/n_T = [T{:}F]$. Let $\det_T:G_T\to T^\times$ be the determinant map. Let $\phi\in X_1(T)$. There is a unique character $\phi^\bk J$ of $\bk J_\theta$ which is trivial on $J^1_\theta = J^1(\beta,\frak a)$ and agrees with $\phi\circ \det_T$ on $G_T\cap \bk J_\theta$. For $\vL\in \scr T(\theta)$, we define 
$$
\phi\odot\vL = \phi^\bk J\otimes \vL. 
$$
The representation $\phi\odot\vL$ is again an extended maximal simple type over $\theta$, and we have defined an action 
$$
\aligned 
X_1(T) \times \scr T(\theta) &\longrightarrow \scr T(\theta), \\ 
(\phi,\vL) &\longmapsto \phi\odot\vL, 
\endaligned 
\tag 5
$$
of the abelian group $X_1(T)$ on the set $\scr T(\theta)$. 
\par 
One is tempted to use the induction relation (4) to transfer the action (5) to one of $X_1(T)$ on $\Aa mF\vT$. However, the resulting action {\it depends on the choice of $\theta$ within its $G$-conjugacy class.} 
\par
We must therefore proceed more circumspectly, via an external definition of the concept of tame parameter field. An endo-class $\vF$ is called {\it totally wild\/} if $\deg\vF = e(\vF) = p^r$, for some integer $r\ge0$. Taking $\vT$ as before, a {\it tame parameter field for $\vT$} is a finite, tamely ramified field extension $E/F$ such that $\vT$ has a totally wild $E/F$-lift, the extension $E/F$ being minimal for this property. Such a field $E$ exists, and it is uniquely determined up to $F$-isomorphism. 
\par 
We choose an m-simple character $\theta$ in $\GL nF$, with tame parameter field $T/F$ and $\text{\it cl\/}(\theta) = \vT$. The field $T$ is $F$-isomorphic to $E$. The pair $(\theta,T)$ determines an endo-class $\vT_T$, which is a totally wild $T/F$-lift of $\vT$: this is the endo-class of the restriction of $\theta$ to the group $H^1_\theta\cap G_T$. If $f:E\to T$ is an $F$-isomorphism, the pull-back $f^*\vT_T$ is a totally wild $E/F$-lift of $\vT$. If we fix a totally wild $E/F$-lift $\varPsi$ of $\vT$, there is a unique choice of $f:E\to T$ such that $f^*\vT_T = \varPsi$. Combining this with (4) and (5), the resulting action of $X_1(E)$ on $\Aa mF\vT$ depends only on $\varPsi$. We accordingly denote it 
$$
\align 
X_1(E) \times \Aa mF\vT &\longrightarrow \Aa mF\vT, \\ 
(\phi,\pi) &\longmapsto \phi\odot_\varPsi \pi. 
\endalign 
$$
The case $m=1$ is of particular interest, since then $\Aa1F\vT$ becomes a {\it principal homogeneous space over $X_1(E)$.} 
\par 
The various actions are related as follows. If $\varPsi$, $\varPsi'$ are totally wild $E/F$-lifts of $\vT$, there exists an automorphism $\gamma$ of $E$ such that $\varPsi' = \varPsi^\gamma$. We then have 
$$
\phi^\gamma\odot_{\varPsi^\gamma} \pi = \phi\odot_{\varPsi} \pi, 
$$ 
for all $\phi\in X_1(E)$ and all $\pi\in \Aa mF\vT$. All of these actions extend the standard action of twisting by  characters of $F^\times$: if $\chi\in X_1(F)$ and $\chi_E = \chi\circ \N EF$, then 
$$
\chi_E\odot_\varPsi \pi = \chi\pi: g\longmapsto \chi(\det g)\,\pi(g). 
$$ 
\indent 
Overall, these discussions attach to $\pi\in \Ao nF$ an endo-class $\vT = \vt(\pi)$, a finite tame extension $E/F$, an integer $m = n/\deg\vT$ and, relative to the choice of a totally wild $E/F$-lift $\varPsi$ of $\vT$, an action $\odot_\varPsi$ of $X_1(E)$ on $\Aa mF\vT$. In particular, $\Aa1F\vT$ is a principal homogeneous space over $X_1(E)$. All of this is, at least superficially, parallel to the structures uncovered on the Galois side. 
\subhead 
7 
\endsubhead 
We may now start to build a connection between the two sides. The foundation is a result from \cite{8}, generalizing the ramification theorem of local class field theory. Let $\scr E(F)$ denote the set of endo-classes of simple characters over $F$. 
\proclaim{Ramification Theorem} 
There is a unique bijection $\Phi_F: \scr W_F\backslash \wP F \to \scr E(F)$ such that $\vt(\upr L\sigma) = \Phi_F(r^1_F(\sigma))$, for all $\sigma\in \wW F$. 
\endproclaim 
This result expresses the relationships between the Langlands Correspon\-dence, tamely ramified automorphic induction and tame lifting for endo-classes. It relies on the explicit formula \cite{16} for the {\it conductor\/} of a pair of representations. We show here that the map $\Phi_F$ is compatible with tamely ramified base field extension: if $K/F$ is a finite tame extension and $\alpha\in \wP F$, then $\Phi_K(\alpha)$ is a $K/F$-lift of $\Phi_F(\alpha)$. We deduce: 
\proclaim{Tame Parameter Theorem} 
Let $\alpha\in \wP F$ and write $\Phi_F(\alpha) = \Phi_F(\scr O_F(\alpha))$. 
\roster 
\item 
The $F$-centralizer field $Z_F(\alpha)/F$ of $\alpha$ is a tame parameter field for $\Phi_F(\alpha)$ and 
$$
\deg\Phi_F(\alpha) = [Z_F(\alpha){:}F]\,\dim\alpha. 
$$ 
\item 
For each $m\ge1$, the Langlands correspondence induces a bijection 
$$
\Gg mF{\scr O_F(\alpha)} \longrightarrow \Aa mF{\Phi_F(\alpha)}. 
\tag 6
$$
\endroster 
\endproclaim 
We therefore concentrate on analyzing the bijection (6) of the theorem. 
\subhead 
8 
\endsubhead 
We return to the presentation (2) of an element $\sigma$ of $\Gg mF{\scr O_F(\alpha)}$. In particular, $\alpha\in \wP F$ has $F$-centralizer field $E/F$. We take the unramified extension $E_m/E$, in $\bar F$, of degree $m$ and set $\vD = \Gal {E_m}E$. The orbit $\scr O_{E_m}(\alpha)$ is just $\{\alpha\}$, so we abbreviate $\Gg1{E_m}{\scr O_{E_m}(\alpha)} = \Gg1{E_m}\alpha$. Note that all representations $\nu\in \Gg1{E_m}\alpha$ are {\it totally wildly ramified.} 
\par
The group $\vD$ acts on $\Gg1{E_m}\alpha$ in a natural way; let $\Gg1{E_m}\alpha^ {\text{$\vD$-reg}}$ denote the subset of $\vD$-regular elements. As in (2) above, we have a canonical bijection 
$$
\Ind_{E_m/F}: \vD\backslash\Gg1{E_m}\alpha^{\text{$\vD$-reg}} @>{\ \ \approx\ \ }>> \Gg mF{\scr O_F(\alpha)}. 
$$ 
If we take $\phi\in X_1(E)$ and write $\phi_m = \phi\circ\N{E_m}E$, this bijection satisfies 
$$
\Ind_{E_m/F}\,\phi_m\otimes \nu = \phi\odot_\alpha \Ind_{E_m/F}\,\nu,  \quad \nu\in \Gg1{E_m}\alpha ^{\text{$\vD$-reg}}, 
$$ 
by definition. 
\par 
Moving to the other side, set $\varPsi = \Phi_E(\alpha)$ and $\varPsi_m = \Phi_{E_m}(\alpha)$. Thus $\varPsi$ is a totally wild $E/F$-lift of $\vT = \Phi_F(\alpha)$, and $E/F$ is a tame parameter field for $\vT$. The endo-class $\varPsi_m$ is the unique $E_m/E$-lift of $\varPsi$. The natural action of $\vD$ on endo-classes over $E_m$ fixes $\varPsi_m$, so $\vD$ acts on $\Aa1{E_m}{\varPsi_m}$. The Langlands correspondence induces a $\vD$-bijection $\Gg1{E_m}\alpha \to \Aa1{E_m}{\varPsi_m}$ such that 
$$
\upr L(\phi_m\otimes \nu) = \phi_m\,\upr L\nu. 
$$
\indent 
The heart of the paper is the explicit construction in \S5 of a canonical bijection 
$$ 
\ind_{E_m/F}: \vD\backslash \Aa1{E_m}{\varPsi_m}^{\text{$\vD$-reg}} @>{\ \ \approx\ \ }>> \Aa mF\vT 
$$
such that 
$$
\ind_{E_m/F}(\phi_m\rho) = \phi\odot_\varPsi \ind_{E_m/F}\,\rho, 
$$
for $\rho\in \Aa1{E_m}{\varPsi_m}^{\text{$\vD$-reg}}$ and $\phi\in X_1(E)$. This is defined in terms of extended maximal simple types, using a developed version of the Glauberman correspondence \cite{21} from the character theory of finite groups. It is elementary and completely explicit in nature, but we say no more of it in this introductory essay. 
\par 
Immediately, there is a unique bijection 
$$ 
\align 
\Bbb N: \Gg mF{\scr O_F(\alpha)} &\longrightarrow \Aa mF\vT, \\ 
\sigma &\longmapsto \upr N\sigma, 
\endalign 
$$ 
such that the following diagram commutes. 
$$
\CD 
\Gg1{E_m}\alpha^{\text{\rm $\vD$-reg}} @>{\Bbb L}>> \Aa1{E_m}{\varPsi_m}^{\text{\rm $\vD$-reg}} 
\hphantom{.} \\ 
@V{\Ind_{E_m/F}}VV @VV{\ind_{E_m/F}}V \\ 
\Gg mF{\scr O_F(\alpha)} @>>{\Bbb N}> \Aa mF\vT. 
\endCD \tag $\star$ 
$$ 
This {\it na\"\i ve correspondence\/} $\Bbb N$ further satisfies 
$$
\upr N(\phi\odot_\alpha \sigma) = \phi\odot_\varPsi \upr N\sigma, 
$$
for $\phi \in X_1(E)$, $\sigma \in \Gg mF{\scr O_F(\alpha)}$, and where $\varPsi = \Phi_E(\alpha)$. 
\subhead 
9 
\endsubhead 
Our main result compares the two bijections $\Gg mF{\scr O_F(\alpha)} \to \Aa mF\vT$ induced by, respectively, the Langlands correspondence $\sigma\mapsto \upr L\sigma$ and the na\"\i ve correspondence $\sigma\mapsto \upr N\sigma$. Let $X_0(E)_m$ denote the group of $\chi\in X_0(E)$ for which $\chi^m=1$. 
\proclaim{Comparison Theorem} 
Let $\alpha \in \wP F$ have $F$-centralizer field $E/F$ and let $m\ge1$ be an integer. There exists a character $\mu = \mu^F_{m,\alpha}\in X_1(E)$ such that 
$$
\upr L\sigma = \mu\odot_{\Phi_E(\alpha)} \upr N\sigma, 
$$ 
for all $\sigma\in \Gg mF{\scr O_F(\alpha)}$. The character $\mu^F_{m,\alpha}$ is uniquely determined modulo $X_0(E)_m$. 
\endproclaim 
We note that, if $\chi\in X_0(E)_m$ and $\pi\in \Aa mF{\Phi_F(\alpha)}$, then $\chi\odot_{\Phi_E(\alpha)}\pi = \pi$. The Comparison Theorem has the following immediate corollary, tying together the two ``interior twisting'' operations. 
\proclaim{Homogeneity Theorem} 
If $\sigma\in \Gg mF{\scr O_F(\alpha)}$ and $\phi\in X_1(E)$, then 
$$
\upr L(\phi\odot_\alpha\sigma) = \phi\odot_{\Phi_E(\alpha)} \upr L\sigma. 
$$
\endproclaim 
The third of our main results is the Types Theorem of 7.6. If $\sigma\in \wW F$ and $\pi = \upr L\sigma$, then $\pi$ contains a maximal simple type. The Types Theorem allows one to read off the structure of this type from the presentation (1) of $\sigma$. 
\subhead 
10 
\endsubhead 
In the essentially tame case, we gave a complete account of the ``discrepancy character'' $\mu^F_{m,\alpha}$, purely in terms of the simple stratum underlying the endo-class $\Phi_F(\alpha)$ \cite{10}, \cite{11} \cite{13}. Even there, it poses challenging problems and we make no such attempt here. In the general case, we determine its restriction to units of $E$: this is the essence of the Types Theorem. Consideration of central characters yields its restriction to $(F^\times)^{m\dim\alpha}$ (7.3.2). The character $\mu^F_{m,\alpha}$ of $E^\times$ is thereby determined up to an unramified factor of degree dividing $me(E|F)\dim\alpha$ (remembering that $\mu^F_{m,\alpha}$ is only defined modulo unramified characters of order dividing $m$). In principle, therefore, it is amenable to description in terms of local constants of a finite number of pairs, using \cite{7}. 
\subhead 
11 
\endsubhead 
We indicate briefly the strategy of the proof. This follows the essentially tame case quite closely, but the technical hurdles are somewhat higher. 
\par 
The diagram $(\star)$ has an analogue relative to the Langlands correspondence: there is a unique bijection 
$$ 
\roman a_{E_m/F}: \vD\backslash \Aa 1{E_m}{\varPsi_m}^{\text{$\vD$-reg}}@>{\ \ \approx\ \ }>> \Aa mF\vT 
$$ 
such that the diagram 
$$
\CD 
\Gg1{E_m}\alpha^{\text{\rm $\vD$-reg}} @>{\Bbb L}>> \Aa1{E_m}{\varPsi_m}^{\text{\rm $\vD$-reg}} 
\\ 
@V{\Ind_{E_m/F}}VV @VV{\roman a_{E_m/F}}V \\ 
\Gg mF{\scr O_F(\alpha)} @>>{\Bbb L}> \Aa mF\vT. 
\endCD \tag $\dagger$ 
$$ 
commutes. We have to compare $\roman a_{E_m/F}$ with $\ind_{E_m/F}$. 
\par 
The field extension $E_m/F$ admits a tower of subfields 
$$
E_m = L_0\supset L_1\supset \dots \supset L_r \supset L_{r+1} = F, 
$$ 
in which 
\roster 
\item $L_r/F$ is unramified; 
\item $L_i/L_{i+1}$ is cyclic and totally ramified of prime degree, $1\le i\le r{-}1$; 
\item the group of $L_1$-automorphisms of $L_0$ is trivial. 
\endroster 
The map $\roman a_{E_m/F}$ has a corresponding factorization 
$$
\roman a_{E_m/F} = \roman A_{L_r/F}\circ \roman A_{L_{r-1}/L_r} \circ \dots \roman A_{L_1/L_2} \circ \roman a_{L_0/L_1}, 
$$
in which each map labelled $\roman A$ is {\it automorphic induction,} in the sense of \cite{28}, \cite{29}, \cite{30}. The first factor, $\roman a_{L_0/L_1}$, is defined by the diagram $(\dagger)$ relative to the base field $L_1$ in place of $F$. 
\par 
Comparison of $\roman a_{L_0/L_1}$ with $\ind_{L_0/L_1}$ is straightforward. One sees easily that these two maps differ by the $\odot$-twist with a character $\nu_{L_0/L_1}\circ \N{E_m}E$, independent of the representations under consideration. However, the method gives no information about this character. From then on, we work inductively, comparing $\ind_{L_0/L_{i+1}}$ with $\roman A_{L_i/L_{i+1}} \circ \ind_{L_0/L_i}$. We have to show that these two maps differ by the $\odot$-twist with a character $\nu_{L_i/L_{i+1}}\in X_1(E)$. The discrepancy character $\mu$ is then the product of these characters $\nu$. 
\par 
We give an exact, rather than informative, expression for each of the characters $\nu_{L_i/L_{i+1}}$, $1\le i\le r$. When the sub-extension in question is totally ramified, that is, when $1\le i\le r{-}1$, the formula (8.9 Corollary) involves an induction constant and transfer factor arising from the automorphic induction equation, along with some symplectic signs coming from the Glauberman correspondence. While we make no attempt to evaluate it in general, a simple exercise shows, for instance, that the character is trivial when $p[L_i{:}L_{i+1}]$ is odd. For the unramified extension $L_r/L_{r+1} = F$, the corresponding character has order at most $2$. We give an explicit formula (10.7 Corollary) involving only elementary quantities and some symplectic signs. 
\par 
Our handling of the automorphic induction equation, for a tame cyclic extension $K/F$, directly generalizes the method of \cite{10}, \cite{13}. It requires the Uniform Induction Theorem of \cite{13}, \cite{29} to control an induction constant. Since we can only rely on rather weaker linear independence properties of characters, we have to compensate by a closer investigation of transfer factors. These do not seem amenable to direct computation, so we have had to evolve some novel methods. Likewise, we have had to develop a more structural approach to the symplectic signs. 
\subhead 
Note on characteristic 
\endsubhead  
We shall often refer to our earlier papers \cite{3}, \cite{4} and their successors, in which we imposed the hypothesis that $F$ be of characteristic zero. The only reason for doing that was the lack of a theory of base change and automorphic induction in positive characteristic. Such material is now available in \cite{30}, so the results of all of our earlier papers in the area apply equally in positive characteristic. 
\subhead 
Notation and conventions 
\endsubhead 
The following will be standard throughout, and used without further elaboration. If $F$ is a non-Archimedean local field, we denote by $\frak o_F$ the discrete valuation ring in $F$, by $\frak p_F$ the maximal ideal of $\frak o_F$, and by $\Bbbk_F$ the residue class field $\frak o_F/\frak p_F$. The characteristic of $\Bbbk_F$ is always denoted $p$. We set $U_F = U^0_F = \frak o_F^\times$ and $U^k_F = 1{+}\frak p_F^k$, $k\ge1$. We write $\bs\mu_F$ for the group of roots of unity in $F$ of order prime to $p$. In particular, $U_F = \bs\mu_F\times U^1_F$ and reduction modulo $\frak p_F$ induces an isomorphism $\bs\mu_F \cong \Bbbk_F^\times$. 
\par 
If $E/F$ is a finite field extension, we denote by $\N EF:E^\times \to F^\times$ and $\roman{Tr}_{E/F}:E\to F$ the norm and trace maps respectively. We write $\roman{Aut}(E|F)$ for the group of $F$-automorphisms of the field $E$. If $E/F$ is Galois, then $\roman{Aut}(E|F) = \Gal EF$. 
\par 
We fix, once for all, a separable algebraic closure $\bar F/F$ and let $\scr W_F$ be the Weil group of $\bar F/F$. We let $\scr I_F$, $\scr P_F$ denote respectively the inertia subgroup and the wild inertia subgroup of $\scr W_F$. If $E/F$ is a finite separable extension with $E\subset \bar F$, we identify the Weil group $\scr W_E$ of $\bar F/E$ with the group of elements of $\scr W_F$ which fix $E$. 
\par 
We denote by $X_0(F)$ (resp\. $X_1(F)$) the group of smooth characters $\scr W_F\to \Bbb C^\times$ which are unramified (resp\. tamely ramified) in the sense of being trivial on $\scr I_F$ (resp\. $\scr P_F$). 
\par 
We let $\bk a_F:\scr W_F \to F^\times$ be the Artin Reciprocity map, normalized to take {\it geometric\/} Frobenius elements of $\scr W_F$ to prime elements of $F$. We use $\bk a_F$ to identify the group of smooth characters of $\scr W_F$ with the group of smooth characters of $F^\times$. In particular, we identify $X_0(F)$ (resp\. $X_1(F)$) with the group of characters of $F^\times$ which are trivial on $U_F$ (resp\. $U^1_F$). If $m\ge1$ is an integer, then $X_0(F)_m$ is the group of $\chi\in X_0(F)$ such that $\chi^m=1$. 
\par 
Let $n\ge1$ be an integer; let $A = \M nF$ (the algebra of $n\times n$ matrices over $F$) and $G = \GL nF$. If $\frak a$ is a hereditary $\frak o_F$-order in $A$, with Jacobson radical $\roman{rad}\,\frak a = \frak p_\frak a$, we write $U_\frak a = U^0_\frak a = \frak a^\times$, and $U^k_\frak a = 1{+}\frak p_\frak a^k$, $k\ge1$. We write $\scr K_\frak a = \{x\in G:x^{-1}\frak ax = \frak a\}$. The group $\scr K_\frak a$ is also the $G$-normalizer of $U_\frak a$. 
\par
If $\vD$ is a  group acting on a set $X$, we denote by $X^\vD$ the set of $\vD$-fixed points in $X$. An element $x$ of $X$ is {\it $\vD$-regular\/} if its $\vD$-isotropy is trivial. We denote by $X^{\text{$\vD$-reg}}$ the set of $\vD$-regular elements of $X$. 
\head\Rm 
1. Representations of Weil groups 
\endhead 
We give a systematic account of the irreducible smooth representations of the locally profinite group $\scr W_F$, in terms of restriction to the wild inertia subgroup $\scr P_F$ of $\scr W_F$. The material of the section, as far as the end of 1.4, is basically well-known but we need to have all of the relevant details in one place. 
\par
Let $E/F$ be a finite field extension, with $E\subset \bar F$. If $\sigma$, resp\. $\rho$, is an irreducible smooth representation of $\scr W_F$, resp\. $\scr W_E$, then both the smoothly induced representation $\Ind_{E/F}\,\rho = \Ind_{\scr W_E}^{\scr W_F}\,\rho$, and the restricted representation $\sigma|_{\scr W_E}$, are finite-dimensional and {\it semisimple, cf\.} \cite{12 (28.7 Lemma)}. Consequently, we may work entirely within the category of finite-dimensional, smooth, semisimple representations of $\scr W_F$, using the elementary  methods of Clifford-Mackey theory. 
\subhead 
1.1 
\endsubhead 
Let $n\ge1$ be an integer. Let $\Go nF$ denote the set of equivalence classes of irreducible smooth representations of $\scr W_F$ of dimension $n$. We set 
$$ 
\wW F = \bigcup_{n\ge1} \Go nF. 
$$ 
If $\sigma$ is a smooth representation of $\scr W_F$, we say that $\sigma$ is {\it unramified\/} (resp\. {\it tamely ramified\/}) if $\roman{Ker}\,\sigma$ contains $\scr I_F$ (resp\. $\scr P_F$). 
\subhead 
1.2 
\endsubhead 
Let $\wP F$ denote the set of equivalence classes of irreducible smooth representations of the profinite group $\scr P_F$. The group $\scr W_F$ acts on $\wP F$ by conjugation. 
\par 
Let $\alpha\in \widehat{\scr P}_F$. Since $\scr P_F$ is a pro-$p$ group, the group $\alpha(\scr P_F)$ is finite. The dimension of $\alpha$ is finite and of the form $p^r$, for an integer $r\ge0$. Let 
$$ 
N_F(\alpha) = \{g\in \scr W_F: \alpha^g\cong\alpha\}, 
$$ 
where $\alpha^g$ denotes the representation $x\mapsto \alpha(gxg^{-1})$ of $\scr P_F$. 
\proclaim{Proposition} 
Let $\alpha\in \wP F$. There is a finite, tamely ramified field extension $E/F$, with $E\subset \bar F$, such that $N_F(\alpha) = \scr W_E$. 
\endproclaim 
\demo{Proof} 
Since $\dim\alpha$ is finite, the smooth representation of $\scr W_F$ {\it compactly\/} induced by $\alpha$ is finitely generated over $\scr W_F$. It therefore admits an irreducible quotient, $\sigma$ say. The dimension of $\sigma$ is finite (see 28.6 Lemma 1 of \cite{12}), so $\sigma|_{\scr P_F}$ is a finite direct sum of irreducible representations, all conjugate to $\alpha$. Moreover, any $\scr W_F$-conjugate of $\alpha$ occurs here. Thus $N_F(\alpha)$ is closed and of finite index in $\scr W_F$.  It follows that $N_F(\alpha)$ is open in $\scr W_F$, of finite index and containing $\scr P_F$. That is, $N_F(\alpha) = \scr W_E$, for a finite, tamely ramified extension $E/F$ inside $\bar F$. \qed 
\enddemo 
In the situation of the proposition, we call $E$ the {\it $F$-centralizer field\/} of $\alpha$, and denote it $Z_F(\alpha)$.  
\subhead 
1.3 
\endsubhead 
We consider extension properties of a given representation $\alpha\in \wP F$. 
\par
Let $E = Z_F(\alpha)$. If $K/F$ is a finite, tamely ramified extension, inside $\bar F$, such that $\alpha$ admits extension to a representation of $\scr W_K$, then surely $K\supset E$. In the opposite direction, we prove: 
\proclaim{Proposition} 
Let $\alpha\in \wP F$, and let $E = Z_F(\alpha)$. 
\roster 
\item 
There exists an irreducible smooth representation $\rho$ of\/ $\scr W_E$ such that 
\itemitem{\rm (a)} $\rho|_{\scr P_F} \cong\alpha$, and 
\itemitem{\rm (b)} $\det\rho|_{\scr I_F}$ has $p$-power order. 
\item 
If $\rho'$ is a smooth representation of\/ $\scr W_E$ such that $\rho'|_{\scr P_F} \cong \alpha$, then there exists a unique character $\chi \in X_1(E)$ such that $\rho'\cong \rho\otimes\chi$. 
\endroster 
\endproclaim 
\demo{Proof} 
We first concentrate on the uniqueness properties. 
\proclaim{Lemma 1} 
\roster 
\item 
Let $\theta$ be an irreducible smooth representation of\/ $\scr I_F$ such that $\theta|_{\scr P_F}$ is irreducible. If\/ $\psi$ is a character of\/ $\scr I_F$, trivial on $\scr P_F$, such that $\theta\otimes\psi \cong \theta$, then $\psi = 1$. 
\item 
Let $\tau$ be an irreducible smooth representation of\/ $\scr W_F$ such that $\tau|_{\scr P_F}$ is irreducible. If\/ $\psi$ is a tamely ramified character of\/ $\scr W_F$, such that $\tau\otimes\psi \cong \tau$, then $\psi = 1$. 
\endroster 
\endproclaim 
\demo{Proof} 
In part (1), we view $\theta$ and $\psi\otimes\theta$ as acting on the same vector space $V$, say. By hypothesis, there is a map $f\in \Aut{\Bbb C}V$ such that $f\circ\theta(g) = \psi(g)\theta(g)\circ f$, for all $g\in \scr I_F$. In particular, $f\circ\theta(x) = \theta(x)\circ f$, for $x\in \scr P_F$. Since $\theta|_{\scr P_F}$ is irreducible, the map $f$ is a non-zero scalar, whence $\psi$ is trivial, as required. The proof of (2) is identical. \qed 
\enddemo 
\proclaim{Lemma 2} 
\roster 
\item 
There exists a unique smooth representation $\rho_\alpha$ of\/ $\scr I_E$ such that 
\itemitem{\rm (a)} $\rho_\alpha|_{\scr P_F} \cong \alpha$, and 
\itemitem{\rm (b)} the character $\det\rho_\alpha$ has finite $p$-power order. 
\item  
If $\rho'$ is an irreducible smooth representation of\/ $\scr I_E$ such that $\rho'|_{\scr P_F}$ con\-tains $\alpha$, then $\rho' \cong \rho_\alpha\otimes \psi$, for a unique character $\psi$ of\/ $\scr I_E/\scr P_F$. 
\endroster 
\endproclaim 
\demo{Proof} 
Let $\tau$ be an irreducible smooth representation of $\scr I_E$ such that $\tau|_{\scr P_F}$ contains $\alpha$. Thus $\tau|_{\scr P_F}$ is a multiple of $\alpha$. If $\scr K = \roman{Ker}\,\tau$, then $\alpha$ surely extends uniquely to a representation $\tau_1$ of $\scr P_F\scr K$ which is trivial on $\scr K$. The restriction $\tau|_{\scr P_F\scr K}$ is thus a multiple of $\tau_1$. The quotient $\scr I_E/\scr P_F\scr K$ is finite cyclic, and $\tau_1$ is stable under conjugation by $\scr I_E$. Therefore $\tau_1$ admits extension to a representation, $\tilde\tau_1$ say, of $\scr I_E$. The irreducible representations of $\scr I_E$ containing $\tau_1$ are then of the form $\chi\otimes\tilde\tau_1$, where $\chi$ ranges over the characters of $\scr I_E$ trivial on $\scr P_F\scr K$. In particular, $\tau|_{\scr P_F\scr K} \cong \tau_1$ and so $\tau|_{\scr P_F}\cong \alpha$. 
\par 
This argument also shows that, if $\tau'$ is an irreducible smooth representation of $\scr I_E$ containing $\alpha$, then $\tau'\cong \tau\otimes\psi$, for a character $\psi$ of $\scr I_E/\scr P_F$. Consider the determinant character 
$$ 
\det(\tau\otimes\psi) = \psi^{\dim\tau}\det\tau. 
\tag 1.3.1 
$$ 
The dimension $\dim\tau = \dim\alpha$ is a power of $p$, the character $\det\tau$ has finite order, while $\psi$ has finite order not divisible by $p$. Thus there is a unique choice of $\psi$ such that $\det(\tau\otimes\psi)$ has $p$-power order. This proves (1). 
\par 
In (2) the existence assertion has already been proved. The uniqueness property is given by Lemma 1. \qed
\enddemo 
We prove the proposition. The uniqueness property of $\rho_\alpha$ (as in Lemma 2) implies that $\rho_\alpha$ is stable under conjugation by $\scr W_E$. The group $\scr W_E/\scr I_E$ is infinite cyclic and discrete, so $\rho_\alpha$ admits extension to a smooth representation $\rho$ of $\scr W_E$, as required for part (1). For the same reason, any irreducible representation of $\scr W_E$ containing $\rho_\alpha$ is of the form $\psi\otimes\rho$, for a character $\psi$ of $\scr W_E/\scr P_F$, as required for part (2). The uniqueness assertion of (2) is given by Lemma 1. \qed 
\enddemo
\subhead 
1.4 
\endsubhead 
To analyze the irreducible smooth representations of $\scr W_F$, we introduce a class of objects generalizing the admissible pairs of \cite{10}. (See below, 1.6 Remark, for the precise connection.) 
\proclaim{Definition} 
An \rom{admissible datum} over $F$ is a triple $(E/F,\rho,\tau)$ satisfying the following conditions. 
\roster 
\item 
$E/F$ is a finite, tamely ramified field extension with $E\subset \bar F$. 
\item 
$\rho$ is an irreducible smooth representation of\/ $\scr W_E$ such that the restriction $\alpha = \rho|_{\scr P_F}$ is irreducible and $Z_F(\alpha) = E$. 
\item 
$\tau$ is an irreducible, smooth, tamely ramified representation of\/ $\scr W_E$. 
\endroster 
\endproclaim 
Two admissible data $(E_i/F,\rho_i,\tau_i)$ over $F$, $i=1,2$, are deemed equivalent if there exist $g\in \scr W_F$ and $\chi\in X_1(E_2)$ such that $E_2 = E_1^g$, $\rho_2 \cong \rho_1^g\otimes\chi$ and $\tau_2 \cong \tau_1^g\otimes\chi^{-1}$. Let $\frak G(F)$ denote the set of equivalence classes of admissible data over $F$. 
\par 
If $(E/F,\rho,\tau)$ is an admissible datum over $F$, we define 
$$
\vS(\rho,\tau) = \Ind_{E/F}\,\rho\otimes\tau. 
\tag 1.4.1 
$$
Surely $\vS(\rho,\tau)$ depends only on the equivalence class of the datum $(E/F,\rho,\tau)$. The main result of the section is the following. 
\proclaim{Theorem} 
If $(E/F,\rho,\tau)$ is an admissible datum over $F$, the representation $\vS(\rho,\tau)$ of\/ $\scr W_F$ is irreducible. The map $\vS:\frak G(F) \to \wW F$ is a bijection. 
\endproclaim 
\demo{Proof} 
Let $(E/F,\rho,\tau)$ be an admissible datum. We prove that $\Ind_{E/F}\,\rho\otimes\tau$ is irreducible. Set $\alpha = \rho|_{\scr P_F}$. The tamely ramified representation $\tau\in \wW E$ is of the form $\tau = \Ind_{K/E}\,\chi$, where $K/E$ is a finite, unramified extension and $\chi\in X_1(K)$ is such that the conjugates $\chi^g$, $g\in \scr W_K\backslash \scr W_E$, are distinct. Writing $\rho_K = \rho|_{\scr W_K}$, we have $\rho\otimes\tau \cong \Ind_{K/E}\,\rho_K\otimes\chi$. Let $x\in \scr W_F$ intertwine $\rho_K\otimes\chi$. It also intertwines $\alpha = (\rho_K\otimes\chi)|_{\scr P_F}$, and therefore lies in $\scr W_E$. The conjugate representation $(\rho_K\otimes \chi)^x = \rho_K\otimes \chi^x$ is equivalent to $\rho_K\otimes\chi$. Part (2) of 1.3 Lemma 1 implies $\chi^x = \chi$, whence $x\in \scr W_K$. Thus 
$$
\Ind_{E/F}\,\rho\otimes\tau \cong \Ind_{K/F}\,\rho_K\otimes\chi 
$$
is irreducible, as required for the first assertion of the theorem. This same argument also proves: 
\proclaim{Lemma} 
Let $\tau$, $\tau'$ be irreducible, tamely ramified representations of\/ $\scr W_E$ of the same dimension. If the representations $\rho\otimes\tau$, $\rho\otimes\tau'$ are intertwined by an element of\/ $\scr W_F$, then $\tau\cong\tau'$. 
\endproclaim  
We next show that the map $\vS:\frak G(F) \to \wW F$ is injective. Let $(E/F,\rho,\tau)$, $(E'/F,\rho',\tau')$ be admissible data with the same image. The representations $\alpha =\rho|_{\scr P_F}$, $\alpha' = \rho'|_{\scr P_F}$ are therefore $\scr W_F$-conjugate. Replacing $(E'/F,\rho',\tau')$ by a conjugate, which does not affect its equivalence class, we may assume $E' = E$ and $\alpha' =\alpha$. By 1.3 Proposition, $\rho' \cong \rho\otimes\psi$, for some tamely ramified character $\psi$ of $\scr W_E$. Adjusting $(E/F,\rho',\tau')$ in its equivalence class, we can take $\rho' = \rho$. The lemma now implies $\tau'\cong \tau$, as required. 
\par
It remains to show that the map $\vS$ is surjective. Let $\sigma\in \wW F$, let $\alpha$ be an irreducible component of $\sigma|_{\scr P_F}$, and set $E = Z_F(\alpha)$. Let $\theta$ denote the natural representation of $\scr W_E$ on the $\alpha$-isotypic subspace of $\sigma$. Let $\theta_0$ be an irreducible component of $\theta$. Any element of $\scr W_F$ which intertwines $\theta_0$ also intertwines $\alpha$ and so lies in $\scr W_E$. This shows that $\Ind_{E/F}\,\theta_0$ is irreducible, hence equivalent to $\sigma$. The Mackey induction formula further shows that $\theta_0 = \theta$. 
\par
By its definition, the restriction of $\theta$ to $\scr P_F$ is a direct sum of copies of $\alpha$.  By part (2) of 1.3 Lemma 2, the restriction of $\theta$ to $\scr I_E$ is a direct sum of representations $\rho_\alpha\otimes\psi$, with $\rho_\alpha$ as in that lemma and various characters $\psi$ of $\scr I_E/\scr P_F$. Moreover, if we choose a character $\psi$ appearing here, an element $x\in \scr W_E$ intertwines $\rho\otimes\psi$ if and only if $\psi^x \cong\psi$. The $\scr W_E$-stabilizer of $\psi$ is of the form $\scr W_K$, where $K/E$ is finite and unramified. Since $\rho_\alpha$ admits extension to a representation $\rho$ of $\scr W_E$ (1.3 Proposition) and $\scr W_K/\scr I_E$ is cyclic, the representation $\rho_\alpha\otimes\psi$ extends to a representation $\rho|_{\scr W_K}\otimes \varPsi$ of $\scr W_K$, occurring in $\theta|_{\scr W_K}$. We then have $\theta \cong \rho\otimes\tau$, where $\tau = \Ind_{K/E}\,\varPsi$. The triple $(E/F,\rho,\tau)$ is an admissible datum and $\vS(\rho,\tau) \cong \sigma$, as required. 
\par
This completes the proof of the theorem. \qed 
\enddemo 
\subhead 
1.5 
\endsubhead 
Let $\scr O\in \scr W_F\backslash \wP F$, say $\scr O$ is the $\scr W_F$-orbit $\scr O_F(\alpha)$ of $\alpha\in \wP F$. We set 
$$ 
d(\scr O) = d_F(\alpha) = [Z_F(\alpha){:}F]\,\dim\alpha. 
\tag 1.5.1 
$$ 
If $m\ge1$ is an integer, we define $\Gg mF{\scr O}$ to be the set of $\sigma\in \Go {md(\scr O)}F$ such that $\sigma|_{\scr P_F}$ contains $\alpha$. From 1.4 Theorem, we deduce: 
\proclaim{Corollary} 
Let $m$ be a positive integer, let $\alpha \in \wP F$ and write $\scr O = \scr O_F(\alpha)$.  For $\sigma\in \wW F$, the following are equivalent: 
\roster 
\item $\sigma\in \Gg mF{\scr O}$; 
\item there is an admissible datum $(E/F,\rho,\tau)$ such that $\rho|_{\scr P_F} \cong \alpha$, $\dim\tau = m$ and $\sigma\cong \vS(\rho,\tau)$; 
\item $\dim\Hom{\scr P_F}\alpha\sigma = m$. 
\endroster 
\endproclaim 
The theorem also reveals a useful family of structures on the set $\wW F$. Taking $\alpha \in \wP F$ and $\scr O = \scr O_F(\alpha)$ as before, set $E = Z_F(\alpha)$. For $\sigma\in \Gg mF{\scr O}$ and $\phi \in X_1(E)$, we define a representation $\phi\odot_\alpha\sigma\in \Gg mF{\scr O}$ as follows. We write $\sigma = \vS(\rho,\tau)$, for an admissible datum $(E/F,\rho,\tau)$ with $\rho|_{\scr P_F} \cong \alpha$. The triple $(E/F,\phi\otimes\rho,\tau)$ is again an admissible datum. We put 
$$
\phi\odot_\alpha\sigma = \vS(\phi\otimes\rho,\tau) = \vS(\rho,\phi\otimes\tau). 
\tag 1.5.2 
$$
This does not depend on the choice of datum $(E/F,\rho,\tau)$ satisfying the stated conditions. The pairing $(\phi,\sigma)\mapsto \phi\odot_\alpha\sigma$ endows $\Gg mF{\scr O}$ with the structure of $X_1(E)$-space, in that 
$$
\phi\phi'\odot_\alpha\sigma = \phi \odot_\alpha (\phi'\odot_\alpha\sigma), 
$$
for $\phi,\phi'\in X_1(E)$, $\sigma \in \Gg mF{\scr O}$. 
\par  This structure {\it does\/} depend on $\alpha$ rather than on $\scr O$ and $E$. For, if $\alpha'\in \scr O$ and $Z_F(\alpha') = E$, then $\alpha' = \alpha^\gamma$, for some $\gamma\in \roman{Aut}(E|F)$. We then get the relation 
$$ 
\phi^\gamma\odot_{\alpha^\gamma} \sigma \cong \phi\odot_\alpha\sigma, 
\tag 1.5.3 
$$ 
for $\phi\in X_1(E)$, $\sigma\in \Gg mF{\scr O}$. (Clearly, the relation (1.5.3) holds equally for any $\gamma\in \scr W_F$.) 
\par 
Working for the moment with base field $E$, the $\scr W_E$-orbit of $\alpha$ is $\{\alpha\}$ and so we write $\Gg mE\alpha$ rather than $\Gg mE{\{\alpha\}}$. In this case, the definition (1.5.2) reduces to  
$$
\phi\odot_\alpha \sigma = \phi\otimes \sigma, 
\tag 1.5.4
$$
for $\sigma\in \Gg mE\alpha$ and $\phi\in X_1(E)$. 
\proclaim{Proposition} 
The map 
$$ 
\align  
\Gg mE\alpha &\longrightarrow \Gg mF{\scr O} \\ 
\sigma&\longmapsto \Ind_{E/F}\,\sigma,  
\endalign 
$$ 
is a bijection satisfying 
$$ 
\Ind_{E/F}\,\phi\otimes\sigma = \phi\odot_\alpha\Ind_{E/F}\,\sigma, \quad \phi\in X_1(E),\ \sigma \in \Gg mE \alpha. 
\tag 1.5.5 
$$ 
In particular, $\Gg1F{\scr O}$ is a principal homogeneous space over $X_1(E)$. 
\endproclaim 
\demo{Proof} 
The first assertion follows from 1.4 Theorem and the definitions. The set $\Gg 1E\alpha$ is a principal homogeneous space over $X_1(E)$, by 1.3 Proposition. The second assertion now follows from the first. \qed 
\enddemo 
\subhead 
1.6 
\endsubhead 
For some purposes, a variant of 1.4 Theorem is technically convenient. 
\par 
Let $\alpha\in \wP F$ and set $E = Z_F(\alpha)$. Let $E_m/E$ be unramified of degree $m$ and set $\vD= \Gal{E_m}E$. The group $\vD$ then acts on the set $\Gg1{E_m}\alpha$. Let $\Gg1{E_m}\alpha^{\vD \text{\rm -reg}}$ denote the subset of $\vD$-regular elements. (As usual, an element of a $\vD$-set is called $\vD$-regular if its $\vD$-isotropy group is trivial.) 
\par
If $\phi\in X_1(E)$, write $\phi_m = \phi\circ\N{E_m}E$. The group $X_1(E)$ thus acts on $\Gg1{E_m}\alpha$ by $(\phi,\zeta)\mapsto \phi_m\otimes\zeta$. The subset $\Gg1{E_m}\alpha^{\vD \text{\rm -reg}}$ is stable under this action. 
\proclaim{Proposition} 
If $\zeta\in \Gg1{E_m}\alpha^{\vD \text{\rm -reg}}$, the representation $\Ind_{E_m/F}\,\zeta$ is irreducible and lies in $\Gg mF{\scr O}$, $\scr O = \scr O_F(\alpha)$. Moreover, 
$$
\Ind_{E_m/F}\,\phi_m\otimes\zeta = \phi\odot_\alpha\Ind_{E_m/F}\,\zeta, 
$$ 
for $\phi\in X_1(E)$ and $\zeta\in \Gg1{E_m}\alpha^{\vD \text{\rm -reg}}$. The induced map 
$$
\Ind_{E_m/F}: \vD\backslash \Gg1{E_m}\alpha^{\vD \text{\rm -reg}}  \longrightarrow \Gg mF{\scr O} 
$$
is bijective. 
\endproclaim 
\demo{Proof} 
Choose $\rho\in \Gg1E\alpha$ and set $\rho_m = \rho|_{\scr W_{E_m}}$. 
\proclaim{Lemma} 
The map 
$$
\align 
X_1(E_m) &\longrightarrow \Gg1{E_m}\alpha, \\ \xi &\longmapsto \xi\otimes\rho_m, 
\endalign 
$$ 
is a bijective $\vD$-map. 
\endproclaim 
\demo{Proof} 
The lemma follows directly from 1.3 Proposition. \qed 
\enddemo 
In particular, $\xi\otimes\rho_m$ is $\vD$-regular if and only if $\xi$ is $\vD$-regular, and this condition is equivalent to $\Ind_{E_m/E}\,\xi$ being irreducible. We have 
$$
\Ind_{E_m/E}\, \xi\otimes\rho_m \cong \rho\otimes \Ind_{E_m/E}\,\xi. 
$$ 
Thus $\Ind_{E_m/E}$ induces a bijection 
$$
\vD\backslash \Gg1{E_m}\alpha^{\vD \text{\rm -reg}} @>{\ \ \approx\ \ }>> \Gg mE{\alpha}. 
$$ 
The proposition now follows from 1.5 Proposition. \qed 
\enddemo 
\remark{Remark} 
Suppose $\dim\alpha = 1$ and let $\zeta\in \Gg1{E_m}\alpha$. Thus $\dim\zeta = 1$ and we may view $\zeta$ as a character of $E_m^\times$, via local class field theory. The representation $\zeta$ lies in $\Gg1{E_m}\alpha^{\vD \text{\rm -reg}}$ if and only if the pair $(E_m/F,\zeta)$ is {\it admissible,} in the sense of \cite{10}. 
\endremark 
\head\Rm 
2. Simple characters and tame parameters
\endhead 
We recall the theory of simple characters \cite{18} and tame lifting \cite{3}. We orient the treatment so as to give a degree of prominence to a certain tamely ramified field extension attached to a simple character $\theta$. The invariance properties of this {\it tame parameter field\/} provide the only novelty here. The main result (2.7) is rather technical, but essential for what follows. 
\par 
We follow the notational schemes of \cite{18} and \cite{3}, making some minor modifications to indicate more clearly the interdependences between various objects. The notation introduced here, particularly that in 2.3, will be standard for the rest of the paper. 
\subhead 
2.1 
\endsubhead 
Let $n\ge1$ be an integer, and set $A = \M nF$, $G = \GL nF$. We consider a simple stratum in $A$ of the form $[\frak a,l,0,\beta]$ (as in \cite{18, Chapter 1}). The positive integer $l$ is determined by the relation $\beta\frak a = \frak p_\frak a^{-l}$, where $\frak p_\frak a$ is the Jacobson radical of $\frak a$, so we tend to omit it from the notation and abbreviate $[\frak a,l,0,\beta] = [\frak a,\beta]$. 
\par 
As in \cite{18, Chapter 3}, the simple stratum $[\frak a,\beta]$ determines compact open subgroups 
$$
H^1(\beta,\frak a)\subset J^1(\beta,\frak a) \subset J^0(\beta,\frak a) 
$$
of $U_\frak a$, such that $J^1(\beta,\frak a) = J^0(\beta,\frak a)\cap U^1_\frak a$. Let $A_\beta$ denote the centralizer of $\beta$ in $A$. Thus $A_\beta$ is isomorphic to a full matrix algebra over the field $F[\beta]$ and $\frak a_\beta = \frak a\cap A_\beta$ is a hereditary $\frak o_{F[\beta]}$-order in $A_\beta$. We set 
$$ 
\bk J(\beta,\frak a) = \scr K_{\frak a_\beta}\cdot  J^1(\beta,\frak a). 
$$ 
\indent 
A smooth character $\psi_F$ of $F$, of level one, attaches to the stratum $[\frak a,\beta]$ a non-empty set $\scr C(\frak a,\beta,\psi_F)$ of distinguished characters of the group $H^1(\beta,\frak a)$, as in \cite{18 (3.2)}. These are known as {\it simple characters.} Since $\psi_F$ will be fixed throughout, we usually omit it from the notation. We assemble a variety of facts about this situation, using \cite{18 (3.3.2, 3.5.1)}. 
\proclaim\nofrills{(2.1.1)} \ 
Let $[\frak a,\beta]$ be a simple stratum in $A$ and let $\theta\in \scr C(\frak a,\beta,\psi_F)$. Let $I_G(\theta)$ denote the set of elements of $G$ which intertwine $\theta$ and let $\bk J_\theta$ be the group of $g\in G$ which normalize $\theta$. 
\roster 
\item 
The group $\bk J_\theta$ is open and compact modulo centre in $G$. Moreover, 
\itemitem{\rm (a)} the group $\bk J_\theta$ has a unique maximal compact subgroup $J^0_\theta$, 
\itemitem{\rm (b)} the group $J^0_\theta$ has a unique maximal, normal pro $p$-subgroup $J^1_\theta$, and 
\itemitem{\rm (c)} the set $I_G(\theta)$ is $J^1_\theta G_\beta J^1_\theta = J^0_\theta G_\beta J^0_\theta$, where $G_\beta$ denotes the $G$ centralizer of $\beta$. 
\item 
If $[\frak a',\beta']$ is a simple stratum in $A$ such that $\theta\in \scr C(\frak a',\beta',\psi_F)$, then 
\itemitem{\rm (a)} $\frak a' = \frak a$, 
\itemitem{\rm (b)} $e(F[\beta']|F) = e(F[\beta]|F)$, $f(F[\beta']|F) = f(F[\beta]|F)$, 
\itemitem{\rm (c)} $J^1_\theta = J^1(\beta',\frak a')$, $J^0_\theta = J^0(\beta',\frak a')$, and 
\itemitem{\rm (d)} $\bk J_\theta = \bk J(\beta',\frak a')$. 
\endroster 
\endproclaim 
In particular, the conclusions (c), (d) of (2) hold for $[\frak a',\beta'] = [\frak a,\beta]$. 
\par 
If $\theta\in \scr C(\frak a,\beta)$, for a simple stratum $[\frak a,\beta]$, we say that $F[\beta]/F$ is a {\it parameter field\/} for $\theta$. The isomorphism class of $F[\beta]/F$ is {\it not\/} necessarily determined by $\theta$, as the following example shows. 
\remark{Example} 
Take $G = \GL pF$, $A = \M pF$, and let $\frak a$ be a {\it minimal\/} hereditary $\frak o_F$-order in $A$. We fix a simple stratum $[\frak a,1,0,\alpha]$ in $A$ (to use the full notation for once). If $E/F$ is a totally ramified extension of degree $p$, there is an embedding $\bsi:E\to A$, and an element $\beta\in E$, such that $\bsi E^\times \subset \scr K_\frak a$ and $[\frak a,1,0,\bsi\beta]$ is equivalent to $[\frak a,1,0,\alpha]$, that is, $\bsi\beta \equiv \alpha\pmod{\frak a}$. In particular, $\scr C(\frak a,\alpha) = \scr C(\frak a,\bsi\beta)$. So, if $\theta\in \scr C(\frak a,\alpha)$, then {\it every\/} totally ramified extension $E/F$, of degree $p$, is a parameter field for $\theta$. 
\endremark 
Returning to the general situation, it is necessary to extend the framework to include the notion of a {\it trivial\/} simple character: a trivial simple character is the trivial character of $U^1_\frak a$, for a hereditary $\frak o_F$-order $\frak a$ in $\M nF$, for some $n\ge1$. If $\theta$ is such a character, we define $\bk J_\theta$, $J^0_\theta$, $J^1_\theta$ as in (2.1.1), to obtain $\bk J_\theta = \scr K_\frak a$, $J^0_\theta = U_\frak a$, $J^1_\theta = U^1_\frak a$. We may treat this case within the same framework, simply by taking the view that $F[\beta] = F$. 
\subhead 
2.2 
\endsubhead 
We consider the class of all simple characters $\theta\in \scr C(\frak a,\beta,\psi_F)$, as $[\frak a,\beta]$ ranges over all simple strata in all matrix algebras $\M nF$, $n\ge1$. As in \cite{3}, this class carries an equivalence relation, called {\it endo-equivalence,\/} extending the relation of conjugacy (see (2.2.4) below). 
\par 
We let $\scr E(F)$ denote the set consisting of all endo-equivalence classes of simple characters, together with a trivial element $\bk 0_F$ which is the class of trivial simple characters. 
\par
If $\vT \in \scr E(F)$ and $\theta\in \vT$, we say that $\theta$ is a {\it realization\/} of $\vT$ and that $\vT$ is the {\it endo-class\/} of $\theta$. We use the notation $\vT = \text{\it cl\/}(\theta)$. From the definition of endo-equivalence \cite{3} and (2.1.1), we obtain: 
\proclaim\nofrills{(2.2.1)} \ 
Let $\vT\in \scr E(F)$, $\vT \neq \bk 0_F$. If $\theta_i \in \scr C(\frak a_i,\beta_i)$ is a realization of $\vT$, for $i =1,2$, then  
$$
e(F[\beta_1]|F) = e(F[\beta_2]|F), \quad f(F[\beta_1]|F) = f(F[\beta_2]|F), 
$$ 
and hence $[F[\beta_1]{:}F] = [F[\beta_2]{:}F]$. 
\endproclaim 
In this situation, we set 
$$
e(\vT) = e(F[\beta_i]|F), \quad f(\vT) = f(F[\beta_i]|F), \quad \deg\vT = [F[\beta_i]{:}F]. 
\tag 2.2.2 
$$ 
By convention, 
$$
e(\bk 0_F) = f(\bk 0_F) = \deg \bk 0_F = 1. 
\tag 2.2.3 
$$ 
If $\vT\in \scr E(F)$, $\vT\neq \bk 0_F$, and if $\theta\in \scr C(\frak a,\beta)$ is a realization of $\vT$, we say that $F[\beta]/F$ is a {\it parameter field\/} for $\vT$. By convention, $F$ is the only parameter field for $\bk 0_F$. In general, however, the parameter field is not uniquely determined (as in the Example of 2.1). 
\par 
Finally, we recall \cite{18 (3.5.11)}, \cite{3 (8.7)}:  
\proclaim\nofrills{(2.2.4)} \ 
For $i=1,2$, let $[\frak a,\beta_i]$ be a simple stratum in $\M nF$ and let $\theta_i\in \scr C(\frak a,\beta_i)$. The following conditions are equivalent: 
\roster 
\item $\theta_1$ intertwines with $\theta_2$ in $G$; 
\item $\theta_1$ is conjugate to $\theta_2$ in $G$; 
\item $\theta_1$ is endo-equivalent to $\theta_2$. 
\endroster 
\endproclaim 
\subhead 
2.3 
\endsubhead 
We recall from \cite{3} some fundamental properties of endo-classes of simple characters, relative to tamely ramified base field extension. We take the opportunity to fix some notation for the rest of the paper. 
\par 
From \cite{3 (7.1, 7.7)}, we obtain the following facts. 
\proclaim\nofrills{(2.3.1)} \ 
Let $[\frak a,\beta]$ be a simple stratum in $A = \M nF$. Let $K/F$ be a subfield of $A$ satisfying the following conditions: 
\roster
\item"\rm (a)" $K/F$ is tamely ramified, 
\item"\rm (b)" $K$ commutes with $\beta$, and 
\item"\rm (c)" the $K$-algebra $K[\beta]$ is a subfield of $A$ such that $K[\beta]^\times\subset \scr K_\frak a$. 
\endroster 
Let $A_K$ be the $A$-centralizer of $K$, and set $G_K = A_K^\times$, $\frak a_K = \frak a\cap A_K$. 
\roster 
\item 
The pair $[\frak a_K,\beta]$ is a simple stratum in $A_K$ and 
$$ 
Y(\beta,\frak a)\cap G_K = Y(\beta,\frak a_K), 
$$
where $Y$ stands for any of $\bk J$, $J^0$, $J^1$, $H^1$. 
\item 
Let $\theta\in \scr C(\frak a,\beta,\psi_F)$, and set $\psi_K = \psi_F\circ \roman{Tr}_{K/F}$. The character 
$$
\theta_K = \theta|_{H^1(\beta,\frak a_K)} 
\tag 2.3.2 
$$
lies in $\scr C(\frak a_K,\beta,\psi_K)$. 
\endroster 
\endproclaim 
We can work in the opposite direction, starting with a finite, tamely ramified field extension $K/F$ and a non-trivial endo-class $\varPhi\in \scr E(K)$. We choose a realization $\vf$ of $\vF$, say $\varphi\in \scr C(\frak b,\gamma,\psi_K)$ where $[\frak b,\gamma]$ is a simple stratum in $B = \End KV$ and $V$ is a finite-dimensional $K$-vector space. Set $A = \End FV$. There is then a unique hereditary $\frak o_F$-order $\frak a$ in $A$, stable under conjugation by $K^\times$ and such that $\frak b = \frak a\cap B$. 
\proclaim\nofrills{(2.3.3)} 
\roster 
\item 
There exists a simple stratum $[\frak b,\gamma']$ in $B$ such that 
\itemitem{\rm (a)} $\scr C(\frak b,\gamma',\psi_K) = \scr C(\frak b,\gamma,\psi_K)$ and 
\itemitem{\rm (b)} $[\frak a,\gamma']$ is a simple stratum in $A$. 
\item 
There exists a unique simple character $\vf^F\in \scr C(\frak a,\gamma',\psi_F)$ such that 
$$ 
\vf^F|_{H^1(\gamma',\frak b)} = \varphi. 
$$ 
\item 
The endo-class $\varPhi^F = \text{\it cl\/}(\vf^F)$ depends only on the endo-class $\varPhi$ of $\varphi$. In particular, it does not depend on the choices of $\gamma'$ and $\psi_F$. 
\item 
The map 
$$
\aligned 
\frak i_{K/F}:\scr E(K) &\longrightarrow \scr E(F), \\ 
\varPhi &\longmapsto \varPhi^F, 
\endaligned 
\tag 2.3.4 
$$
is surjective with finite fibres. When $K/F$ is Galois, the fibres of $\frak i_{K/F}$ are also the orbits of\/ $\Gal KF$ in $\scr E(K)$. 
\item 
If $F\subset K'\subset K$, then $\frak i_{K/F} = \frak i_{K'/F}\circ \frak i_{K/K'}$. 
\endroster 
\endproclaim 
These assertions combine 7.10 Theorem and 9.13 Corollary of \cite{3}. 
\par 
Let $\vT\in \scr E(F)$. The elements of the finite set $\frak i_{K/F}^{-1}\vT$ are called the {\it $K/F$-lifts of\/} $\vT$. The only $K/F$-lift of $\bk 0_F$ is $\bk 0_K$. One deduces from \cite{3 (11.2)} the following method of computing the basic invariants (2.2.2) of the lifts of an endo-class. 
\proclaim\nofrills{(2.3.5)} \ 
Let $\vT\in \scr E(F)$ have parameter field $P/F$. Let $K/F$ be a finite, tamely ramified field extension. Write $
K\otimes_FP = P_1\times P_2\times \dots \times P_r$, where each $P_i/K$ is a field extension. There is a canonical bijection 
$$ 
\align 
\{P_1,P_2,\dots,P_r\} &\longrightarrow \frak i_{K/F}^{-1}\vT, \\ 
P_i&\longmapsto \vT_i, 
\endalign 
$$ 
such that $P_i/K$ is a parameter field for $\vT_i$. This bijection is natural with respect to isomorphisms $K/F \to K'/F'$ of field extensions. 
\endproclaim 
In the context of (2.3.1), the field $K[\beta]$ is $K$-isomorphic to one, and only one, of the simple components $P_i$ of $K\otimes_FF[\beta]$. The endo-class, in $\scr E(K)$, of the character $\theta_K$ of (2.3.2)  is the $K/F$-lift of $\text{\it cl\/}(\theta)$ corresponding to the component $P_i$. 
\subhead 
2.4 
\endsubhead 
Let $\vT\in \scr E(F)$. We say that $\vT$ is {\it totally wild\/} if $\deg\vT = e(\vT) = p^r$, for an integer $r\ge0$. Equivalently, $\vT$ is totally wild if any parameter field for $\vT$ is totally wildly ramified over $F$. 
\proclaim{Proposition} 
Let $\vT\in \scr E(F)$ have a parameter field $P/F$, and let $T/F$ be the maximal tamely ramified sub-extension of $P/F$. Let $K/F$ be a finite, tamely ramified field extension. The following conditions are equivalent: 
\roster 
\item 
$\vT$ has a totally wild $K/F$-lift; 
\item 
there exists an $F$-embedding $T\to K$. 
\endroster 
If $P'/F$ is a parameter field for $\vT$ and $T'/F$ is the maximal tamely ramified sub-extension of $P'/F$, then $T'$ is $F$-isomorphic to $T$. 
\endproclaim 
\demo{Proof} 
All assertions follow readily from (2.3.5). \qed 
\enddemo 
We therefore define a tame parameter field for $\vT\in \scr E(F)$ to be a finite tame extension $E/F$ such that $\vT$ has a totally wild $E/F$-lift, and such that $E$ is minimal for this property. The proposition shows that $\vT$ has a tame parameter field, and it is uniquely determined up to $F$-isomorphism (unlike the general parameter field). 
\par
The group $\roman{Aut}(T|F)$ acts on $\scr E(T)$. Using (2.3.5) again, the proposition yields: 
\proclaim{Corollary} 
Let $\vT\in \scr E(F)$ have tame parameter field $T/F$. Let $\varPsi\in \scr E(T)$ be a totally wild $T/F$-lift of $\vT$. The map $\gamma\mapsto \varPsi^\gamma$ is a bijection between $\roman{Aut}(T|F)$ and the set of totally wild $T/F$-lifts of $\vT$. 
\endproclaim 
\subhead 
2.5 
\endsubhead 
Some new terminology will be convenient. Let $[\frak a,\beta]$ be a simple stratum in $A = \M nF$; we say that $[\frak a,\beta]$ is {\it m-simple\/} if $\frak a$ is maximal among $\frak o_F$-orders in $A$ stable under conjugation by $F[\beta]^\times$. 
\proclaim{Lemma 1} 
Let $[\frak a,\beta]$ be a simple stratum in $A$. The following conditions are equivalent. 
\roster 
\item 
$[\frak a,\beta]$ is m-simple; 
\item 
the order $\frak a$ is principal and the ramification index $e(F[\beta]|F)$ is equal to the $\frak o_F$-period of $\frak a$;  
\item 
the hereditary $\frak o_{F[\beta]}$-order $\frak a_{F[\beta]} = \frak a\cap A_{F[\beta]}$ is maximal. 
\endroster 
\endproclaim 
\demo{Proof} 
See \cite{18 (1.2.4)}. \qed
\enddemo 
A simple character $\theta$ is deemed m-simple if $\theta\in \scr C(\frak a,\beta)$ where $[\frak a,\beta]$ is m-simple. If $\vT\in \scr E(F)$, an m-realization of $\vT$ is a realization which is m-simple. 
\proclaim{Lemma 2}  
Two m-simple characters in $G = \GL nF$ are endo-equivalent if and only if they are conjugate in $G$. 
\endproclaim 
\demo{Proof} 
This is a special case of (2.2.4). \qed 
\enddemo 
Let $[\frak a,\beta]$ be an m-simple stratum in $A$ and set $P = F[\beta]$. By condition (3) of Lemma 1, the order $\frak a_P = \frak a \cap A_P$ is isomorphic to $\M s{\frak o_P}$, where $s = n/[P{:}F]$. We further have 
$$ 
\aligned 
\bk J(\beta,\frak a) &= P^\times J^0(\beta,\frak a), \\ 
J^0(\beta,\frak a)/J^1(\beta,\frak a) &\cong \GL s{\Bbbk_P}. 
\endaligned 
\tag 2.5.1 
$$ 
\subhead 
2.6 
\endsubhead 
We need a more concrete version of the notion of tame parameter field. Let $\theta$ be a simple character in $\GL nF$. A {\it tame parameter field\/} for $\theta$ is a subfield $T/F$ of $\M nF$ of the following form: there is a simple stratum $[\frak a,\beta]$ in $\M nF$ such that $\theta\in \scr C(\frak a,\beta)$ and $T/F$ is the maximal tamely ramified sub-extension of $F[\beta]/F$. The proposition of 2.4 implies that a tame parameter field for $\theta$ is also a tame parameter field for $\text{\it cl\/}(\theta)$, and so is uniquely determined up to $F$-isomorphism. There is a more precise version of this uniqueness property. 
\proclaim{Proposition} 
Let $\theta$ be an m-simple character in $G = \GL nF$. Let $T_1/F$, $T_2/F$ be tame parameter fields for $\theta$. 
\roster 
\item There exists $j\in J^1_\theta$ such that $T_2 = T_1^j$. 
\item 
An element $y\in J^1_\theta$ normalizes $T_1^\times$ if and only if it centralizes $T_1^\times$. 
\endroster 
In particular, there is a unique $F$-isomorphism $T_1\to T_2$ implemented by conjugation by an element of $J^1_\theta$. 
\endproclaim 
\demo{Proof} 
Before starting the proof, it will be useful to recall a standard technical result proved in the same way as, for example, \cite{3 (15.19)}.  
\proclaim{Conjugacy Lemma} 
Let $\scr G$ be a pro $p$-group and let $\alpha$ be an automorphism of $\scr G$, of finite order relatively prime to $p$. Suppose $\scr G$ admits a descending chain $\{\scr G_i\}_{i\ge1}$ of open, normal, $\alpha$-stable subgroups $\scr G_i$ such that $\scr G_i/\scr G_{i+1}$ is a finite abelian $p$-group and $[\scr G_i,\scr G_j] \subset \scr G_{i+j}$, for all $i,j\ge 1$. In the group $\langle\alpha\rangle \ltimes \scr G$, any element $\alpha g$, $g\in \scr G$, is $\scr G$-conjugate to one of the form $\alpha h$, where $h\in \scr G$ commutes with $\alpha$. 
\endproclaim 
We prove the proposition. We choose a simple stratum $[\frak a,\beta_i]$ in $A = \M nF$ such that $\theta\in \scr C(\frak a,\beta_i,)$ and $T_i \subset P_i = F[\beta_i]$. Thus $J^1_\theta = J^1(\beta_i,\frak a)$, and so on. We abbreviate $J^1 = J^1_\theta$. 
\par 
Let $K_i/F$ be the maximal unramified sub-extension of $T_i/F$. Consider the ring $\frak J = \frak J(\beta_i,\frak a)$ (in the notation of \cite{18 (3.1.8)}). If we write $\frak J^1 = \frak J\cap \frak p_\frak a$, then $J^1 = 1{+}\frak J^1$ and the ring $\frak j = \frak J/\frak J^1$ is isomorphic to $\M s{\Bbbk_{P_i}}$, where $s[P_i{:}F] = n$. The centre $\frak z$ of $\frak j^\times$ is therefore $\Bbbk_{P_i}^\times$. Since $\frak o_{P_i}\cap \frak J^1 = \frak p_{P_i}$, reduction modulo $\frak p_{P_i}$ induces an isomorphism $\bs\mu_{P_i} \cong \frak z$. Since $P_i/K_i$ is totally ramified and $\frak J^1\cap \frak o_{K_i} = \frak p_{K_i}$, the group $\bs\mu_{K_i}$ also maps isomorphically to $\frak z$ under reduction modulo $\frak J^1$. So, if $\zeta_1$ is a generator of $\bs\mu_{K_1}$, there exists $\zeta_2\in \bs\mu_{K_2}$ such that $\zeta_2 = \zeta_1 j$, for some $j\in J^1$. The roots of unity $\zeta_1$, $\zeta_2$ have the same order, so $\zeta_2$ generates $\bs\mu_{K_2}$. 
\par 
We apply the Conjugacy Lemma to $\scr G = J^1$ (with its standard filtration) and the automorphism $\alpha$ induced by conjugation by $\zeta_1$. We deduce that $\zeta_2 = \zeta_1 j$ is $J^1$-conjugate to an element $\zeta_3 = \zeta_1 j_1$, where $j_1\in J^1$ commutes with $\zeta_1$. Thus $j_1 = \zeta_1^{-1}\zeta_3$ has finite order dividing that of $\zeta_1$, and this order is not divisible by $p$. Since $J^1$ is a pro-$p$ group, we deduce that $j_1 = 1$. 
\par 
Therefore, after applying a $J^1$-conjugation to $T_2$ (which does not affect the conditions imposed on $T_2$), we may assume that $K_1 = K_2 = K$, say. Let $A_K$ be the centralizer of $K$ in $\M nF$, set $G_K = A_K^\times$, and put $\frak a_K = \frak a\cap A_K$. As in (2.3.1), the pair $[\frak a_K,\beta_i]$ is a simple stratum in $A_K$ and $\scr C(\frak a_K,\beta_1) = \scr C(\frak a_K,\beta_2)$. Moreover, $\bk J(\beta_i,\frak a_K) = \bk J(\beta_i,\frak a)\cap G_K = \bk J_K$, say. Also, $J^1(\beta_i,\frak a_K) = J^1\cap G_K = J^1_K$, say. 
\par  
The extensions $T_1/K$, $T_2/K$ are isomorphic and totally tamely ramified. Let $\vp$ be a prime element of $T_1$ such that $\vp^{[T_1{:}K]} = \vp_K$, for some prime element $\vp_K$ of $K$. Since we are in the m-simple case, the centre of the group $\bk J_K/J^1_K$ is  $P_i^\times J_K^1/J_K^1$. The groups $T_1^\times$, $T_2^\times$ have the same image in $\bk J(\beta,\frak a_K)/ J^1(\beta,\frak a_K)$, and we can argue exactly as before to complete the proof of (1). 
\par 
In part (2), write $T = T_1$ and take $K$ as before. Let $j\in J^1$ normalize $T^\times$. Surely $\roman{Ad}\,j$ acts trivially on the group $\bs\mu_K = \bs\mu_T$, whence $j$ centralizes $K^\times$. That is, $j\in J^1_K = J^1(\beta,\frak a_K)$. Taking a prime element $\vp$ as before, $j\vp j^{-1} = \zeta\vp$, for some $\zeta\in \bs\mu_K$. However, $j\vp j^{-1}\equiv \vp \pmod{J^1_K}$, whence $\zeta = 1$, as required. \qed 
\enddemo 
The condition of m-simplicity imposed on $\theta$ is unnecessary. However, it shortens the arguments a little and is the only case we need. 
\remark{Remark} 
In the situation of the proposition, there may be a field extension $T_0/F$, with $T_0^\times \subset \bk J_\theta$, which is isomorphic to a tame parameter field for $\theta$ without being one. See 10.1 Remark for a context in which such examples arise. 
\endremark 
\subhead 
2.7 
\endsubhead 
We connect our two uses of the phrase ``tame parameter field''. As a matter of notation, if we have an isomorphism $\bsi:F\to F'$ of fields, we denote by $\bsi_*$ the induced bijection $\scr E(F) \to \scr E(F')$. 
\proclaim{Proposition} 
Let $\vT\in \scr E(F)$ have tame parameter field $E/F$, and let $\varPsi$ be a totally wild $E/F$-lift of $\vT$. Let $d = \deg\vT$ and $n = rd$, for an integer $r\ge1$. Let $\bsi$ be an $F$-embedding of $E$ in $A = \M nF$. There exists an m-simple character $\theta$ in $G = \GL nF$, with the following properties: 
\roster 
\item"\rm (a)" $\text{\it cl\/}(\theta) = \vT$; 
\item"\rm (b)" the field $\bsi E$ is a tame parameter field for $\theta$; 
\item"\rm (c)" if $\theta_{\bsi E}$ is defined as in \rom{(2.3.2),} then $\text{\it cl\/}(\theta_{\bsi E}) = \bsi_*\varPsi$. 
\endroster 
These conditions determine $\theta$ uniquely, up to conjugation by an element of the $G$-centralizer $G_{\bsi E}$ of $\bsi E^\times$. 
\endproclaim 
\demo{Proof} 
Let $\theta\in \scr C(\frak a,\beta)$ be an m-realization of $\vT$, and let $T/F$ be the maximal tamely ramified sub-extension of $F[\beta]/F$. The fields $T$, $\bsi E$ are $F$-isomorphic and hence $G$-conjugate. So, replacing $\theta$ by a $G$-conjugate, we can assume $T = \bsi E$. The endo-class $\text{\it cl\/}(\theta_T)$ is then a totally wild lift of $\vT$. Let $g\in G$ normalize $\bsi E$. We have the relation $\text{\it cl\/}(\theta_{\bsi E}^g) = \text{\it cl\/}(\theta_{\bsi E})^g$. Invoking 2.4 Corollary, we may choose $\theta$ to achieve $\text{\it cl\/}(\theta_{\bsi E}) = \bsi_*\varPsi$. This character $\theta$ satisfies the conditions (a)--(c). 
\par
For the uniqueness assertion, let $\theta$, $\theta'$ be m-simple characters in $G$ satisfying conditions (a)--(c) relative to the field $\bsi E$. Condition (c) says that the m-simple characters $\theta_{\bsi E}$, $\theta'_{\bsi E}$ in $G_{\bsi E}$ are endo-equivalent over $\bsi E$, and so are $G_{\bsi E}$-conjugate (2.5 Lemma 2). We may as well assume, therefore, that $\theta_{\bsi E} = \theta'_{\bsi E}$. Condition (a) implies, via 2.5 Lemma 2, that $\theta' = \theta^g$, for some $g\in G$. The element $g$ conjugates $\bsi E$ to a tame parameter field for $\theta'$. Replacing $g$ by $gj$, for some $j\in J^1_{\theta'}$, we may assume that $(\bsi E)^g = \bsi E$ (2.6 Proposition). In other words, $g$ normalizes $\bsi E$. Since $\theta'$ satisfies (c),  2.4 Corollary implies that $g$ commutes with $\bsi E$, as required. \qed 
\enddemo  
\proclaim{Corollary} 
For $k=1,2$, let $\bsi^{(k)}$ be an $F$-embedding of $E$ in $A$ and let $\theta^{(k)}$ be an m-simple character in $G$, of endo-class $\vT$ and for which $\bsi^{(k)}E/F$ is a tame parameter field. Let $\Phi_k$ denote the endo-class of $\theta^{(k)}_{\bsi^{(k)}E}$. The pairs $(\bsi^{(k)},\theta^{(k)})$ are conjugate in $G$ if and only if 
$$
(\bsi^{(1)}_*)^{-1}\,\Phi_1 = (\bsi^{(2)}_*)^{-1}\,\Phi_2. 
$$
\endproclaim 
\demo{Proof} 
Replacing one pair by a $G$-conjugate, we can assume $\bsi^{(1)} = \bsi^{(2)} = \bsi$ say. By the proposition, the characters $\theta^{(j)}$ give rise to the same endo-class over $E$ if and only if they are $G_{\bsi E}$-conjugate. \qed 
\enddemo 
\head\Rm 
3. Action of tame characters 
\endhead 
We take a (possibly trivial) m-simple character $\theta$ in $G = \GL nF$, and let $T/F$ be a tame parameter field for $\theta$. We analyze various families of irreducible representations of $\bk J_\theta$ which contain $\theta$ and show that these families carry canonical actions of the group $X_1(T)$. We follow throughout the notational conventions introduced in 2.3. 
\subhead 
3.1 
\endsubhead 
Let $\theta$ be an m-simple character in $G = \GL nF$ and let $I_G(\theta)$ denote the set of elements of $G$ which intertwine $\theta$. 
\proclaim{Definition} 
A character of $\bk J_\theta$ is {\rm $\theta$-flat\/} if it is trivial on $J^1_\theta$ and is intertwined by every element of $I_G(\theta)$. 
\endproclaim 
For example, if $\chi\in X_1(F)$, then $\chi\circ\det|_{\bk J_\theta}$ is a $\theta$-flat character of $\bk J_\theta$. Let $X_1(\theta)$ denote the group of $\theta$-flat characters of $\bk J_\theta$. 
\proclaim{Proposition} 
Let $T/F$ be a tame parameter field for $\theta$ and let\/ $\det_T:G_T\to T^\times$ be the determinant map. 
\roster 
\item 
Let $\phi\in X_1(T)$. There is a unique $\theta$-flat character $\phi^\bk J$ of $\bk J = \bk J_\theta$ such that $\phi^\bk J(x) = \phi(\det_T x)$, for all $x\in \bk J\cap G_T$. 
\item 
The map 
$$
\aligned 
X_1(T) &\longrightarrow X_1(\theta), \\ 
\phi&\longmapsto \phi^\bk J, 
\endaligned 
\tag 3.1.1 
$$
is a surjective homomorphism of abelian groups. Its kernel is the group $X_0(T)_s$ of unramified characters $\chi$ of\/ $T^\times$ satisfying $\chi^s=1$, where $s = n/\deg\text{\it cl\/}(\theta)$. 
\item 
Let $T'/F$ be a tame parameter field for $\theta$. For $\xi\in X_1(T')$, define $\xi^\bk J\in X_1(\theta)$ as in \rom{(1),} using $T'$ in place of\/ $T$. If $j\in J^1_\theta$ satisfies $T' = T^j$, then 
$$
\big(\phi^j\big)^\bk J = \phi^\bk J,\quad \phi\in X_1(T). 
$$
\endroster 
\endproclaim 
\demo{Proof} 
By definition, there is a simple stratum $[\frak a,\beta]$ in $A$ such that $\theta\in \scr C(\frak a,\beta)$ and $T\subset F[\beta]$. Set $P = F[\beta]$, and use the standard notation of 2.3. Thus $\bk J_\theta = P^\times U_{\frak a_P} J^1_\theta$ (2.5.1), and $P^\times U_{\frak a_P}\cap J^1_\theta = P^\times U_{\frak a_P} \cap U^1_\frak a = U^1_{\frak a_P}$. Let $\det_P:G_P\to P^\times$ denote the determinant map.
\par 
Let $\chi\in X_1(P)$ and define a character $\pre{\bk J}\chi$ of $\bk J_\theta$ by 
$$
\pre{\bk J}\chi(xj) = \chi({\det}_P\,x), \quad x\in P^\times U_{\frak a_P},\ j\in J^1_\theta. 
$$
We show that $\pre{\bk J}\chi\in X_1(\theta)$. Surely $\pre{\bk J}\chi$ is trivial on $J^1_\theta$. The set $I_G(\theta)$ is $J^1_\theta G_P J^1_\theta$. If $y\in G_P$, then $y^{-1} \bk J_\theta y\cap \bk J_\theta = P^\times(y^{-1} J^0_\theta y\cap J^0_\theta)$. The characters $\pre{\bk J}\chi$, $(\pre{\bk J}\chi)^y$ agree on $P^\times$, and we are so reduced to showing that $y$ intertwines $\pre{\bk J}\chi|_{J^0_\theta}$. This will follow from: 
\proclaim{Lemma} 
Set $K = U_{\frak a_P}U^1_\frak a$, and define a character $\bar\chi$ of $K$ by 
$$
\bar\chi(hu) = \chi({\det}_P\,h), \quad h\in U_{\frak a_P},\ u\in U^1_\frak a. 
$$
The character $\bar\chi$ is intertwined by every element $y$ of $G_P$. 
\endproclaim 
\demo{Proof} 
Let $y\in G_P$. We show 
$$
K\cap K^y = (U_{\frak a_P}\cap U_{\frak a_P}^y)(U^1_\frak a\cap {U^1_\frak a}^y). 
\tag 3.1.2 
$$ 
The characters $\bar\chi$, $\bar\chi^y$ agree on each factor, so the lemma will follow. The relation (3.1.2) is equivalent to the additive relation 
$$
(\frak a_P + \frak p_\frak a)\cap (\frak a_P + \frak p_\frak a)^y = \frak a_P\cap \frak a_P^y + \frak p_\frak a\cap \frak p_\frak a^y, 
\tag 3.1.3 
$$ 
where $\frak p_\frak a$ is the Jacobson radical of $\frak a$. Techniques to deal with such identities are developed in the early chapters of \cite{18}, but we give an {\it ad hoc\/} argument of the same kind. We view $A_P$ as $\End PV$, where $V$ is a $P$-vector space of dimension $s$. Let $L$ be an $\frak a_P$-lattice in $V$. We choose an $\frak o_P$-basis of $L$ to identify $A_P$ with $\M sP$ and $\frak a_P$ with $\M s{\frak o_P}$. Let $A(P) = \End FP$ and let $\frak a(P)$ be the unique hereditary $\frak o_F$-order in $A(P)$ stable under conjugation by $P^\times$. The Jacobson radical of $\frak a(P)$ is then $\vp\frak a(P)$, where $\vp$ is a prime element of $\frak o_P$. The order $\frak a$ is $\M s{\frak a(P)}$, and its radical is $\vp \frak a = \M s{\vp \frak a(P)}$. The desired properties (3.1.2), (3.1.3) depend only on the double coset $U_{\frak a_P} y U_{\frak a_P}$. We may therefore take $y$ diagonal, with all its eigenvalues powers of $\vp$. The relation (3.1.3) follows immediately. \qed 
\enddemo 
We have shown that $\chi\mapsto \pre{\bk J}\chi$ gives a homomorphism $X_1(P) \to X_1(\theta)$. Further, $\pre{\bk J}\chi=1$ if and only if $\chi$ is trivial on $\det_P P^\times U_{\frak a_P} = (P^\times)^sU_P$, that is, if and only if $\chi \in X_0(P)_s$. 
\par 
On the other hand, if $\xi\in X_1(\theta)$, then the restriction of $\xi$ to $\bk J_\theta \cap G_P  = P^\times U_{\frak a_P}$ is intertwined by every element of $G_P$. This implies that $\xi|_{\bk J_\theta\cap G_P}$ factors through $\det_P$ ({\it cf\.} \cite{18 (2.4.11)}) and so $\xi$ lies in the image of $X_1(P)$. The map $X_1(P)\to X_1(\theta)$, $\chi\mapsto \pre{\bk J}\chi$, is therefore surjective and induces an isomorphism $X_1(P)/ X_0(P)_s \to X_1(\theta)$. 
\par 
We next observe that $\bk J_\theta = (\bk J_\theta\cap G_T)\,J^1_\theta$, and $G_T\cap J^1_\theta \subset U^1_{\frak a_T}$. The field extension $P/T$ is totally wildly ramified, so the map $\phi\mapsto \phi\circ\N PT$ gives an isomorphism $X_1(T) \to X_1(P)$ carrying $X_0(T)_s$ to $X_0(P)_s$. For $\phi\in X_1(T)$, the character $\pre{\bk J}(\phi\circ\N PT)$ of $\bk J_\theta$ satisfies the defining property for $\phi^\bk J$, and we have proven the first two parts of the proposition. Part (3) is an immediate consequence of the definitions. \qed 
\enddemo 
\remark{Remark} 
The identification of $X_1(\theta)$ with a quotient of $X_1(P)$ (as in the proof of the proposition) is convenient for calculation. It cannot be regarded as canonical, since $P$ is not unambiguously defined by $\theta$. However, the proposition shows that the identification of $X_1(\theta)$ with $X_1(T)/X_0(T)_s$ is canonical. 
\endremark 
\subhead 
3.2 
\endsubhead 
We continue with the m-simple character $\theta$ of 3.1. Let $\eta$  be the unique irreducible representation of $J^1_\theta$ which contains $\theta$, as in \cite{18 (5.1.1)}. We refer to $\eta$ as the $1${\it -Heisenberg representation over $\theta$.} 
\proclaim{Definition} 
A {\rm full Heisenberg representation over $\theta$} is a representation $\kappa$ of $\bk J_\theta$ such that 
\roster 
\item"\rm (a)" $\kappa|_{J^1_\theta} \cong \eta$ and 
\item"\rm (b)" $\kappa$ is intertwined by every element of $I_G(\theta)$. 
\endroster 
\endproclaim 
We denote by $\scr H(\theta)$ the set of equivalence classes of full Heisenberg representations over $\theta$.  We observe that, if $\kappa\in \scr H(\theta)$ and $\phi\in X_1(\theta)$, then $\phi\otimes \kappa$ also lies in $\scr H(\theta)$. Thus $\scr H(\theta)$ comes equipped with a canonical action of $X_1(\theta)$. 
\proclaim{Proposition} 
Let $\theta$ be an m-simple character in $G = \GL nF$. 
\roster 
\item 
There exists a full Heisenberg representation over $\theta$. 
\item 
If $\kappa\in \scr H(\theta)$ and $\phi\in X_1(\theta)$, then $\phi\otimes\kappa\in \scr H(\theta)$ and the pairing 
$$
\align 
X_1(\theta) \times \scr H(\theta) &\longrightarrow \scr H(\theta), \\ 
(\phi,\kappa) &\longmapsto \phi\otimes\kappa, 
\endalign
$$ 
endows $\scr H(\theta)$ with the structure of a principal homogeneous space over $X_1(\theta)$. 
\endroster 
\endproclaim 
\demo{Proof} 
The proof of part (1) requires a different family of techniques which we do not use elsewhere. We therefore give it separate treatment in 3.3 below. Here we deduce part (2). 
\par 
We need an intermediate step. Let $[\frak a,\beta]$ be a simple stratum such that $\theta\in \scr C(\frak a,\beta)$, and set $P = F[\beta]$. We abbreviate $J^1 = J^1_\theta$, and so on. Thus $I_G(\theta) = J^1G_PJ^1$. Let $\scr H^0(\theta)$ be the set of equivalence classes of representations $\mu$ of $J^0$ satisfying $\mu|_{J^1} \cong \eta$ and such that $I_G(\mu) = I_G(\theta)$: in the language of \cite{18}, the elements of $\scr H^0(\theta)$ are the ``$\beta$-extensions of $\eta$''. 
\par 
Let $X_1^0(\theta)$ denote the group of characters of $J^0$ which are trivial on $J^1$ and are intertwined by every element of $I_G(\theta)$. 
\proclaim{Lemma} 
\roster 
\item 
Let $\chi\in X_1^0(\theta)$. The restriction of $\chi$ to $U_{\frak a_P} = J^0\cap G_P$ factors through $\det_P$, and $\chi$ is the restriction to $J^0$ of some $\tilde\chi\in X_1(\theta)$. 
\item 
The set $\scr H^0(\theta)$ is a principal homogeneous space over $X_1^0(\theta)$. Distinct elements of $\scr H^0(\theta)$ do not intertwine in $G$. 
\endroster 
\endproclaim 
\demo{Proof} 
In (1), 
the restriction of $\chi$ to $U_{\frak a_P}$ is intertwined by every element of $G_P$, and so must factor through $\det_P$. The second assertion follows from the description of elements of $X_1(\theta)$ given in the proof of 3.1 Proposition. Part (2) is part (iii) of \cite{18 (5.2.2)}. \qed 
\enddemo 
Let $\kappa, \kappa'\in \scr H(\theta)$. By the lemma, there exists $\chi\in X_1(\theta)$ such that $\chi\otimes\kappa$ agrees with $\kappa'$ on $J^0$. Thus there exists a character $\phi$ of $\bk J$, trivial on $J^0$, such that $\phi\chi\otimes\kappa \cong \kappa'$. The character $\phi\chi$ lies in $X_1(\theta)$, so the set $\scr H(\theta)$ comprises a single $X_1(\theta)$-orbit. 
\par 
Next, let $\kappa\in \scr H(\theta)$, $\chi\in X_1(\theta)$, and suppose that $\kappa\cong \chi\otimes\kappa$. The lemma implies that $\chi$ is trivial on $J^0$. We view $\kappa$ and $\chi\otimes\kappa$ as acting on the same vector space $V$. By hypothesis, there is an automorphism $f$ of $V$ such that $\chi(x)\kappa(x)\circ f = f\circ \kappa(x)$, for $x\in \bk J$. In particular, $\eta(y)\circ f = f\circ \eta(y)$, $y\in J^1$. Thus $f$ is a scalar, whence $\chi = 1$, as required. 
\par 
We have shown that the set $\scr H(\theta)$, if non-empty, is a principal homogeneous space over $X_1(\theta)$. \qed 
\enddemo 
Using 3.1 Proposition, we can equally view $\scr H(\theta)$ as an $X_1(T)$-space. We use the notation 
$$
\xi\odot\kappa = \xi^\bk J\otimes\kappa, 
\tag 3.2.1 
$$
for $\xi\in X_1(T)$ and $\kappa\in \scr H(\theta)$. 
\proclaim{Corollary} 
With respect to the action \rom{(3.2.1),} the set $\scr H(\theta)$ is a principal ho\-mogeneous space over $X_1(T)/X_0(T)_s$, where $s = n/\deg\text{\it cl\/}(\theta)$. 
\endproclaim 
\demo{Proof} 
The corollary follows directly from the proposition and 3.1 Proposition. \qed 
\enddemo 
\remark{Remark} 
Let $\tau$ be an irreducible representation of $\bk J_\theta$ such that $\tau|_{H^1_\theta}$ contains $\theta$. The proof of the proposition shows that $\tau\in \scr H(\theta)$ if and only if $\tau|_{J^0_\theta} \in \scr H^0(\theta)$. 
\endremark 
\subhead 
3.3 
\endsubhead 
We prove part (1) of 3.2 Proposition. Specifically, we are given an m-simple stratum $[\frak a,\beta]$ in $A = \M nF$ and a simple character $\theta\in \scr C(\frak a,\beta)$. Let $\eta$ be the $1$-Heisenberg representation over $\theta$. Writing $G = \GL nF$ and $P = F[\beta]$, so that $G_P$ is the $G$-centralizer of $P^\times$, we have $I_G(\theta) = J^0G_PJ^0$. We have to prove: 
\proclaim{Proposition} 
There exists a representation $\kappa$ of $\bk J_\theta$ such that $\kappa|_{J^1_\theta} \cong \eta$ and which is intertwined by every element of $G_P$. 
\endproclaim 
\demo{Proof}
As recalled in 3.2, there exists a representation $\mu$ of $J^0 = J^0_\theta$, extending $\eta$ and which is intertwined by every element of $G_P$. This surely admits extension to a representation of $\bk J = P^\times J^0$. If $[P{:}F] = n$, then $G_P = P^\times$ and there is nothing to prove. 
\par
To deal with the general case, we first recall from \cite{18} and \cite{3} one method of constructing the representation $\mu$. We abbreviate $\frak b = \frak a_P$ and let $L$ be a $\frak b$-lattice in $V = F^n$. If we choose a prime element $\vp$ of $P$, then $\{\vp^jL:j\in \Bbb Z\}$ is the $\frak o_F$-lattice chain in $V$ defining the hereditary order $\frak a$. Let $s = n/[P{:}F]$. We choose a decomposition 
$$
L = L_0 \oplus L_0 \oplus\ \cdots\ \oplus L_0 
$$ 
of $L$, in which $L_0$ is an $\frak o_P$-lattice of rank one. We accordingly write  
$$
V = V_0\oplus V_0 \oplus\ \cdots\ \oplus V_0. 
$$ 
Let $M$ be the subgroup of $G$ preserving this decomposition of $V$. In particular, $M$ is a Levi component of a parabolic subgroup $Q_+$ of $G$. Let $Q_-$ be the $M$-opposite of $Q_+$, and let $N_\pm$ be the unipotent radical of $Q_\pm$. The obvious projections $V\to V_0$ lie in $\frak a$, so we get an Iwahori decomposition 
$$
U^1_\frak a = U^1_\frak a\cap N_-\ \cdot\ U^1_\frak a\cap M\ \cdot\ U^1_\frak a\cap N_+. 
\tag 3.3.1 
$$
Moreover, $U^1_\frak a\cap M$ is the direct product of $s$ groups $U^1_{\frak a_0}$, where $\frak a_0$ is the hereditary $\frak o_F$-order in $\End F{V_0}$ defined by the lattice chain $\{\vp^jL_0:j\in \Bbb Z\}$. The pair $[\frak a_0,\beta]$ is an m-simple stratum in $\End F{V_0}$. 
\par 
The group $U_\frak a$ does not admit an Iwahori decomposition in general, but we do have 
$$
\align 
U_\frak a\cap Q_+ &= U_\frak a\cap M\ \cdot\ U_\frak a\cap N_+, \\ 
U_\frak a\cap M &= U_{\frak a_0}\times U_{\frak a_0} \times\ \cdots\ \times U_{\frak a_0}. 
\endalign 
$$ 
\proclaim{Lemma 1} 
Let $H^1 = H^1_\theta = H^1(\beta,\frak a)$ and similarly for $J^1$, $J^0$. 
\roster 
\item 
There are Iwahori decompositions 
$$
\align 
H^1 &= H^1\cap N_-\ \cdot\ H^1\cap M\ \cdot\ H^1\cap N_+, \\ 
H^1\cap M &= H^1(\beta,\frak a_0) \times H^1(\beta,\frak a_0) \times\ \cdots \ \times H^1(\beta,\frak a_0), 
\intertext{and} 
J^1 &= J^1\cap N_-\ \cdot\ J^1\cap M\ \cdot\ J^1\cap N_+, \\ 
J^1\cap M &= J^1(\beta,\frak a_0) \times J^1(\beta,\frak a_0) \times\ \cdots \ \times J^1(\beta,\frak a_0). 
\endalign 
$$ 
\item 
The characters $\theta|_{H^1\cap N_-}$, $\theta|_{H^1\cap N_+}$ are trivial, while there is a unique character $\theta_0\in \scr C(\frak a_0,\beta)$ such that 
$$
\theta|_{H^1\cap M} = \theta_0\otimes \theta_0\otimes\ \cdots\ \otimes \theta_0. 
$$ 
\item 
The set 
$$
J^1_+ = H^1\cap N_-\ \cdot\ J^1\cap Q_+ 
$$
is a group. Let $\eta_0$ be the $1$-Heisenberg representation of $J^1(\beta,\frak a_0)$ over $\theta_0$. There is a unique representation $\eta_+$ of $J^1_+$, which is trivial on $H^1\cap N_-$ and $J^1\cap N_+$, and which restricts to $\eta_0\otimes\eta_0\otimes\dots \otimes \eta_0$ on $J^1\cap M$. This representation satisfies 
$$ 
\Ind_{J^1_+}^{J^1}\,\eta_+ \cong \eta. 
$$ 
\endroster 
\endproclaim 
\demo{Proof} 
This is an instance of the discussion in \cite{3 \S10 Example 10.9}. \qed 
\enddemo 
\remark{Remark} 
The simple characters $\theta$, $\theta_0$ are endo-equivalent. 
\endremark 
For the next step, we take $\mu_0 \in \scr H^0(\theta_0)$ (in the notation of 3.2). The group $H^1\cap N_-\cdot J^0\cap Q_+$ admits a unique representation $\mu_+$ which is trivial on $H^1\cap N_-$ and $J^0\cap N_+$, and which restricts to $\mu_0\otimes\mu_0\otimes \dots \otimes \mu_0$ on $J^0\cap M$. 
\proclaim{Lemma 2} 
\roster 
\item 
The representation $\tilde\mu_+$ of $J^1\cap N_-\cdot J^0\cap Q_+$, induced by $\mu_+$, is irreducible, and satisfies $\tilde\mu_+|_{J^1} \cong \eta$. 
\item 
There is a unique representation $\mu$ of $J^0$ extending $\tilde\mu_+$. 
\item 
The representation $\mu$ lies in $\scr H^0(\theta)$. 
\endroster 
\endproclaim 
\demo{Proof} 
Assertion (1) is straightforward. The restriction of $\tilde\mu_+$ to $H^1\cap N_-\cdot J^1\cap M\cdot J^0\cap N_+$ is the representation denoted $\tilde\eta_\roman M$ in \cite{18 (5.1.14)}. It admits extension to a representation $\mu'$ of $J^0$, and any such extension lies in $\scr H^0(\theta)$ \cite{18 (5.2.4)}. Inspection of intertwining implies that $\mu'|_{J^0\cap M}$ is of the form $\mu'_0\otimes\dots\otimes\mu'_0$, for some $\mu'_0 \in \scr H^0(\theta_0)$. By 3.2 Lemma, we may choose $\mu'$ such that $\mu'_0 \cong \mu_0$. This proves (2) and (3). \qed 
\enddemo 
As remarked at the beginning of the proof, there exists $\kappa_0\in \scr H(\theta_0)$ such that $\kappa_0|_{J^0(\beta,\frak a_0)} \cong \mu_0$. We use $\kappa_0$ to define a representation $\kappa_+$ of the group $H^1\cap N_-\cdot P^\times J^0\cap Q_+$, extending $\mu_+$ and agreeing with $\kappa_0\otimes \dots \otimes \kappa_0$ on $P^\times J^0\cap M$. We then induce $\kappa_+$ to get a representation $\tilde\kappa_+$ of $J^1\cap N_-\cdot P^\times J^0\cap Q_+$ extending $\tilde\mu_+$. 
\par 
The representation $\mu$ of Lemma 2 is stable under conjugation by $P^\times$ and so admits extension to $\bk J$. This extension is uniquely determined up to tensoring with a character of $\bk J/J^0 \cong P^\times/U_P$. We may therefore choose this extension, call it $\kappa$, to agree with $\tilde\kappa_+$ on $J^1\cap N_-\cdot P^\times J^0\cap Q_+$. 
\par 
To complete the proof, we have to show that $\kappa$ is intertwined by any element $y$ of $I_G(\theta) = J^0 G_P J^0$. It will only be the coset $J^0 yJ^0$ which intervenes, so we may assume $y\in G_P$. Indeed, $J^0G_PJ^0 = J^0(G_P\cap M)J^0$, so we take $y\in G_P\cap M$. 
\par 
The element $y$ intertwines $\eta$ {\it with multiplicity one,} in the following sense: there is a non-zero $J^1\cap {J^1}^y$-homomorphism $\phi:\eta \to \eta^y$ and $\phi$ is uniquely determined, up to scalar multiple \cite{18 (5.1.8)}. The element $y$ also intertwines $\mu$, by definition, and it therefore does so with multiplicity one. That is, $\phi$ provides a $J^0\cap {J^0}^y$-homomorphism $\mu\to \mu^y$. Surely $y$ intertwines $\kappa_+$, hence also $\tilde\kappa_+$. That implies $\phi$ to be a homomorphism relative to the group $(J^1\cap N_-\cdot P^\times J^0\cap Q_+)\cap (J^1\cap N_-\cdot P^\times J^0\cap Q_+)^y$. So, $\phi$ is a homomorphism relative to the group generated by 
$$
J^0\cap {J^0}^y \cup \big((J^1\cap N_-\cdot P^\times J^0\cap Q_+)\cap (J^1\cap N_-\cdot P^\times J^0\cap Q_+)^y\big).  
$$ 
This surely contains $(J^0\cap {J^0}^y)P^\times = \bk J\cap \bk J^y$. We have shown that $y$ intertwines $\kappa$, as required. \qed 
\enddemo 
The proof of the proposition yields rather more. First, it produces a map 
$$
\aligned 
f_s:\scr H(\theta_0) &\longrightarrow \scr H(\theta), \\ 
\kappa_0 &\longmapsto \kappa, 
\endaligned 
\tag 3.3.2 
$$ 
using the same notation as in the proof. If $T/F$ is the maximal tamely ramified sub-extension of $P/F$, then $\scr H(\theta_0)$ is a principal homogeneous space over $X_1(T)$ and $\scr H(\theta)$ is a principal homogeneous space over $X_1(T)/X_0(T)_s$.  A simple check yields: 
\proclaim{Corollary 1} 
Let $T/F$ be the maximal tamely ramified sub-extension of $P/F$. The map \rom{(3.3.2)} satisfies 
$$
f_s(\phi\odot\kappa_0) = \phi\odot f_s(\kappa_0), \quad \kappa_0\in \scr H(\theta_0),\ \phi\in X_1(T). 
$$
It is surjective, and its fibres are the orbits of $X_0(T)_s$ in $\scr H(\theta_0)$. 
\endproclaim 
The intertwining analysis in the proof further implies: 
\proclaim{Corollary 2} 
Distinct elements of $\scr H(\theta)$ do not intertwine in $G$. 
\endproclaim 
\subhead 
3.4 
\endsubhead 
We recall some fundamental results from \cite{18}, re-arranged to suit our present purposes. Let $\bk J$ be an open, compact modulo centre, subgroup of $G$ and $\vL$ an irreducible smooth representation of $\bk J$. 
\proclaim{Definition} 
The pair $(\bk J,\vL)$ is an {\rm extended maximal simple type\/} in $G$ if there is an m-simple character $\theta$ in $G$ such that 
\roster 
\item $\bk J = \bk J_\theta$, 
\item the pair $(J^0_\theta,\vL|_{J^0_\theta})$ is a maximal simple type in $G$, and   
\item $\vL|_{H^1_\theta}$ is a multiple of $\theta$. 
\endroster 
\endproclaim 
For the phrase ``maximal simple type'', see 5.5.10 and 6.2 of \cite{18} or 3.6 below. 
\remark{Remark} 
Let $\theta$ be an m-simple character in $\GL nF$. An irreducible representation $\tau$ of $\bk J_\theta$ containing $\theta$ is an extended maximal simple type over $\theta$ if and only if $I_G(\tau) = \bk J_\theta$. 
\endremark 
If $\theta$ is an m-simple character in $G$, we denote by $\scr T(\theta)$ the set of equivalence classes of extended maximal simple types $(\bk J_\theta,\vL)$ with $\vL|_{H^1_\theta}$ a multiple of $\theta$. We recall \cite{18 (6.1.2)} that distinct elements of $\scr T(\theta)$ do not intertwine in $G$. 
\par
Let $P/F$ be a parameter field for $\theta$ and $T/F$ the tame parameter field it contains. Let $\vL\in \scr T(\theta)$. If $\phi\in X_1(T)$ we define $\phi^\bk J\in X_1(\theta)$ as in 3.1, and form the representation 
$$
\phi\odot\vL = \phi^\bk J\otimes\vL. 
\tag 3.4.1 
$$ 
\proclaim{Proposition} 
Let $\theta$ be an m-simple character in $G = \GL nF$, and let $T/F$ be a tame parameter field for $\theta$. 
\roster 
\item 
If $\phi\in X_1(T)$ and $\vL\in \scr T(\theta)$, the representation $\phi\odot\vL$ of \rom{(3.4.1)} lies in $\scr T(\theta)$. The canonical pairing 
$$
\aligned 
X_1(T)\times \scr T(\theta) &\longrightarrow \scr T(\theta), \\ 
(\phi,\vL) &\longmapsto \phi\odot\vL. 
\endaligned 
\tag 3.4.2 
$$ 
defines a structure of $X_1(T)$-space on the set $\scr T(\theta)$. 
\item 
Suppose that $\deg\text{\it cl\/}(\theta) = n$. The sets $\scr H(\theta)$, $\scr T(\theta)$ are then identical and the action \rom{(3.4.2)} endows $\scr T(\theta)$ with the structure of principal homoge\-neous space over $X_1(T)$. 
\endroster 
\endproclaim 
\demo{Proof} 
The character $\phi^\bk J$ in (3.4.1) is trivial on $J^1_\theta$ while its restriction to $G_P\cap \bk J_\theta$ is of the form $\chi\circ\det_P$, for some $\chi\in X_1(P)$ (as in the proof of 3.1 Proposition). Part (1) then follows from the definition of maximal simple type and 3.2 Remark. In part (2), the first assertion follows from the definition of maximal simple type and the second is a case of 3.2 Corollary. \qed 
\enddemo 
\subhead 
3.5 
\endsubhead 
Consider the case where the m-simple character $\theta$ is {\it trivial.} That is, $\theta$ is the trivial character of $U^1_\frak m$, where $\frak m$ is a {\it maximal\/} $\frak o_F$-order in $A = \M nF$. By definition \cite{18 (5.5.10)}, the set $\scr T(\theta)$ consists of equivalence classes of representations $\vL$ of $\bk J_\theta = F^\times U_\frak m$, trivial on $U^1_\frak m$, such that $\vL|_{U_\frak m}$ is the inflation of an irreducible {\it cuspidal\/} representation of $U_\frak m/U^1_\frak m\cong \GL n{\Bbbk_F}$. 
\par 
The set $\scr T(\theta)$ carries a natural action of the group $X_1(F)$, denoted $(\chi,\vL)\mapsto \chi\vL$, where $\chi\vL:x\mapsto \chi(\det x)\vL(x)$, $x\in \bk J_\theta$. Since $\theta$ has parameter field $F$, this action is the same as that given by (3.4.2). 
\par 
Any two choices of $\theta$ (that is, of $\frak m$) are $G$-conjugate. Conjugation by an element $g$ of $G$ induces an $X_1(F)$-bijection $\scr T(\theta) \to \scr T(\theta^g)$, so the $X_1(F)$-set $\scr T(\theta)$ essentially depends only on the dimension $n$. We shall therefore write $\scr T(\theta) = \scr T_n(\bk 0_F)$ in this case. 
\par 
To describe the elements of $\scr T_n(\bk 0_F)$ concretely, it is therefore enough to treat the standard example $\frak m = \M n{\frak o_F}$. We take an unramified field extension $E/F$ of degree $n$, and set $\vD = \Gal EF$. The group $\vD$ acts on $X_1(E)$. We denote by $X_1(E)^{\vD\text{\rm -reg}}$ the set of $\vD$-regular elements of $X_1(E)$ and by $X_1(E)^\vD$ the set of $\vD$-fixed points. We identify $E$ with a subfield of $A$, via some $F$-embedding, such that $E^\times \subset \bk J_\theta = F^\times U_\frak m$. Any two such embeddings are conjugate under $U_\frak m$, so the choice will be irrelevant. 
\proclaim{Proposition} 
Let $\phi\in X_1(E)^{\vD\text{\rm -reg}}$. 
\roster 
\item 
There exists a unique representation $\lambda_\phi\in \scr T_n(\bk 0_F)$ such that 
$$
\roman{tr}\,\lambda_\phi(z\zeta) = (-1)^{n-1} \phi(z) \sum_{\delta\in \vD} \phi^\delta(\zeta), 
\tag 3.5.1 
$$ 
for $z\in F^\times$ and every $\vD$-regular element $\zeta$ of $\bs\mu_E$. 
\item 
The representation $\lambda_\phi$ depends, up to equivalence, only on the $\vD$-orbit of $\phi$ and the map $\phi\mapsto \lambda_\phi$ induces a canonical bijection 
$$ 
\vD\backslash X_1(E)^{\vD\text{\rm -reg}} @>{\ \ \approx\ \ }>> \scr T_n(\bk 0_F). 
\tag 3.5.2 
$$
\item 
If $\chi\in X_1(F)$, write $\chi_E = \chi\circ\N EF \in X_1(E)$. The map $\chi\mapsto \chi_E$ is an isomorphism $X_1(F)/X_0(F)_n \to X_1(E)^\vD$ and, if $\phi\in X_1(E)^{\vD\text{\rm -reg}}$, then 
$$ 
\lambda_{\chi_E\phi} \cong \chi\lambda_\phi,\quad \chi\in X_1(F). 
\tag 3.5.3 
$$
\endroster 
\endproclaim 
\demo{Proof} 
Parts (1) and (2) re-state, for example, \cite{10 (2.2 Proposition)}. Part (3) is elementary. \qed 
\enddemo 
For a more detailed examination of this classical construction and further references, see \cite{13 \S2}. 
\subhead 
3.6 
\endsubhead 
We return to a {\it non-trivial\/} m-simple character $\theta$ in $G = \GL nF$ with parameter field $P/F$ of degree $n/s = \deg \text{\it cl\/}(\theta)$. We let $T/F$ be the tame parameter field inside $P$. We re-interpret the definition of extended maximal simple type by defining a pairing 
$$
\scr T_s(\bk 0_T) \times \scr H(\theta) \longrightarrow \scr T(\theta). 
\tag 3.6.1 
$$ 
To do this, we take $\lambda \in \scr T_s(\bk 0_T)$. Let $P_s/P$ be unramified of degree $s$ and choose a $P$-embedding of $P_s$ in $A$ such that $P_s^\times \subset \bk J_\theta$. Let $T_s/F$ be the maximal tamely ramified sub-extension of $P_s/F$. Thus $T_s/T$ is unramified of degree $s$, and we may identify $\vD = \Gal{P_s}P$ with $\Gal{T_s}T$. There is then a character $\xi\in X_1(T)^{\vD \text{\rm -reg}}$ such that $\lambda \cong \lambda_\xi$. 
\par
We next set $\xi_P = \xi\circ\N{P_s}{T_s} \in X_1(P_s)^{\vD \text{\rm -reg}}$. We so obtain a representation $\lambda_{\xi_P} \in \scr T_s(\bk 0_P)$. Since $\frak a_P \cong \M s{\frak o_P}$, we may think of $\lambda_{\xi_P}$ as a representation of $P^\times U_{\frak a_P}$ trivial on $U^1_{\frak a_P}$. As $\bk J_\theta = P^\times U_{\frak a_P}J^1_\theta$ and $P^\times U_{\frak a_P} \cap J^1_\theta = U^1_{\frak a_P}$, we may extend $\lambda_{\xi_P}$ to a representation of $\bk J_\theta$ trivial on $J^1_\theta$. We denote this representation $\lambda^\bk J_\xi$. 
\par 
We may specify the representation $\lambda_\xi^\bk J$ without reference to the parameter field $P/F$ as follows. We observe that $P^\times J^1_\theta$ is the inverse image, in $\bk J_\theta$, of the centre of the group $\bk J_\theta/J^1_\theta$. 
\proclaim{Lemma} 
Let $T_s/T$ be unramified of degree $s$, such that $T_s^\times\subset \bk J_\theta$. Let $\vD = \Gal{T_s} T$. The representation $\lambda_\xi^\bk J$ has the following properties: 
\roster 
\item"\rm (a)" $\lambda_\xi^\bk J$ is trivial on $J^1_\theta$; 
\item"\rm (b)" the restriction of $\lambda^\bk J_\xi$ to $P^\times J^1_\theta$ is a multiple of the character $(\xi|_{T^\times})^\bk J\in X_1(\theta)$ (cf\. \rom{(3.1.1)}); 
\item"\rm (c)" 
if $\zeta\in \bs\mu_{T_s}$ is $\vD$-regular, then 
$$
\roman{tr}\,\lambda_\xi^\bk J(\zeta) = (-1)^{s-1}\sum_{\delta\in \vD} \xi^\delta(\zeta^{p^r}), 
$$
where $p^r = [P{:}T] = \deg\text{\it cl\/}(\theta)/[T{:}F]$. 
\endroster 
\endproclaim 
We now form the representation 
$$
\lambda_\xi \ltimes \kappa = \lambda_\xi^\bk J \otimes \kappa
\tag 3.6.2 
$$
of $\bk J_\theta$. Re-assembling all the definitions and recalling 3.2 Proposition, we find: 
\proclaim{Proposition} 
Let $\kappa\in \scr H(\theta)$, $\lambda\in \scr T_s(\bk 0_T)$. The representation $\lambda \ltimes \kappa$
of $\bk J_\theta$ lies in $\scr T(\theta)$. 
\roster 
\item 
If $\phi\in X_1(T)$, $\kappa\in \scr H(\theta)$ and $\lambda\in \scr T_s(\bk 0_T)$ then 
$$
\phi\odot(\lambda \ltimes \kappa) \cong (\phi\,\lambda)\ltimes \kappa \cong \lambda\ltimes (\phi\odot\kappa). 
$$
\item 
For any $\kappa\in \scr H(\theta)$, the map 
$$
\align 
\scr T_s(\bk 0_T) &\longrightarrow \scr T(\theta), \\ 
\lambda &\longmapsto \lambda \ltimes \kappa, 
\endalign 
$$ 
is a bijection. 
\endroster 
\endproclaim 
\demo{Proof} 
The first assertion, and the surjectivity of the map in part (2), are the definition of maximal simple type. 
Part (1) is also immediate from the definitions, while part (2) follows from part (1) and 3.2 Proposition. \qed 
\enddemo 
\remark{Remark} 
Observe the similarity between part (2) of the proposition and 1.6 Lemma. 
\endremark 
\head\Rm 
4. Cuspidal representations 
\endhead 
Let $\theta$ be an m-simple character in $G = \GL nF$, with tame parameter field $T/F$. In \S3, we described a canonical action of $X_1(T)$ on the set $\scr T(\theta)$ of equivalence classes of extended maximal simple types in $G$ defined by $\theta$. We wish to translate this into a structure on the set of equivalence classes of irreducible cuspidal representations of $G$ which contain $\theta$. To do this, we must work with the endo-class $\vT = \text{\it cl\/}(\theta)$ and a tame parameter field $E/F$ for $\vT$, taking appropriate care with identifications. 
\subhead 
4.1 
\endsubhead 
If $n\ge1$ is an integer, let $\Ao nF$ denote the set of equivalence classes of irreducible {\it cuspidal\/} representations of $\GL nF$. We summarize the basic results, from 6.2 and 8.4 of \cite{18}, concerning the classification of the elements of $\Ao nF$ via simple types. 
\par 
A representation $\pi\in \Ao nF$ contains an m-simple character $\theta$ in $G = \GL nF$ and, up to $G$-conjugacy, only one. The endo-class $\text{\it cl\/}(\theta)$ therefore depends only on the equivalence class of $\pi$. We accordingly write 
$$
\vt(\pi) = \text{\it cl\/}(\theta) .  
\tag 4.1.1 
$$ 
\indent 
If $\vT\in \scr E(F)$ and if $m\ge1$ is an integer, we put $n = m\deg\vT$ and define 
$$
\Aa mF\vT = \{\pi\in \Ao nF: \vt(\pi) = \vT\}. 
\tag 4.1.2 
$$ 
If we choose an m-realization $\theta$ of $\vT$ in $G = \GL nF$, a representation $\pi\in \Ao nF$ lies in $\Aa mF\vT$ if and only if $\pi$ contains $\theta$ ({\it cf\.} 2.5 Lemma 2). If $\vL \in \scr T(\theta)$, the representation $\pi _\vL = \cind_{\bk J_\theta}^G\,\vL$ is irreducible, cuspidal and contains $\theta$: that is, $\pi_\vL \in \Aa mF\vT$. The map 
$$
\aligned 
\scr T(\theta) &\longrightarrow \Aa mF\vT, \\ 
\vL &\longmapsto \cind_{\bk J_\theta}^G\vL, 
\endaligned 
\tag 4.1.3  
$$ 
is a bijection. 
\remark{Remark} 
Let $\pi\in \Ao nF$ contain a simple character $\vf$. One may show that $\vf$ is m-simple and $\text{\it cl\/} (\vf) = \vt(\pi)$. We have no use here for this fact, but a proof may be found in \cite{15}. 
\endremark 
\subhead 
4.2 
\endsubhead 
Let $\vT\in \scr E(F)$ have degree $d$, let $E/F$ be a tame parameter field for $\vT$, and let $\varPsi$ be  a totally wild $E/F$-lift of $\vT$. 
\par 
Let $m\ge1$ be an integer, and set $n = md$. We take a pair $(\theta,\bsi)$ consisting of an m-realization $\theta$ of $\vT$ in $G = \GL nF$ and an $F$-embedding $\bsi$ of $E$ in $\M nF$ such that $\bsi E/F$ is a tame parameter field for $\theta$ satisfying $\bsi_*\varPsi = \text{\it cl\/}(\theta_{\bsi E})$. Such a pair $(\theta,\bsi)$ exists and is uniquely determined, up to conjugation in $G$ by an element commuting with $\bsi E$ (2.7 Proposition). We first use $\bsi$ to translate the natural action (3.4.2) of $X_1(\bsi E)$ on $\scr T(\theta)$ into an action of $X_1(E)$. The induction relation  (4.1.3) then turns this into an action of $X_1(E)$ on $\Aa mF\vT$. In all, we get a {\it canonical\/} pairing 
$$
\aligned 
X_1(E) \times \Aa mF\vT &\longrightarrow \Aa mF\vT, \\
(\phi,\pi) &\longmapsto \phi\odot_\varPsi\pi, 
\endaligned 
\tag 4.2.1 
$$
satisfying $(\phi\phi')\odot_\varPsi \pi = \phi\odot_\varPsi(\phi'\odot_\varPsi\pi)$. This action of $X_1(E)$ does depend on $\varPsi$, but on no other choices. 
\par 
In the case $E=F$, we have $\vT=\varPsi$ and the $\odot_\varPsi$-action of $X_1(F)$ on $\Aa mF\vT$ is the standard one $(\phi,\pi)\mapsto \phi\pi$, where $\phi\pi: x\mapsto \phi(\det x)\,\pi(x)$. 
\par
Another special case is worthy of note. As a direct consequence of 3.4 Proposition, we have: 
\proclaim{Proposition} 
The action $\odot_\varPsi$ endows $\Aa 1F\vT$ with the structure of a principal homogeneous space over $X_1(E)$. 
\endproclaim 
\remark{Remark} 
In the general situation, let $\varPsi'$ be some other totally wild $E/F$-lift of $\vT$. Thus there exists $\gamma\in \roman{Aut}(E|F)$ such that $\varPsi' = \varPsi^\gamma$ (2.4 Corollary). We then get the relation 
$$
\phi^\gamma\odot_{\varPsi^\gamma}\pi = \phi\odot_\varPsi\pi, \quad \phi\in X_1(E),\ \pi \in \Aa mF\vT. 
\tag 4.2.2 
$$ 
Indeed, this relation holds for any $F$-isomorphism $\gamma:E\to E^\gamma$. (Observe the apparent parallel with the discussion in 1.5, especially (1.5.3).) 
\endremark 
\head\Rm 
5. Algebraic induction maps 
\endhead 
In \S1, we saw that the irreducible representations of the Weil group $\scr W_F$ can be constructed in a uniform manner, using only very particular induction processes (as in 1.6 Proposition). In this section, we define analogous ``algebraic induction maps'' on the other side. 
\subhead 
5.1 
\endsubhead 
To start, we recall something of the Glauberman correspondence \cite{21}, as developed in Appendix 1 to \cite{4}. In this sub-section, we use a system of notation distinct from that of the rest of the paper. 
\par 
Let $G$ be a finite group, and let $\roman{Irr}\,G$ denote the set of equivalence classes of irreducible representations of $G$. Let $A$ be a finite soluble group and $A\to \roman{Aut}\,G$ a homomorphism with image of order relatively prime to that of $G$. Let $G^A$ denote the group of $A$-fixed points in $G$. The group $A$ acts on the set $\roman{Irr}\,G$, so let $\roman{Irr}^A\,G$ denote the set of fixed points. The Glauberman correspondence is a canonical bijection 
$$
\eurm g_G^A: \roman{Irr}^A\,G \longrightarrow \roman{Irr}\,G^A, 
$$ 
which depends only on the image of $A$ in $\roman{Aut}\,G$, and is transitive in the following sense. Let $B$ be a normal subgroup of $A$. Thus $A/B$ acts on $G^B$, and similarly on representations. We then have 
$$
\eurm g_G^A = \eurm g_{G^B}^{A/B}\circ \eurm g_G^B. 
\tag 5.1.1 
$$
To describe the correspondence $\eurm g_G^A$, it is therefore enough to consider the case where $A$ is {\it cyclic.} 
\proclaim\nofrills{(5.1.2)} \ 
Suppose $A$ is cyclic, of order relatively prime to $|G|$. Let $\rho\in \roman{Irr}^A\,G$. 
\roster 
\item  There is a unique representation $\tilde\rho$ of $A\ltimes G$ such that $\tilde\rho|_G\cong \rho$ and $\det\tilde\rho|_A = 1$. 
\item 
There is a unique $\rho^A \in \roman{Irr}\,G^A$ such that 
$$
\roman{tr}\,\rho^A(h) = \eps\,\roman{tr}\,\tilde\rho(ah), 
\tag 5.1.3 
$$
for all $h\in G^A$, any generator $a$ of $A$, and a constant $\eps$. 
\item 
The constant $\eps$ has the value $\pm1$. 
\endroster 
\endproclaim 
In this context, the representation $\rho^A$ is $\eurm g_G^A(\rho)$. The constant $\eps$ depends on both $\rho$ and the image of $A$ in $\roman{Aut}\,G$: different cyclic operator groups on $G$ may have the same fixed points and give the same character correspondence, but for different constants $\eps$. An instance of this occurs in 5.4 below. 
\subhead 
5.2 
\endsubhead 
We return to our standard situation and notation. 
\par 
A {\it complementary subgroup\/} of $F^\times$ is a closed subgroup $C$ of $F^\times$ such that the product map $C\times U^1_F \to F^\times$ is bijective. A complementary subgroup $C$ is necessarily of the form $C = C_F(\vp_F) = \ags{\vp_F,\bs\mu_F}$, for some prime element $\vp_F$ of $F$. If $\vp_F$, $\vp_F'$ are prime elements of $F$, then $C_F(\vp'_F) = C_F(\vp_F)$ if and only if $\vp'_F = \alpha\vp_F$, for some $\alpha\in \bs\mu_F$. 
\proclaim{Lemma} 
Let $E/F$ be a finite, tamely ramified field extension, and let $\vp_F$ be a prime element of $F$. 
\roster 
\item 
There exists a unique complementary subgroup $C_E(\vp_F)$ of $E$ such that $\vp_F\in C_E(\vp_F)$. 
\item  
Let $F\subset L\subset E$. The group $C_E(\vp_F)\cap L^\times$ is the complementary subgroup $C_L(\vp_F)$ of $L^\times$ containing $\vp_F$. 
\endroster 
\endproclaim  
\demo{Proof} 
Let $\vp$ be a prime element of $E$, and set $e = e(E|F)$. Thus $\vp^e = \zeta\vp_Fu$, for some $\zeta\in \bs\mu_E$ and $u\in U^1_E$. Since $p$ does not divide $e$, there exists a unique element $v$ of $U^1_E$ such that $v^e = u$. Replacing $\vp$ by $v^{-1}\vp$, we may assume $\vp^e\in \bs\mu_E\vp_F$. The group $\ags{\vp,\bs\mu_E}$ is a complementary subgroup of $E^\times$ containing $\vp_F$. It is clearly the unique such subgroup with this property. Assertion (2) is now immediate. \qed 
\enddemo 
\subhead 
5.3 
\endsubhead 
We define a family of ``algebraic induction maps''. The construction proceeds in two steps, reflecting the discussion in 5.1.  For the first, we work in the following special situation. 
\proclaim{Notation} 
\roster 
\item 
$\theta$ is a simple character in $G = \GL nF$, 
\item 
$\theta\in \scr C(\frak a,\beta)$, for a simple stratum $[\frak a,\beta]$ in $A = \M nF$ such that $P = F[\beta]/F$ is \rom{totally ramified of degree} $n$; 
\item $E/F$ is the maximal tamely ramified sub-extension of $P/F$; 
\item $[E{:}F] = e$ and $[P{:}E] = p^r$. 
\endroster 
\endproclaim 
The simple character $\theta$ is m-simple, because of condition (2). Consequently, the set $\scr T(\theta) = \scr H(\theta)$ is a principal homogeneous space over $X_1(E)$ (3.4 Proposition). 
\par 
Let $\eta$ be the 1-Heisenberg representation of $J^1 = J^1(\beta,\frak a)$ over $\theta$. The natural conjugation action of $E^\times$ on $J^1$ induces an action of the group $M_{E/F} = E^\times/F^\times U^1_E$ on the finite $p$-group $J^1/\roman{Ker}\,\theta$. The group $M_{E/F}$ is cyclic, of order $e = [E{:}F]$ (which is relatively prime to $p$). The group of $M_{E/F}$-fixed points in $J^1/\roman{Ker}\,\theta$ is $J^1_E/\roman{Ker}\,\theta_E$ ({\it cf\.} 4.1 Lemma 1 of \cite{4}). 
\par 
Taking complementary subgroups $C_F(\vp_F) \subset C_E(\vp_F)$ as in 5.2, the canonical map $C_E(\vp_F) \to M_{E/F}$ induces an isomorphism $C_E(\vp_F)/C_F(\vp_F) \cong M_{E/F}$. We may therefore switch, as convenient, between $C_E(\vp_F)$ and $M_{E/F}$ as operator groups. 
\par 
As in 2.3, $G_E$ denotes the $G$-centralizer of $E^\times$. Using the standard notation of that paragraph, let $\eta_E$ denote the 1-Heisenberg representation of $J^1_E = J^1\cap G_E = J^1(\beta,\frak a_E)$ over $\theta_E = \theta|_{H^1_E}$. For the next result, we initially view $\eta$ as a representation of the finite $p$-group $J^1/\roman{Ker}\, \theta$, and similarly for $\eta_E$. 
\proclaim{Lemma} 
\roster 
\item 
There exists a unique representation $\tilde\eta$ of $M_{E/F}\ltimes J^1/\roman{Ker}\,\theta$ extending $\eta$ and such that $\det\tilde\eta|_{M_{E/F}} = 1$. 
\item 
There is a constant $\eps = \eps_{E/F} = \pm1$ such that 
$$
\roman{tr}\,\tilde\eta(zj) = \eps\,\roman{tr}\,\eta_E(j), 
\tag 5.3.1
$$
for every $j\in J^1_E$ and every generator $z$ of the cyclic group $M_{E/F}$. 
\item 
The sign $\eps$ depends only on on $\theta$, and not on the choice of tame parameter field $E/F$. 
\endroster 
\endproclaim 
\demo{Proof} 
The first assertion is elementary ({\it cf\.} (5.1.2)). In (2), the Glauberman correspondence gives a unique irreducible representation $\zeta$ of $J^1_E/\roman{Ker}\,\theta_E$ such that $\roman{tr}\,\zeta(j) = \eps\,\roman{tr}\,\tilde\eta(zj)$, for $z$ and $j$ as above. We inflate $\tilde\eta$, $\zeta$ to representations of $C_E(\vp_F)\ltimes J^1$, $J^1_E$ respectively. Replacing $j$ by $jh$, $h\in H^1_E$, we find that $\roman{tr}\,\zeta(jh) =\theta_E(h)\,\roman{tr}\,\zeta(j)$. Thus $\zeta$ contains $\theta_E$, whence $\zeta\cong \eta_E$, as required. 
\par 
The sign $\eps = \eps_{E/F}$ can be computed using \cite{2 (8.6.1)}. We form the finite, elementary abelian $p$-group $V = V_\theta = J^1/H^1$. This carries a nondegenerate alternating form over $\Bbb F_p$, as in \cite{18 (5.1)}. The conjugation action of $E^\times$ on $J^1$ induces an action of $M_{E/F}$ on $V$, and $V$ provides a symplectic representation of $M_{E/F}$ over $\Bbb F_p$. In the notation of \cite{13, \S3}, $\eps_{E/F}$ is given by 
$$
\eps_{E/F} = t_{M_{E/F}}(V). 
\tag 5.3.2 
$$ 
That is, $\eps_{E/F}$ is an invariant of the $M_{E/F}$-module $V$ ({\it cf\.} \cite{13 (3.3)}), and $V$ is determined directly from $[\frak a,\beta]$ (or from $\theta$). Any two choices of $E$ are $J^1$-conjugate (2.6 Proposition), and such a conjugation has no effect on the module structure of $V$. \qed 
\enddemo 
\remark{Remark} 
The relation (5.3.2) shows that $\eps_{E/F}$ may be computed purely in terms of a simple stratum underlying $\theta$ (and does not depend on the choice of this stratum). 
\endremark 
\proclaim{Proposition} 
Define $\eps = \eps_{E/F} = \pm1$ by \rom{(5.3.1).} Let $\vL\in \scr T(\theta)$. There exists a unique $\vL_E\in \scr T(\theta_E)$ such that 
$$
\roman{tr}\,\vL_E(x) = \eps_{E/F}\,\roman{tr}\,\vL(x), 
\tag 5.3.3 
$$ 
for all $x\in \bk J_E$ such that $\ups_E(\det_E x) = \ups_F(\det x)$ is relatively prime to $[E{:}F]$. The map 
$$
\aligned 
\scr T(\theta) &\longrightarrow \scr T(\theta_E), \\ 
\vL&\longmapsto \vL_E, 
\endaligned 
\tag 5.3.4
$$
is a bijection, and an isomorphism of $X_1(E)$-spaces. 
\endproclaim 
\demo{Proof} 
If $\vL \in \scr T(\theta)$, a representation $\vL_E$ of $\bk J_E$ satisfying (5.3.3) must lie in $\scr T(\theta_E)$, and is uniquely determined by $\vL$. Moreover, if we have such a pair and if $\phi\in X_1(E)$, then (in the notation of 3.2)
$$
\roman{tr}\,\phi\odot\vL_E(x) = \eps_{E/F} \,\roman{tr}\,\phi\odot\vL(x), 
$$
for all $x$ as before. In this situation, each of the sets $\scr T(\theta)$, $\scr T(\theta_E)$ is a principal homogeneous space over $X_1(E)$. The result will follow, therefore, if we produce one representation $\vL\in \scr T(\theta)$ for which there exists $\vL_E$ satisfying (5.3.3). 
\par 
To do this, we choose a prime element $\vp_F$ of $F$ and let $C_E = C_E(\vp_F)$ be the unique complementary subgroup of $E^\times$ containing $\vp_F$. Let $\pre p{\bk J} = \pre p{\bk J}(\vp_F)$ denote the inverse image in $\bk J$ of the unique Sylow pro-$p$ subgroup of the profinite group $\bk J/\ags{\vp_F}$. Thus $\bk J = C_E \cdot \pre p{\bk J}$. Also, $\pre p{\bk J_E} = \pre p{\bk J}\cap G_E$ is the inverse image in $\bk J_E$ of the Sylow pro-$p$ subgroup of $\bk J_E/\ags{\vp_F}$. 
\par
We choose $\vL\in \scr T(\theta)$ such that $\vL(\vp_F) = 1$ and $\det\vL$ is trivial on $C_E$. Since $\dim\vL$ is a power of $p$ and $\bs\mu_F$ is central, these conditions imply that $\vL$ is trivial on $\bs\mu_F$, hence also on $C_F = C_E\cap F^\times$. We may therefore view $\vL$ as the inflation of a representation of $C_E/C_F\ltimes \pre p{\bk J}/\ags{\vp_F}$. The Glauberman correspondence gives a representation $\tilde\vL_E$ of $\pre p{\bk J_E}/\ags{\vp_F}$ such that 
$$
\roman{tr}\,\tilde\vL_E(x) = \eps\,\roman{tr}\,\vL(cx),\quad x\in \pre p{\bk J_E}/\ags{\vp_F}, 
$$
for any generator $c$ of the cyclic group $C_E/C_F$. Moreover, $\tilde\vL_E|_{J^1_E} \cong \eta_E$. We extend $\tilde\vL_E$, by triviality, to a representation of $C_E/C_F \times \pre p{\bk J_E}/\ags{\vp_F}$ and then inflate it to a representation $\vL_E$ of $C_E\cdot\pre p{\bk J_E} = \bk J_E$. This representation satisfies (5.3.3). \qed 
\enddemo 
We interpret the proposition in terms of cuspidal representations of linear groups, following the discussion in \S4. 
\proclaim{Corollary 1} 
Let $\vT = \text{\it cl\/}(\theta)$, $\vT_E = \text{\it cl\/}(\theta_E)$. The map \rom{(5.3.4)} induces a canonical bijection 
$$
\aligned 
\Aa1F\vT &\longrightarrow \Aa1E{\vT_E}, \\ 
\pi &\longmapsto \pi_E, 
\endaligned 
\tag 5.3.5 
$$
such that 
$$
(\phi\odot_{\vT_E}\pi)_E = \phi\pi_E, 
\tag 5.3.6 
$$ 
for $\phi\in X_1(E)$ and $\pi\in \Aa1F\vT$. 
\endproclaim 
\demo{Proof} 
We start from the endo-classes $\vT_E\in \scr E(E)$ and $\vT\in \scr E(F)$. Thus $\vT_E$ is totally wild, $\vT = \frak i_{E/F}\vT_E$ ({\it cf\.} (2.3.4)), and $E/F$ is a tame parameter field for $\vT$. 
\par
We identify $E$ with a subfield of $\M nF$. Following 2.7 Proposition, we choose a realization $\theta$ of $\vT$ in $G = \GL nF$ such that $E/F$ is a tame parameter field for $\theta$ and $\text{\it cl\/}(\theta_E) = \vT_E$. Thus $\theta$ is uniquely determined up to conjugation by an element of $G_E$. 
\par 
Let $\pi\in \Aa1F\vT$. By definition, $\pi$ contains a representation $(\bk J_\theta,\vL) \in \scr T(\theta)$, giving $\cind_{\bk J_\theta}^G\,\vL \cong \pi$. Using the map (5.3.4), we form the representation $\vL_E\in \scr T(\theta_E)$ and set $\pi_E = \cind_{\bk J_{\theta_E}}^{G_E} \vL_E$. The $G_E$-conjugacy class of $\vL_E$ is then independent of the choice of $\theta$, so $\pi\mapsto \pi_E$ gives a canonical map $\Aa1F\vT\to \Aa1E{\vT_E}$. The bijectivity of this map and property (5.3.6) then follow from the proposition. \qed 
\enddemo 
We let 
$$ 
\ind_{E/F}:\Aa1E{\vT_E} \longrightarrow \Aa1F\vT 
\tag 5.3.7 
$$
denote the inverse of the map (5.3.5). Thus $\ind_{E/F}$ is a bijection satisfying 
$$
\ind_{E/F}\,\phi\rho = \phi\odot_{\vT_E} \ind_{E/F}\,\rho, \quad \rho \in \Aa1E{\vT_E},\ \phi\in X_1(E). 
\tag 5.3.8 
$$ 
It also has an important naturality property. 
\proclaim{Corollary 2} 
If $\gamma:E\to E^\gamma$ is an isomorphism of local fields, then 
$$
\ind_{E^\gamma/F^\gamma}\,\rho^\gamma = \big(\ind_{E/F}\,\rho\big)^\gamma, 
$$ 
for all $\rho\in \Aa1E{\vT_E}$. 
\endproclaim 
\demo{Proof} 
To construct the representation $\ind_{E/F}\,\rho$ more directly, we first choose an $F$-embedding of $E$ in $\M nF$. We then choose a realization $\vf$ of $\vT_E$ in $G_E$, and form the simple character $\vf^F$ in $G$, as in (2.3.3). Let $\lambda\in \scr T(\vf)$ occur in $\rho$. There is a unique $\lambda^F\in \scr T(\vf^F)$ such that $(\lambda^F)_E \cong \lambda$. The representation $\ind_{E/F}\,\rho$ is then equivalent to $\cind\,\lambda^F$. Changing the embedding of $E$ in $\M nF$ or the realization $\vf$ only replaces $\lambda^F$ by a $G$-conjugate. It follows that $\big(\ind_{E^\gamma/F^\gamma}\,\rho^\gamma\big)^{\gamma^{-1}} \cong \ind_{E/F}\,\rho$, as required. \qed 
\enddemo 
\subhead 
5.4 
\endsubhead 
For the second step in the construction, we pass to the general situation. The notation that follows will be standard for the rest of the section. 
\proclaim{Notation} 
\roster 
\item 
$[\frak a,\beta]$ is an m-simple stratum in $A = \M nF$ and $\theta\in \scr C(\frak a,\beta,\psi_F)$; 
\item 
$P_0 = F[\beta]$ and $m = n/[P_0{:}F]$; 
\item 
$E_0/F$ is the maximal tamely ramified sub-extension of $P_0/F$ and $K_0/F$ is its maximal unramified sub-extension; 
\item 
$n_0 = [E_0{:}F]$, $e = e(E_0|F)$ and $p^r = [P_0{:}E_0]$; 
\item 
$P/P_0$ is unramified of degree $m$ and $P^\times \subset \bk J(\beta,\frak a) = \bk J_\theta$; 
\item 
$K/F$ is the maximal unramified sub-extension of $P/F$ and $E = E_0K/F$ its maximal tamely ramified sub-extension; 
\item $\vG = \Gal KF$ and $\vD = \Gal P{P_0} \cong \Gal E{E_0} \cong \Gal K{K_0}$. 
\endroster 
\endproclaim 
As usual, we put $G = \GL nF$. Any two extensions $P/P_0$, satisfying the conditions (5), are $J^0_\theta$-conjugate. As we will see (Comments, 5.7), this means that the choice of $P/P_0$ is irrelevant.  We use the notational scheme of 2.3. In particular, $G_K$ is the $G$-centralizer of $K^\times$ and $\theta_K \in \scr C(\frak a_K,\beta,\psi_K)$. The extension $E/K$ is thus a tame parameter field for $\theta_K$. 
\par 
Our immediate aim is to define a canonical $X_1(E_0)$-bijection 
$$
\ell_{K/F}: \scr H(\theta) @>{\ \ \approx\ \ }>> \scr H(\theta_K)^\vD. 
$$ 
This will be achieved in 5.6 Proposition. 
\par 
We apply the Glauberman correspondence to the action of $\bs\mu_K$ on $J^1= J^1(\beta,\frak a)$. 
\proclaim{Lemma} 
Let $\tilde\eta$ be the representation of $\bs\mu_K\ltimes J^1$ such that $\tilde\eta|_{J^1}\cong \eta$ and $\det\tilde\eta |_{\bs\mu_K} = 1$. There is a constant $\eps^0_{\bs\mu_K} = \pm1$ and a character $\eps^1_{\bs\mu_K}:\bs\mu_K\to \{\pm1\}$ such that 
$$
\roman{tr}\,\eta_K(x) = \eps^0_{\bs\mu_K}\eps^1_{\bs\mu_K}(\zeta)\,\roman{tr}\,\tilde\eta(\zeta x), 
\tag 5.4.1 
$$
for $x\in J^1_K$ and any $\vG$-regular element $\zeta$ of $\bs\mu_K$. 
\endproclaim 
\demo{Proof} 
Just as in 5.3 Lemma, there is a constant $\eps(\bs\mu_K) = \pm1$ such that 
$$
\roman{tr}\,\eta_K(x) = \eps(\bs\mu_K)\,\roman{tr}\,\tilde\eta(\zeta x), 
$$
for $x\in J^1_K$ and any {\it generator\/} $\zeta$ of $\bs\mu_K$. More generally, let $\zeta\in \bs\mu_K$ be $\vG$-regular. The set of $\zeta$-fixed points in $J^1$ is then $J^1_K$ and 
$$
\roman{tr}\,\eta_K(x) = \eps(\ags\zeta)\,\roman{tr}\,\tilde\eta(\zeta x), \quad x\in J^1_K, 
$$ 
where $\eps(\ags\zeta) = \pm1$ depends on the subgroup of $\bs\mu_K$ generated by $\zeta$. According to \cite{2 (8.6.1)} (but using the notation of \cite{13 \S3}), we have $\eps(\ags\zeta) = t_{\ags\zeta}(V)$, where $V = J^1/H^1$. However, \cite{13 (Proposition 3.6)} gives 
$$
\eps(\ags\zeta) = \eps^0_{\bs\mu_K}\,\eps^1_{\bs\mu_K}(\zeta), 
$$ 
where 
$$
\eps^j_{\bs\mu_K} = t^j_{\bs\mu_K}(V),\quad j=0,1, 
\tag 5.4.2  
$$ 
as required. \qed 
\enddemo 
As in 5.3, the relation (5.4.2) shows that the constant $\eps^0_{\bs\mu_K}$ and the character $\eps^1_{\bs\mu_K}$ are determined by a simple stratum underlying $\theta$. 
\remark{Remark} 
The character $\eps^1_{\bs\mu_K}$ is, by its definition in terms of the conjugation action of $\bs\mu_K$ on $V$, trivial on $\bs\mu_F$. Moreover, the definition of $t^1_{\bs\mu_K}(V)$ in \cite{13 (3.4 Definition 3)} shows that $\eps^1_{\bs\mu_K}$ is trivial when $p=2$. 
\endremark 
\subhead 
5.5 
\endsubhead 
We choose a prime element $\vp_F$ of $F$. As in 5.2 Lemma, let $C_{E_0}(\vp_F)$ (resp\. $C_E(\vp_F)$) be the unique complementary subgroup of $E^\times_0$ (resp\. $E^\times$) containing $\vp_F$. 
\par 
We consider the group $P_0^\times J^1$, described more intrinsically as the inverse image in $\bk J$ of the centre of the group $\bk J/J^1$. The group $P_0^\times J^1/\ags{\vp_F}$ is profinite, and has a unique Sylow pro-$p$ subgroup. We denote by $\pre p{\bk J}$ the inverse image of this subgroup in $\bk J$. {\it The definition of $\pre p{\bk J}$ depends on the choice of $\vp_F$:\/} we accordingly use the notation $\pre p{\bk J} = \pre p{\bk J}(\vp_F)$ when it is necessary to emphasize this fact. The quotient $\pre p{\bk J}(\vp_F)/J^1$ is cyclic, generated by an element $\vp_0$ of $P_0$ such that $\vp_0^{p^r} \equiv \vp_F \pmod{U^1_{P_0}}$. 
\par 
Using the quotient $P^\times J^1_K/\ags{\vp_F}$, we may similarly define a subgroup $\pre p{\bk J_K} = \pre p{\bk J_K}(\vp_F)$ of $\bk J_K$. This satisfies $\pre p{\bk J_K} = \pre p{\bk J}\cap \bk J_K$. 
\par 
We will need to apply the Conjugacy Lemma of 2.6 in this new context. We must therefore verify: 
\proclaim{Lemma} 
The pro-$p$ group $\pre p{\bk J}/\ags{\vp_F}$ admits a $P^\times$-stable filtration satisfying the conditions of \rom{(2.6).} 
\endproclaim 
\demo{Proof} 
We have to refine the standard filtration of the group $J^1$. A generator $\vp_0$ of $\pre p{\bk J}/J^1$ acts on each step $J^k/J^{k+1}$ of the standard filtration as a unipotent automorphism commuting with the natural action of $P^\times$. We therefore insert into the standard filtration the extra steps $\roman{Ker}\,(\vp_0{-}1)^r|_{J^k/J^{k+1}}$. \qed 
\enddemo 
\subhead 
5.6 
\endsubhead 
We start by noting a consequence of 3.2 Proposition. 
\proclaim{Lemma 1} 
Let $\kappa_1,\kappa_2\in \scr H(\theta)$. If $\kappa_1|_{P^\times J^1}\cong \kappa_2|_{P^\times J^1}$, then $\kappa_1\cong \kappa_2$. 
\endproclaim 
We apply the Glauberman correspondence to representations of $\pre p{\bk J}$. We first identify a family of special elements of $\scr H(\theta)$. 
\proclaim{Lemma 2} 
Let $\vp_F$ be a prime element of $F$. There exists $\kappa^0\in \scr H(\theta)$ such that 
\roster 
\item $\vp_F\in \roman{Ker}\,\kappa^0$ and 
\item $C_E(\vp_F)\subset \roman{Ker}\,\det\kappa^0$. 
\endroster 
These conditions determine $\kappa^0$ uniquely, up to tensoring with a character $\phi\in X_1(\theta)$, trivial on $J^0$ and such that $\phi^{p^r} = 1$. Moreover, $\bs\mu_F\subset \roman{Ker}\,\kappa^0$. 
\endproclaim 
\demo{Proof} 
The existence of $\kappa^0$ follows from 3.2 Proposition, and its uniqueness property from Lemma 1. The group $\bs\mu_F$ is central in $G$. It has order relatively prime to $p$, while $\dim\kappa^0$ is a power of $p$. The final assertion thus follows from condition (2). \qed 
\enddemo 
\proclaim{Lemma 3} 
Let $\kappa^0\in \scr H(\theta)$ satisfy the conditions of \rom{Lemma 2} relative to the element $\vp_F$. There exists a unique $\kappa^0_K\in \scr H(\theta_K)$ such that 
\roster 
\item $\vp_F\in \roman{Ker}\,\kappa^0_K$, 
\item $C_E(\vp_F)\subset \roman{Ker}\,\det\kappa^0_K$, and 
\item 
if $h\in \pre p{\bk J_K}(\vp_F)$, then 
$$
\roman{tr}\,\kappa^0_K(h) = \eps^0_{\bs\mu_K}\,\eps^1_{\bs\mu_K}(\zeta)\,\roman{tr}\,\kappa^0(\zeta h), 
\tag 5.6.1 
$$
for every $\vG$-regular element $\zeta$ of $\bs\mu_K$. 
\endroster 
\endproclaim 
\demo{Proof} 
The Glauberman correspondence, together with (5.4.1), yields a unique representation $\kappa^0_K$ of $\pre p{\bk J_K}$ satisfying the desired relation (5.6.1) and condition (1). This representation admits extension to a representation of $\bk J_K = P^\times J^1_K$, and the extension may be chosen to satisfy condition (2). In this case, any representation of $\bk J_K$ extending $\eta_K$ lies in $\scr H(\theta_K) = \scr T(\theta_K)$, so the uniqueness follows from Lemma 1. \qed 
\enddemo 
\remark{Remark}
As in Lemma 2, condition (2) of Lemma 3 implies $\bs\mu_K \subset \roman{Ker}\,\kappa^0_K$. In (5.6.1) therefore, we have $\roman{tr}\,\kappa^0_K(\zeta h) = \roman{tr}\,\kappa^0_K(h)$. 
\endremark  
The set $\scr H(\theta_K)$ carries an action of $X_1(E)$, making it a principal homogeneous space over $X_1(E)$, as in 3.2 Corollary. Consequently, it carries an action of $X_1(E_0)$ via the canonical map $X_1(E_0) \to X_1(E)$ given by $\chi\mapsto \chi_E = \chi\circ\N E{E_0}$. If we view $\vD$ as $\Gal E{E_0}$, this canonical map induces an isomorphism $X_1(E_0)/X_0(E_0)_m$ with the group $X_1(E)^\vD$ of $\vD$-fixed points in $X_1(E)$. 
\par 
We consider the group $\vD = \Gal P{P_0}$. For each $\delta\in \vD$, there exists $j_\delta\in G$ such that $j_\delta^{-1}xj_\delta = x^\delta$, $x\in P^\times$. The same relation holds for $x\in K^\times$ or $E^\times$. The definition of $P$ shows we may choose $j_\delta\in U_{\frak a_{P_0}}\subset J^0_\theta$. We use these elements $j_\delta$ to extend the canonical action of $\vD$ on $K$ to one on $G_K$, by setting $g^\delta = j_\delta^{-1}gj_\delta$, $g\in G_K$, $\delta\in \vD$. Since $j_\delta\in J^0$, it conjugates $\theta$ to itself and therefore normalizes $\bk J_K = \bk J_\theta\cap G_K$. Moreover, conjugation by $j_\delta$ fixes $\theta_K$ and so $\vD$ acts on the set $\scr H(\theta_K)$ via $\kappa \mapsto \kappa^\delta = \kappa^{j_\delta}$. If $\phi\in X_1(E)$, we also have 
$$
(\phi\odot\lambda)^\delta \cong \phi^\delta\odot\lambda^\delta, \quad \phi\in X_1(E),\ \lambda\in \scr H(\theta_K), 
\tag 5.6.2 
$$ 
where $\vD$, viewed as $\Gal E{E_0}$, acts on $X_1(E)$ in the natural way. 
\proclaim{Lemma 4} 
\roster 
\item 
The representation $\kappa^0_K$ of \rom{Lemma 3} lies in $\scr H(\theta_K)^\vD$, and 
\item 
$\scr H(\theta_K)^\vD$ is a principal homogeneous space over $X_1(E_0)/X_0(E_0)_m$. 
\endroster 
\endproclaim 
\demo{Proof} 
Since $j_\delta\in J^0$, the function $\roman{tr}\,\kappa^0$ is invariant under conjugation by $j_\delta$. Assertion (1) now follows from the definition of $\kappa^0_K$, and (2) from (5.6.2). \qed 
\enddemo 
\proclaim{Proposition} 
Let $\kappa^0\in \scr H(\theta)$ satisfy the conditions of \rom{Lemma 2,} and define $\kappa^0_K\in \scr H(\theta_K)$ as in \rom{Lemma 3.} 
\roster 
\item 
There is a unique $X_1(E_0)$-map 
$$
\ell_{K/F}:\scr H(\theta) \longrightarrow \scr H(\theta_K) 
\tag 5.6.3 
$$
such that 
$$
\ell_{K/F}(\kappa^0) = \kappa^0_K. 
\tag 5.6.4 
$$ 
\item 
The map $\ell_{K/F}$ is injective, and its image is the set $\scr H(\theta_K)^\vD$ of $\vD$-fixed points in $\scr H(\theta_K)$. 
\item  
The definition of $\ell_{K/F}$ is independent of the choices of $\vp_F$ and the representation $\kappa^0$ satisfying the conditions of \rom{Lemma 2,} relative to $\vp_F$. 
\endroster 
\endproclaim 
\demo{Proof} 
Parts (1) and (2) follow from the fact that $\scr H(\theta)$ and $\scr H(\theta_K)^\vD$ are both principal homogeneous spaces over $X_1(E_0)/X_0(E_0)_m$. 
\par 
We prove part (3). We work at first relative to fixed choices of $\vp_F$ and $\kappa^0$. Taking $\kappa\in \scr H(\theta)$, we write $\kappa = \phi\odot\kappa^0$, for $\phi\in X_1(E_0)$ uniquely determined modulo $X_0(E_0)_m$ ({\it cf\.} 3.2 Corollary). We define 
$$
\ell_{K/F}(\kappa) = \ell_{K/F}(\phi\odot\kappa^0) = \phi_E\odot \kappa^0_K. 
$$ 
This is the unique $X_0(E_0)$-map $\scr H(\theta) \to \scr H(\theta_K)$ with the property (5.6.4). It is clearly injective with image $\scr H(\theta_K)^\vD$. 
\par
We have to check that this definition of $\ell_{K/F}$ is independent of choices. First, if we keep $\vp_F$ fixed, the definition is independent of the choice of $\kappa^0$, as follows from the uniqueness statement in Lemma 2. \par
Next, let us replace $\vp_F$ by $u\vp_F$, for some $u\in U^1_F$. The group $\pre p{\bk J}(\vp_F)$ is then unchanged, but $C_E$ is replaced by $C'_E = \ags{v\vp_{E_0},\bs\mu_K}$, where $v\in U^1_F$ satisfies $v^e=u$. There is a character $\phi\in X_0(E_0)$, of $p$-power order modulo $X_0(E_0)_m$, such that $\phi\odot\kappa^0$ satisfies the conditions of Lemma 2 relative to $\vp_F'$. The element of $\scr H(\theta_K)$ corresponding to $\phi\odot\kappa^0$ via Lemma 3 is $\phi_E\odot\kappa^0_K = \ell_{K/F}(\phi\odot\kappa^0)$, as required. 
\par
Suppose now that $\vp_F$ is replaced by $\alpha\vp_F$, for some $\alpha\in \bs\mu_F$. The group $C_E$ is then unchanged. If the infinite cyclic group $\pre p{\bk J}(\vp_F)/J^1$ is generated by $\vp_0$, then $\pre p{\bk J}(\alpha\vp_F)/J^1$ is generated by $\alpha_1\vp_0$, where $\alpha_1\in \bs\mu_F$ satisfies $\alpha_1^{p^r} = \alpha$. However, we have remarked that $\bs\mu_F\subset \roman{Ker}\,\kappa^0$. In particular, $\alpha\in \roman{Ker}\,\kappa^0$ and the same representation $\kappa^0$ therefore satisfies the conditions of Lemma 2 relative to $\alpha\vp_F$. Let $h\in \pre p{\bk J_K}(\alpha\vp_F)$; there exists $\ve\in \bs\mu_F$ such that $\ve h\in \pre p{\bk J_K}(\vp_F)$. The relation (5.6.1) reads 
$$ 
\multline 
\eps^0_{\bs\mu_K}\,\eps^1_{\bs\mu_K}(\zeta)\,\roman{tr}\, \kappa^0(\zeta h) = \eps^0_{\bs\mu_K}\,\eps^1_{\bs\mu_K}(\zeta\ve)\,\roman{tr}\, \kappa^0(\zeta \ve h) \\ 
= \roman{tr}\,\kappa^0_K(\zeta \ve h) = \roman{tr}\,\kappa^0_K(h) = \roman{tr}\,\kappa^0_K(\zeta h), 
\endmultline 
$$
since $\kappa^0_K$ is trivial on $\bs\mu_K$ and $\kappa^0$, $\eps^1_{\bs\mu_K}$ are both trivial on $\bs\mu_F$. In other words, the representation $\kappa^0_K$ is also unchanged and the result follows.  \qed 
\enddemo 
\subhead 
5.7 
\endsubhead 
We consider the set $\scr H(\theta_K)^{\vD \text{\rm -reg}}$ of {\it $\vD$-regular\/} elements of $\scr H(\theta_K) = \scr T(\theta_K)$. 
\proclaim{Lemma 1} 
A representation $\vL\in \scr H(\theta_K)$ is $\vD$-regular if and only if it is equivalent to $\xi\odot \ell_{K/F}(\kappa)$, for some $\kappa\in \scr H(\theta)$ and a $\vD$-regular element $\xi$ of $X_1(E)$. 
\endproclaim 
\demo{Proof} 
This follows from 5.6 Proposition and (5.6.2). \qed 
\enddemo 
Let $X_1(E)^{\text{\rm $\vD$-reg}}$ denote the set of $\vD$-regular elements of $X_1(E)$. A character $\xi\in X_1(E)^{\text{\rm $\vD$-reg}}$ determines a representation $\lambda^\bk J_\xi$ of $\bk J/J^1$ as in 3.6. In the notation of 3.6, we may form the representation $\lambda^\bk J_\xi \otimes \kappa = \lambda_\xi\ltimes\kappa \in \scr T(\theta)$, for any $\kappa\in \scr H(\theta)$. 
\proclaim{Lemma 2} 
For any $\kappa\in \scr H(\theta)$, the map $\xi \mapsto \lambda_\xi\ltimes\kappa$ induces an $X_1(E_0)$-isomorphism $\vD\backslash X_1(E)^{\text{\rm $\vD$-reg}} \to \scr T(\theta)$. 
\endproclaim 
\demo{Proof} 
This follows from Lemma 1 and 3.6 Proposition (2). \qed 
\enddemo 
We define a map $\ind_{K/F}: \scr T(\theta_K)^{\text{\rm $\vD$-reg}} \to \scr T(\theta)$ as follows. Choose $\kappa\in \scr H(\theta)$ and let $\vL\in \scr T(\theta_K)^{\text{\rm $\vD$-reg}}$. Using Lemma 1, we may write $\vL = \xi\odot\ell_{K/F}(\kappa)$, for some $\xi\in X_1(E)^{\text{\rm $\vD$-reg}}$. We set 
$$
\ind_{K/F}\,\vL = \lambda_\xi\ltimes\kappa. 
\tag 5.7.1 
$$ 
The definition is independent of the choice of $\kappa\in \scr H(\theta)$ ({\it cf\.} 3.6 Proposition (1)). We so obtain a {\it canonical\/} bijection 
$$
\ind_{K/F}: \vD\backslash \scr T(\theta_K)^{\text{\rm $\vD$-reg}} @>{\ \ \approx\ \ }>> \scr T(\theta). 
\tag 5.7.2 
$$ 
By construction, it satisfies 
$$
\ind_{K/F}\,\big(\phi_E\odot \vL\big) \cong \phi\odot \ind_{K/F}\,\vL, 
\tag 5.7.3 
$$
for all $\vL \in \scr T(\theta_K)^{\text{\rm $\vD$-reg}}$ and $\phi\in X_1(E_0)$. 
\remark{Comments} 
All of the preceding constructions are in terms of a field $P/P_0$, chosen as in 5.4 Notation (5). The field $K/F$ is then defined as the maximal unramified sub-extension of $P/F$. As remarked at the time, any two choices of $P$, hence of $K$, are $J^0_\theta$-conjugate. Indeed, they are conjugate by an element of $J^0_\theta$ commuting with $P_0$. All of the constructions are invariant under such conjugations: if $j\in J^0_\theta\cap G_{P_0}$, then conjugation by $j$ gives a bijection $\scr H(\theta_K) \to \scr H(\theta_{K^j})$ preserving the actions of $X_1(E_0)$. Also, conjugation by $j$ fixes the representations $\lambda_\xi^{\bk J}$ and $\kappa$. Thus $\ind_{K^j/F}\, \vL^j \cong \ind_{K/F}\,\vL$, and so the constructions are independent of the choices of $P$ and $K$. 
\endremark 
\subhead 
5.8 
\endsubhead 
Following 4.2, we translate this machinery to the context of cuspidal representations of the groups $G$ and $G_K$. To do this, we must first ``externalize'' our notation. From this point of view, the hypotheses listed in 5.4 give us: 
\proclaim{Data} 
\roster 
\item 
a finite unramified field extension $K/F$, 
\item 
an endo-class $\vT_K\in \scr E(K)$ with totally ramified parameter field $E/K$, 
\item 
a group $\vD$ of $F$-automorphisms of $E$, of order $m$ and such that $E/E^\vD$ is unramified and 
\item 
a totally wild $E/K$-lift $\vT_E$ of $\vT_K$, 
\endroster 
such that 
\roster 
\item"(5)" $E_0 = E^\vD/F$ is a tame parameter field for $\vT = \frak i_{K/F}\vT_K$. 
\endroster 
\endproclaim 
In particular, $\vT_{E_0} = \frak i_{E/E_0}\vT_E$ is a totally wild $E_0/F$-lift of $\vT$ and $\vT_E$ is the unique $E/E_0$-lift of $\vT_{E_0}$. Thus $\vT_E$ is fixed by $\vD$. We identify $\vD$ with $\Gal K{K\cap E_0}$ by restriction, and $\vD$ then fixes $\vT_K$. 
\par 
Set $n = m\deg\vT$ and $n_K = n/[K{:}F]$. Let $G_K = \GL{n_K}K$. We choose an m-realization $\theta_K$ of $\vT_K$ in $G_K$. If $T/K$ is a tame parameter field for $\theta_K$, we choose a $K$-isomorphism $\bsi:E\to T$ such that $\bsi_*(\vT_K) = \text{\it cl\/}(\theta_K)$. This configuration is determined up to conjugation by an element of the $G_K$-centralizer $G_{\bsi E}$ of $\bsi E^\times$ (as in 2.7). As in (4.1.3), we have the bijection $\scr T(\theta_K) \to \Aa1K{\vT_K}$ given by $\vL\mapsto \cind_{\bk J_{\theta_K}}^{G_K} \vL$. This map takes the natural $\odot$-action of $X_1(E)$ on $\scr H(\theta_K) = \scr T(\theta_K)$ to its $\odot_{\vT_E}$-action on $\Aa1K{\vT_K}$. 
\par
Next, we follow (2.3.3) to produce a simple character $\theta$ in $G = \GL nF$ such that $H^1_\theta\cap G_K = H^1_{\theta_K}$ and $\theta|_{H^1_{\theta_K}} = \theta_K$. Thus $\theta$ is an m-realization of $\vT$, with tame parameter field $\bsi E_0/F$. The character $\theta_K$ determines the character $\theta$ uniquely \cite{3 (7.15)}. We view $\vD$ as acting on $G_K$ via conjugation by elements of $J^0_\theta$, as in 5.6. This induces the natural action of $\vD$ on $\Aa1K{\vT_K}$ and the induction map $\scr T(\theta_K)\to \Aa1K{\vT_K}$ is a $\vD$-map. We likewise have the induction map $\scr T(\theta) \to \Aa mF\vT$, transforming the $\odot$-action of $X_1(E_0)$ on $\scr T(\theta)$ into its $\odot_{\vT_{E_0}}$-action on $\Aa mF\vT$. 
\proclaim{Proposition} 
The map $\ind_{K/F}$ of \rom{(5.7.1)} induces a canonical bijection 
$$
\ind_{K/F}:\vD\backslash \Aa1K{\vT_K}^{\text{\rm $\vD$-reg}} @>{\ \ \approx\ \ }>> \Aa mF\vT 
\tag 5.8.2 
$$
satisfying 
$$
\ind_{K/F}\,\big(\phi_E\odot_{\vT_E} \tau\big) = \phi\odot_{\vT_{E_0}} \ind_{K/F}\,\tau, 
\tag 5.8.3 
$$
for $\tau \in \Aa1K{\vT_K}^{\text{\rm $\vD$-reg}}$ and $\phi\in X_1(E_0)$. 
\par
If $\gamma:K\to K^\gamma$ is an isomorphism of local fields, then 
$$
\ind_{K^\gamma/F^\gamma}\,\tau^\gamma \cong \big(\ind_{K/F}\,\tau\big)^\gamma, 
\tag 5.8.4 
$$ 
for all $\tau\in \Aa1K{\vT_K}^{\text{\rm $\vD$-reg}}$. 
\endproclaim 
\demo{Proof}
The triple $(\theta_K, \bsi)$ is determined up to conjugation by an element of $G_{\bsi E}$, as already noted. Since $\theta$ is determined by $\theta_K$, the triple $(\theta_K,\bsi,\theta)$ is likewise uniquely determined, up to conjugation by $G_{\bsi E}$. The map $\Aa1K{\vT_K}^{\text{$\vD$-reg}} \to \Aa mF\vT$, induced by the map $\ind_{K/F}:\scr T(\theta_K)^{\text{$\vD$-reg}} \to \scr T(\theta)$, is therefore independent of the choice of $(\theta_K,\bsi,\theta)$. 
\par 
That the map $\ind_{K/F}$ of (5.7.1) induces a bijection (5.8.2) follows from (5.7.2). The property (5.8.3) follows from (5.7.3). 
\par
Finally, consider the map $\vL\mapsto \big(\ind_{K^\gamma/F^\gamma}\,\vL^\gamma\big)^{\gamma^{-1}}$, $\vL\in \scr T(\theta_K)^{\text{$\vD$-reg}}$. This is a version of $\ind_{K/F}$ defined relative to different choices of $(\theta_K,\bsi,\theta)$. Such choices have no effect on the map (5.8.1), whence (5.8.4) follows. \qed 
\enddemo 
\subhead 
5.9 
\endsubhead 
We continue in the same situation, as laid out at the beginning of 5.4 but using the viewpoint of 5.8. The group $\vD$ acts on $E$ as $\Gal E{E_0}$. It fixes the endo-class $\vT_E$, so it acts on $\Aa1E{\vT_E}$. Let $\Aa1E{\vT_E}^{\text{\rm $\vD$-reg}}$ be the subset of $\vD$-regular elements. From (5.3.6), applied with base field $K$, we get a canonical bijection $\ind_{E/K}: \Aa1E{\vT_E} \to \Aa1K{\vT_K}$. By 5.3 Corollary 2, this is a $\vD$-map. We therefore define a map 
$$ 
\ind_{E/F}: \Aa1E{\vT_E}^{\text{\rm $\vD$-reg}} \longrightarrow \Aa mF\vT 
$$ 
by 
$$ 
\ind_{E/F} = \ind_{K/F}\circ \ind_{E/K} . 
\tag 5.9.1 
$$ 
We combine 5.3 Corollary 1 with 5.8 Proposition to obtain: 
\proclaim{Theorem} 
Let $\vT\in \scr E(F)$ have tame parameter field $E_0/F$. Let $E/E_0$ be unramified of degree $m$, and set $\vD = \Gal E{E_0}$. Let $\vT_{E_0}$ be a totally wild $E_0/F$-lift of $\vT$, let $\vT_E$ be the unique $E/E_0$-lift of $\vT_{E_0}$. The map 
$$
\ind_{E/F}:\vD\backslash \Aa1E{\vT_E}^{\text{\rm $\vD$-reg}} \longrightarrow \Aa mF\vT 
$$
is a canonical bijection satisfying 
$$
\ind_{E/F}(\phi_E\,\pi) = \phi\odot_{\vT_{E_0}} \ind_{E/F}\,\pi, 
$$
for $\pi\in \Aa1E{\vT_E}^{\text{\rm $\vD$-reg}}$, $\phi\in X_1(E_0)$, and $\phi_E = \phi\circ\N E{E_0}$. 
\par 
The map $\ind_{E/F}$ is natural with respect to isomorphisms of the field $E$. 
\endproclaim 
\head\Rm 
6. Some properties of the Langlands correspondence 
\endhead 
We make some preliminary connections between the machinery developed in \S2--\,\S4 and the representation theory of the Weil group set out in \S1. We use the same notation as before. In particular, $\Go nF$ is the set of equivalence classes of irreducible, smooth, $n$-dimensional representations of $\scr W_F$ and $\wW F = \bigcup_{n\ge1} \Go nF$. It will be convenient to have the analogous notation $\wG F = \bigcup_{n\ge1} \Ao nF$. The Langlands correspondence therefore gives a bijection 
$$ 
\align 
\Bbb L:\wW F &\longrightarrow \wG F, \\
\sigma&\longmapsto \upr L\sigma. 
\endalign 
$$ 
\subhead 
6.1 
\endsubhead 
To a representation $\pi\in \wG F$ we attach the endo-class $\vt(\pi)\in \scr E(F)$ of an m-simple character occurring in $\pi$, as in (4.1.1). On the other hand, let 
$$
r^1_F:\wW F \longrightarrow \scr W_F\backslash \wP F 
\tag 6.1.1 
$$
be the map taking an irreducible smooth representation $\sigma$ of $\scr W_F$ to the $\scr W_F$-orbit $\scr O_F(\alpha)$ of an irreducible component $\alpha$ of $\sigma|_{\scr P_F}$. 
\proclaim{Ramification Theorem} 
There is a unique map 
$$
\Phi_F: \scr W_F\backslash \wP F \longrightarrow \scr E(F) 
$$
such that the following diagram commutes. 
$$
\CD 
\wW F @>{\Bbb L}>> \wG F \\ 
@V{r^1_F}VV @VV{\vt}V \\ 
\scr W_F\backslash \wP F @>>{\Phi_F}> \scr E(F) 
\endCD 
$$ 
The map $\Phi_F$ is bijective. If $\gamma:F\to F^\gamma$ is an isomorphism of topological fields, extended in some way to an isomorphism $\scr W_F \to \scr W_{F^\gamma}$, then 
$$
\Phi_{F^\gamma}(\scr O_{F^\gamma}(\alpha^\gamma)) = \Phi_F(\scr O_F(\alpha))^\gamma. 
\tag 6.1.2 
$$ 
\endproclaim 
\demo{Proof} 
All assertions except the last are 8.2 Theorem of \cite{8}. The last one follows from the uniqueness of $\Phi_F$ and the corresponding property of the Langlands correspondence. \qed 
\enddemo 
If $\alpha\in \wP F$, we usually write $\Phi_F(\alpha)$ rather than $\Phi_F(\scr O_F(\alpha))$. 
\par 
When applied to one-dimensional representations of $\scr P_F$, the Ramification Theorem reduces to the standard ramification theorem of local class field theory. 
\subhead 
6.2 
\endsubhead 
If $K/F$ is finite and tamely ramified, then $\scr P_K = \scr P_F$ and there is an obvious surjective map 
$$ 
\aligned 
\bk i_{K/F}: \scr W_K\backslash \wP F &\longrightarrow \scr W_F\backslash \wP F, \\ 
\scr O_K(\alpha) &\longmapsto \scr O_F(\alpha). 
\endaligned 
\tag 6.2.1 
$$ 
We also have the map $\frak i_{K/F}:\scr E(K) \to \scr E(F)$ of (2.3.4), such that $\frak i_{K/F}^{-1}\vT$ is the set of $K/F$-lifts of $\vT\in \scr E(F)$. 
\proclaim{Proposition} 
Let $K/F$ be a finite, tamely ramified field extension with $K\subset \bar F$. The diagram 
$$ 
\CD 
\scr W_K\backslash\wP K @>{\Phi_K}>> \scr E(K) \\ 
@V{\bk i_{K/F}}VV @VV{\frak i_{K/F}}V \\ 
\scr W_F\backslash\wP F @>>{\Phi_F}> \scr E(F) 
\endCD 
$$ 
is commutative. 
\endproclaim 
\demo{Proof} 
The assertion is transitive with respect to the finite tame extension $K/F$. We may therefore assume that $K/F$ is {\it of prime degree.} 
\par
We deal first with the case where $K/F$ is {\it cyclic.} Let $\alpha\in \widehat{\scr P}_F$, and let $\tau$ be an irreducible representation of $\scr W_K$ such that $\tau|_{\scr P_F}$ contains $\alpha$, that is, $r^1_K(\tau) = \scr O_K(\alpha)$. Set $\phi = \upr L\tau$, so that $\vt(\phi) = \Phi_K(\alpha)$ by the Ramification Theorem. Consider the representation $\sigma = \Ind_{K/F}\,\tau$. 
\par 
Suppose $\sigma$ is not irreducible: equivalently, $\tau^\gamma \cong \tau$ for all $\gamma\in \vG = \Gal KF$. We have  
$$
\sigma = \bigoplus_\chi \chi\otimes \sigma_1, 
$$ 
where $\sigma_1$ is irreducible and $\chi$ ranges over the characters of the group $\scr W_F/\scr W_K = \vG$. All characters $\chi$ appearing here are trivial on $\scr P_K = \scr P_F$, so $r^1_F(\sigma_1)$ is the orbit $\scr O_F(\alpha) = \bk i_{K/F}\scr O_K(\alpha)$. On the other hand, the Langlands correspondence takes $\sigma$ to the representation $\pi = \roman A_{K/F}(\phi)$ automorphically induced by $\phi$ \cite{17 (3.1)}, and this takes the form of a Zelevinsky sum \cite{31} 
$$
\pi = \sqsum\chi\,\chi\pi_1, 
$$
where $\pi_1 = \upr L\sigma_1$ and $\chi$ ranges as before. We have 
$$
\vt(\pi_1) = \Phi_F(r^1_F(\sigma_1)) = \Phi_F(\scr O_F(\alpha)) = \Phi_F\circ\bk i_{K/F}(\scr O_K(\alpha)) 
$$ 
while, according to 7.1 Corollary of \cite{8}, 
$$ 
\vt(\pi_1) = \frak i_{K/F}(\vt(\phi)) = \frak i_{K/F}(\Phi_K(r^1_K(\tau))) = \frak i_{K/F}\circ\Phi_K(\scr O_K(\alpha)). 
$$ 
Therefore 
$$ 
\Phi_F\circ\bk i_{K/F}(\scr O_K(\alpha)) = \vt(\pi_1) = \frak i_{K/F}\circ\Phi_K(\scr O_K(\alpha)), 
$$ 
which proves the proposition in the present case. 
\par 
Suppose now that $\sigma$ is irreducible. The representation $\pi = \upr L\sigma$ is again automorphically induced, $\pi = \roman A_{K/F}(\phi)$ and, in the same way, $\vt(\pi) = \frak i_{K/F}(\vt(\phi))$. The result follows as before. 
\par 
The proposition thus holds when $K/F$ is cyclic of prime degree. By transitivity, it holds when the tame extension $K/F$ is Galois of any finite degree. 
\par
We are reduced to the case where $K/F$ is of prime degree but not cyclic. In particular, $K/F$ is totally tamely ramifed. Let $E/F$ be the normal closure of $K/F$. From the first part of the proof, we know that the result holds for the tamely ramified Galois extensions $E/F$, $E/K$. A diagram chase shows that it holds for $K/F$. \qed 
\enddemo 
\subhead 
6.3 
\endsubhead 
We use 6.2 Proposition to refine the Ramification Theorem. As in (1.5.1), we set 
$$
d_F(\alpha) = [Z_F(\alpha){:}F]\,\dim\alpha, \quad \alpha \in \wP F. 
$$ 
\proclaim{Tame Parameter Theorem} 
Let $\alpha \in \wP F$ have dimension $p^r$, $r\ge0$. If $E = Z_F(\alpha)$, then 
\roster 
\item 
$\deg\Phi_F(\alpha) = d_F(\alpha)$, and 
\item 
$E/F$ is a tame parameter field for $\Phi_F(\alpha)$. 
\endroster 
\endproclaim 
\demo{Proof} 
We recall from \cite{18 (6.2.5)} the following relation. 
\proclaim{Lemma 1} \  
Let $\pi \in \Ao nF$, and denote by $t(\pi)$ the number of characters $\chi\in X_0(F)$ for which $\chi\pi\cong\pi$. We then have $t(\pi) = n/e(\vt(\pi))$. 
\endproclaim 
Similarly, if $\sigma \in \wW F$, we define $t(\sigma)$ to be the number of $\chi \in X_0(F)$ such that $\chi\otimes\sigma \cong \sigma$. Since $\upr L(\chi\otimes \sigma) = \chi\cdot \upr L\sigma$, we have 
$$
t(\sigma) = t(\upr L\sigma). 
\tag 6.3.1 
$$ 
Using 1.4 Theorem to write $\sigma = \vS(\rho,\tau)$, for an admissible datum $(E/F,\rho,\tau)$, we find   
$$ 
t(\sigma) = f(E|F)\,\dim\tau. 
\tag 6.3.2 
$$ 
\indent 
We apply Lemma 1, first in the case where $E=F$. By 1.3 Proposition, there exists $\sigma\in \wW F$ such that $\sigma|_{\scr P_F} \cong \alpha$. The representation $\sigma$ satisfies $t(\sigma) = 1$ (6.3.2) and $\dim\sigma = \dim\alpha = p^r$. If $\pi = \upr L\sigma$, it follows that $t(\pi) = 1$ and then that $e(\vt(\pi)) = \dim\sigma = p^r$. We deduce that $\deg\vt(\pi)= p^r$ and, if $F[\beta]$ is a parameter field for $\vt(\pi)$, then $F[\beta]/F$ is totally wildly ramified. The tame parameter field for $\vt(\pi) = \Phi_F(\alpha)$ is therefore $F$, as required. 
\par 
We pass to the general case $E\neq F$. 
\proclaim{Lemma 2} \ 
The field $E$ contains a tame parameter field $T/F$ for $\Phi_F(\alpha)$, and $\deg\Phi_F(\alpha) = p^r[T{:}F]$. 
\endproclaim 
\demo{Proof} 
According to 6.2 Proposition, we have $\Phi_F(\alpha)= \frak i_{E/F}\Phi_E(\alpha)$. The fibre $\bk i^{-1}_{E/F}(\scr O_F(\alpha))$ is given by the disjoint union 
$$
\bk i^{-1}_{E/F}\scr O_F(\alpha) = \bigcup_{g\in \scr W_E\backslash \scr W_F/\scr W_E} \scr O_E(\alpha^g). 
$$
By the first case above, the term $\Phi_E (\alpha)$ is totally wild of degree $p^r$, but is also an $E/F$-lift of $\Phi_F(\alpha)$. The lemma now follows from 2.4 Proposition. \qed 
\enddemo 
Set $s = p^r[T{:}F]$. There exists an irreducible cuspidal representation $\tau$ of $\GL sF$ with $\vt(\tau) = \Phi_F(\alpha)$. Define $\nu\in \wW F$ by $\upr L\nu = \tau$. Thus $\dim\nu = s$ while $\nu|_{\scr P_F}$ contains the representation $\alpha$. However, an irreducible smooth representation of $\scr W_F$ containing $\alpha$ has dimension divisible by $p^r[E{:}F]$ (1.4 Theorem). Therefore $s=p^r[E{:}F]$. By the lemma, $T\subset E$, so $E=T$. \qed 
\enddemo 
\subhead 
6.4 
\endsubhead 
Let $\scr O \in \scr W_F\backslash \wP F$, let $\vT\in \scr E(F)$ and let $m\ge1$ be an integer. We define sets $\Gg mF{\scr O}$, $\Aa mF\vT$ as in 1.5, (4.1.2) respectively. 
\par
If $\alpha\in \wP F$ and $\scr O = \scr O_F(\alpha)$, we set $d(\scr O) = d_F(\alpha)$. Thus $\Gg mF{\scr O} \subset \Go nF$, where $n = md(\scr O)$,  and $\Go nF$ is the disjoint union 
$$
\Go nF = \bigcup_{(\scr O,m)} \Gg mF{\scr O}, 
\tag 6.4.1 
$$
where $(\scr O,m)$ ranges over all pairs $\scr O\in \scr W_F\backslash \wP F$, $m\in \Bbb Z$, such that $n = md(\scr O)$. Likewise, $\Ao nF$ is the disjoint union 
$$
\Ao nF = \bigcup_{(\vT,m)} \Aa mF\vT, 
\tag 6.4.2 
$$
where $\vT\in \scr E(F)$, $m\in \Bbb Z$ and $n=m\deg\vT$. 
\proclaim{Corollary} 
Let $\scr O \in \scr W_F\backslash \wP F$ and put $\vT = \Phi_F(\scr O) \in \scr E(F)$. For every integer $m\ge1$, the Langlands correspondence induces a bijection 
$$
\Gg mF{\scr O} @>{\ \ \approx\ \ }>> \Aa mF\vT. 
\tag 6.4.3 
$$
\endproclaim 
\demo{Proof} 
Let $\rho\in \Gg mF{\scr O}$. By the Ramification Theorem, $\upr L\rho\in \Aa{m'}F\vT$, where $m'\deg\vT = md(\scr O)$. Part (1) of the Tame Parameter Theorem says $\deg\vT = d(\scr O)$. Consequently 
$$
\upr L{\big(\Gg mF{\scr O}\big)} \subset \Aa mF{\Phi_F(\scr O)}, \quad \scr O\in \scr W_F \backslash\wP F,\ m\ge1. 
$$
The corollary now follows from (6.4.1), (6.4.2). \qed 
\enddemo 
\head\Rm 
7. A na\"\i ve correspondence and the Langlands correspondence 
\endhead 
We use the machinery of \S5, and the relationships uncovered in \S6, to define a canonical bijection 
$$
\align 
\Bbb N: \wW F &\longrightarrow \wG F, \\ \sigma &\longmapsto \upr N\sigma. 
\endalign 
$$ 
We state our main results, comparing this ``na\"\i ve correspondence'' with the Langlands correspondence. 
\subhead 
7.1 
\endsubhead 
Suppose first that $\sigma\in \wW F$ is {\it totally wildly ramified,\/} i.e., $\sigma|_{\scr P_F}$ is irreducible. For such a representation $\sigma$, we set 
$$
\upr N\sigma = \upr L\sigma. 
\tag 7.1.1 
$$ 
Basic properties of the Langlands correspondence \cite{17 (3.1)} imply: 
\proclaim{Proposition} 
Let $\sigma\in \wW F$ be totally wildly ramified. 
\roster 
\item 
If $\chi$ is a character of $F^\times$, then 
$$
\upr N(\chi\otimes\sigma) = \chi\cdot \upr N\sigma. 
\tag 7.1.2 
$$
\item 
If $\gamma:F\to F^\gamma$ is an isomorphism of local fields, then 
$$
\upr N(\sigma^\gamma) = (\upr N\sigma)^\gamma, 
\tag 7.1.3 
$$
\endroster 
\endproclaim 
\subhead 
7.2 
\endsubhead 
Let $\sigma \in \wW F$. As in 6.4, there is a unique pair $(m,\scr O)$, where $m\ge 1$ and $\scr O\in \scr W_F\backslash \wP F$, such that $\sigma \in \Gg mF{\scr O}$. Choose $\alpha\in \scr O$, and set $E = Z_F(\alpha)$. Let $E_m/E$ be unramified of degree $m$ and put $\vD = \Gal {E_m}E$. As in 1.6 Proposition, there is a $\vD$-regular representation $\rho\in \Gg1{E_m}\alpha$ such that $\sigma \cong \Ind_{E_m/F}\, \rho$. The $\vD$-orbit of $\rho$ is thereby uniquely determined. 
\par
By (7.1.3) and 6.4 Corollary, the representation $\upr N\rho$ is a $\vD$-regular element of $\Aa1{E_m} {\Phi_{E_m}(\alpha)}$. We use the map $\ind_{E_m/F}$ of  (5.9.1) to define 
$$
\upr N\sigma = \text{\it ind}_{E_m/F}\,\upr N\rho. 
\tag 7.2.1 
$$ 
\proclaim{Proposition} 
For each integer $m\ge1$ and each $\scr O\in \scr W_F\backslash \wP F$, the assignment \rom{(7.2.1)} induces a bijection 
$$
\aligned 
\Gg mF{\scr O} &\longrightarrow \Aa mF{\Phi_F(\scr O)}, \\ 
\sigma &\longmapsto \upr N\sigma, 
\endaligned 
\tag 7.2.2 
$$ 
which does not depend on the choice of $\alpha\in \scr O$ used in its definition. Further,  
\roster 
\item 
if $x\mapsto x^\gamma$ is an isomorphism $F \to F^\gamma$ of local fields, then 
$$
\upr N(\sigma^\gamma) = \big(\upr N\sigma\big)^\gamma,  \quad \sigma \in \Gg mF{\scr O}, 
$$ 
and 
\item 
if $\alpha\in \scr O$ and $E = Z_F(\alpha)$, then  
$$
\upr N(\phi\odot_\alpha\sigma) = \phi\odot_{\Phi_E(\alpha)} \upr N\sigma,  
$$ 
for all $\sigma\in \Gg mF{\scr O}$ and all $\phi\in X_1(E)$. 
\endroster 
\endproclaim 
\demo{Proof} 
The map $\Ind_{E_m/F}$ induces a canonical bijection 
$$ 
\vD\backslash \Gg1{E_m }\alpha^{\vD \text{\rm -reg}} \longrightarrow \Gg mF{\scr O} 
$$ 
(1.6 Proposition). The map $\text{\ind}_{E_m/F}$ induces a canonical bijection 
$$ 
\vD\backslash \Aa1{E_m }{\Phi_{E_m}(\alpha)}^{\vD \text{\rm -reg}} \longrightarrow \Aa mF{\Phi_F(\alpha)} 
$$ 
(5.9 Theorem). The na\"\i ve correspondence induces a $\vD$-bijection $\Gg1{E_m }\alpha \to \Aa1{E_m}{\Phi_{E_m}(\alpha)}$ (7.1 Proposition). The map (7.2.2) is therefore a bijection. If $\alpha'\in \scr O$, there exists $\gamma\in \scr W_F$ such that $\alpha' = \alpha^\gamma$, whence (by 5.9 Theorem, (6.1.2) and 7.1 Proposition) the definition of $\upr N\sigma$ does not depend on $\alpha$. Assertions (1) and (2) likewise follow from 5.9 Theorem. \qed 
\enddemo 
In the notation of the proposition, the central character $\omega_\pi$ of $\pi = \upr L\sigma$ is given by 
$$ 
\omega_\pi = \det\sigma. 
\tag 7.2.3 
$$ 
If we write $\sigma = \Ind_{E_m/F}\,\rho$, where $\rho\in \Gg1{E_m}\alpha$ is $\vD$-regular, then 
$$
\det\sigma = (d_{E_m/F})^{p^r}\,\det\rho|_{F^\times}, 
\tag 7.2.4 
$$ 
where $d_{E_m/F} \in X_1(F)$ is the discriminantal character of the extension $E_m/F$, 
$$
d_{E_m/F} = \det \Ind_{E_m/F}\,1_{E_m}, 
$$ 
and $p^r = \dim\rho = \dim\alpha$. On the other hand, setting $\pi' = \upr N\sigma$, we get 
$$
\omega_{\pi'} = \det\rho|_{F^\times}. 
\tag 7.2.5 
$$ 
That is, $\omega_\pi/\omega_{\pi'} = (d_{E_m/F})^{p^r}$, a relation which will be useful later. 
\remark{Remark} 
Suppose that $\sigma\in \wW F$ is {\it essentially tame,\/} that is, $r^1_F(\sigma) = \scr O_F(\alpha)$, where $\dim\alpha =1$. The definition of $\upr N\sigma$ is then equivalent to that of \cite{10 (2.3 Theorem)}, the first step (7.1.1) in that case being given by local class field theory. 
\endremark 
\subhead 
7.3 
\endsubhead 
We state our main result. 
\proclaim{Comparison Theorem} 
Let $m\ge 1$ be an integer and $\alpha\in \wP F$. Let $E = Z_F(\alpha)$ and put $\scr O = \scr O_F(\alpha) \in \scr W_F\backslash\wP F$. There exists a character $\mu = \mu^F_{m,\alpha} \in X_1(E)$ such that 
$$
\upr L\sigma = \mu\odot_{\Phi_E(\alpha)} \upr N\sigma, 
\tag 7.3.1 
$$
for all $\sigma \in \Gg mF{\scr O}$. The character $\mu^F_{m,\alpha}$ is thereby uniquely determined modulo $X_0(E)_m$. 
\endproclaim 
The proof of the Comparison Theorem will occupy the rest of the paper. 
\remark{Remark 1} 
Let $\sigma\in \Gg mF{\scr O}$, as in the statement of the theorem. Let $E_m/E$ be unramified of degree $m$ and set $\vD = \Gal {E_m}E$ (as in the definition (7.2.1)). The theorem asserts that the character $\mu^F_{m,\alpha}\circ \N {E_m}E \in X_1(E_m)^\vD$ is uniquely determined. It is often more convenient to use this version. 
\endremark 
\remark{Remark 2}
Let $\dim\alpha = p^r$, let $G_E = \GL{mp^r}E$ and let $\det_E:G_E\to E^\times$ be the determinant map. Setting $\pi = \upr L\sigma$ and $\pi' = \upr N\sigma$ as before, (7.2.3--5) yield 
$$ 
\mu^F_{m,\alpha}\circ{\det}_E|_{F^\times} = \omega_\pi/\omega_{\pi'} = (d_{E_m/F})^{p^r}. 
$$ 
In other words, 
$$
\mu^F_{m,\alpha}(x)^{mp^r} = d_{E_m/F}(x)^{p^r}, \quad x\in F^\times. 
\tag 7.3.2 
$$
\endremark 
\subhead 
7.4 
\endsubhead 
Combining the Comparison Theorem with 7.2 Proposition, we obtain a valuable corollary. 
\proclaim{Homogeneity Theorem} 
If $\alpha\in \wP F$ and $E = Z_F(\alpha)$, then 
$$
\upr L(\phi\odot_\alpha\sigma) = \phi\odot_{\Phi_E(\alpha)} \upr L\sigma, 
$$
for all $\sigma\in \wW F$ such that $r^1_F(\sigma) = \scr O_F(\alpha)$ and all $\phi\in X_1(E)$. 
\endproclaim 
\demo{Proof} 
Abbreviating $\scr O = \scr O_F(\alpha)$, let $\sigma\in \Gg mF{\scr O}$ for some $m\ge1$. Let $\mu = \mu^F_{m,\alpha}$, in the notation of the Comparison Theorem. For $\phi \in X_1(E)$, we have 
$$ 
\align 
\upr L(\phi\odot_\alpha\sigma) &= \mu\odot_{\Phi_E(\alpha)} \phi\odot_{\Phi_E(\alpha)} \upr N\sigma \\ 
&= \phi \odot_{\Phi_E(\alpha)} \mu\odot_{\Phi_E(\alpha)} \upr N\sigma = \phi\odot_{\Phi_E(\alpha)} \upr L\sigma, 
\endalign 
$$
as required. \qed 
\enddemo 
In general, the Homogeneity Theorem is weaker than the Comparison The\-orem. However, the two are equivalent in the case $m=1$ since the sets $\Gg1F{\scr O}$, $\Aa1F{\Phi_F(\alpha)}$ are then principal homogeneous spaces over $X_1(E)$ (1.5 Proposition, 4.2 Proposition). 
\subhead 
7.5 
\endsubhead 
We comment on a well-known special case \cite{24} or \cite{13 (2.4)}. 
\par
Let $\rho \in \wW F$ be {\it tamely ramified\/} of dimension $n$, that is, $\rho\in \Gg nF{\bk 1_F}$. Thus $\rho = \rho_\xi = \Ind_{K/F}\,\xi$, where $K/F$ is unramified of degree $n$ and $\xi\in X_1(K)$ is $K/F$-regular. As in 3.5, we form the extended maximal simple type $\lambda_\xi\in \scr T_n(\bk 0_F)$. The representation $\upr N\rho_\xi$ is then the representation $\pi_\xi\in \Ao nF$ containing $\lambda_\xi$. 
\par
On the other hand, let $\chi_2$ denote the unramified character of $K^\times$ of order $2$. The representation $\upr L\rho_\xi$ is then given by 
$$
\upr L\rho_\xi = \pi_{\xi'}, \quad \text{where} \quad \xi'= \chi_2^{n-1}\xi 
\tag 7.5.1 
$$
\cite{13 (2.4 Theorem 2)}. In the notation of the Comparison Theorem therefore, 
$$
\mu^F_{n,\bk 1_F}\circ\N KF = \chi_2^{n-1}. 
\tag 7.5.2 
$$ 
\subhead 
7.6 
\endsubhead 
We do not give an explicit formula for the ``discrepancy character'' $\mu^F_{m,\alpha}$. However, we will be able to describe its restriction to units and so obtain a succinct account of the maximal simple type occurring in $\upr L\sigma$, for any $\sigma\in \wW F$. 
\par
We follow the description of representations of the Weil group given by 1.4 Theorem. So, let $(E/F,\rho,\tau)$ be an admissible datum, in which $\dim\tau = m\ge1$, and set $\vS = \vS(\rho,\tau) = \Ind_{E/F}\, \rho\otimes \tau$. Following 1.3 Proposition, we may choose the datum $(E/F,\rho,\tau)$ so that $\det\rho|_{\scr I_E}$ has $p$-power order. Set $\alpha = \rho|_{\scr P_F}$, and put $\vT = \Phi_F(\alpha)$. 
\par
We choose an m-realization $\theta$ of $\vT$ in $G = \GL nF$, where $n = \dim\vS$. We identify $E$ with a tame parameter field for $\theta$, in such a way that the simple character $\theta_E$ (as in 2.3) has endo-class $\Phi_E(\alpha)$ ({\it cf\.} 2.7 Proposition). Thus $E$ is contained in a parameter field $P/F$ for $\theta$. We let $P_m/P$ be unramified of degree $m$, with $P_m^\times\subset \bk J_\theta$. Let $K_m/F$ be the maximal unramified sub-extension of $P_m/F$. 
\par 
Let $\lambda_\tau \in \scr T_m(\bk 0_E)$ be the extended maximal simple type occurring in $\upr L\tau$: this has been described in 7.5. Using this notation, we describe a maximal simple type occurring in $\vS$.  
\proclaim{Types Theorem} 
\roster 
\item 
There exists $\nu\in \scr H(\theta)$ such that $\roman{tr}\,\nu$ is constant on the set of $K_m/F$-regular elements of $\bs\mu_{K_m}$. This condition determines $\nu|_{J^0_\theta}$ uniquely. 
\item 
The representation $\upr L\vS$ contains an element $\vL$ of $\scr T(\theta)$ of the form 
$$
\vL = \psi \odot \lambda_\tau \ltimes\nu, 
$$
where $\psi\in X_1(E)$ satisfies the following condition. If $e(E|F)$ is odd, then $\psi$ is unramified, while $\psi|_{U_E}$ has order\/ $2$ if $e(E|F)$ is even. 
\endroster 
\endproclaim 
We prove this result in \S12. 
\subhead 
7.7 
\endsubhead 
We describe briefly a variant on the Types Theorem. In the same situation as 7.6, let $\sigma_1 = \Ind_{E/F}\,\rho\in \Gg 1F\vT$, and set $\pi_1 = \upr L\sigma_1$. In particular, $\dim\sigma_1 = n_1 = n/m$. Let $\vt = \upr L\tau \in \Aa mE{\bk 0_E}$. 
\par
We choose an m-simple character $\theta_1$ in $\GL{n_1}F$ of endo-class $\vT$. The representation $\pi_1$ then contains an element $\kappa_1$ of $\scr H(\theta_1) = \scr T(\theta_1)$, while $\vt$ contains some $\lambda \in \scr T^0_m(\bk 0_E)$. Using the notation of 3.3 Corollary 1, the representation $\kappa_1$ gives rise to a representation $\kappa_m = f_m(\kappa_1) \in \scr H(\theta)$. We set $\vL' = \lambda\ltimes \kappa_m \in \scr T(\theta)$ and form $\pi_m = \cind_{\bk J_\theta}^G\,\vL' \in \Aa mF\vT$, where $G = \GL nF$. One sees easily that the representation $\pi_m$ is of the form 
$$
\pi_m = \phi\odot_{\Phi_E(\alpha)} \upr L\vS, 
\tag 7.7.1 
$$
for some $\phi\in X_1(E)$ depending only on $m$ and $\alpha$. 
\subhead 
7.8 
\endsubhead 
The first step (7.1.1) in the definition of the na\"\i ve correspondence is surely the natural one. However, any family of maps satisfying (7.1.2) and (7.1.3) would lead to the same outcome: only the character $\mu^F_{m,\alpha}$ would change by a constant factor depending on this choice. 
\par 
More particularly, in the context of 7.2, each of $\Gg1{E_m}\alpha$, $\Aa1 {E_m}{\Phi_{E_m}(\alpha)}$ carries an action of the group $\vD\ltimes X_1(E_m)$, and each is $\vD \ltimes X_1(E_m)$-isomorphic to $X_1(E_m)$. Any two $\vD\ltimes X_1(E_m)$-bijections of $\Gg1{E_m}\alpha$ with $\Aa1 {E_m}{\Phi_{E_m}(\alpha)}$ therefore differ by a constant element of $X_1(E_m)^\vD$. 
\par
There is indeed a plausible alternative to (7.1.1), better preserving the spirit of explicitness. Let $\rho\in \wW F$ be totally wildly ramified of dimension $p^r$, $r\ge1$. In \cite{4}, using methods not dissimilar to those of this paper, we constructed a totally ramified cuspidal representation $\upr W{\!\rho}$ of $\GL{p^r}F$. The ``totally wild'' correspondence $\rho\mapsto \upr W{\!\rho}$ has properties (7.1.2) and (7.1.3). There is an unramified character $\chi_\alpha$ of $F^\times$, of order dividing $p^{r-1}$ and depending only on $\alpha = \rho|_{\scr P_F}$, such that $\upr W{\!\rho} = \chi_\alpha\cdot\upr L\rho$. In the definition of the na\"\i ve correspondence $\sigma\mapsto \upr N\sigma$, we could equally well start by setting $\upr N\sigma = \upr W{\!\sigma}$ when $\sigma$ is totally wildly ramified. The relationship between the Langlands correspondence and the totally wild correspondence is discussed in \cite{5}, \cite{6}. 
\head\Rm 
8. Totally ramified representations 
\endhead 
We prove the Comparison Theorem for {\it totally ramified\/} representations. Thus, in terms of the statement in 7.3, we have $m=1$ and $E/F$ is totally tamely ramified of degree $e$. We consequently revert to the notation listed at the beginning of 5.3. 
\subhead 
8.1 
\endsubhead 
We work relative to a simple character $\theta\in \scr C(\frak a,\beta)$ in $G = \GL nF$ such that $\vT = \text{\it cl\/}(\theta) = \Phi_F(\alpha)$ and $\deg\vT=n$. Thus $P = F[\beta]/F$ is totally ramified of degree $n$. We identify $E/F$ with the maximal tamely ramified sub-extension of $P/F$, in such a way that $\vT_E = \text{\it cl\/}(\theta_E) = \Phi_E(\alpha)$ ({\it cf\.} 2.7). We use the standard abbreviations $\bk J = \bk J(\beta,\frak a) = \bk J_\theta$ and so on. We follow the notational conventions, relative to centralizers of subfields of $\M nF$, laid out in 2.3. 
\par 
Set $\vG = \roman{Aut}(E|F)$. Let $N_E$ denote the $G$-normalizer $N_G(E^\times)$ of $E^\times$. If $\gamma\in \vG$, there exists $g_\gamma\in N_E$ such that $g_\gamma^{-1}xg_\gamma = x^\gamma$, for every $x\in E$. The element $g_\gamma$ is uniquely determined modulo $G_E$, and $g_\gamma\mapsto \gamma$ provides a canonical isomorphism $N_E/G_E\cong \vG$. If $\tau$ is a smooth representation of $G_E$, the equivalence class of the representation $\tau^\gamma:x\mapsto \tau(g_\gamma xg_\gamma^{-1})$ then depends only on that of $\tau$, and not on the choice of $g_\gamma$. 
\par 
We define the constant $\eps = \eps_{E/F} = \pm1$ by (5.3.1) ({\it cf\.} (5.3.2)). For $\pi\in \Aa1F\vT$, we define $\pi_E\in \Aa1E{\vT_E}$ as in 5.3 Corollary 1. As usual, we denote by $\LR G$ the set of {\it elliptic regular\/} elements of $G$. 
\proclaim{Proposition} 
Let $\pi\in \Aa1F\vT$. If $h\in G_E\cap \LR G$ and $\ups_F(\det h)$ is relatively prime to $n$, then 
$$
\roman{tr}\,\pi(h) = \eps_{E/F}\sum_{\gamma\in  \vG} \roman{tr}\,\pi_E^\gamma(h). 
\tag 8.1.1
$$ 
\endproclaim 
\demo{Proof} 
We write $\pi = \cind_{\bk J}^G\,\vL$, where $\vL\in \scr T(\theta)$. Thus $\pi_E = \cind_{\bk J_E}^{G_E}\,\vL_E$, where $\vL_E\in \scr T(\theta_E)$ is given by (5.3.3). 
\par
We examine more closely the elements $h$ of the statement. 
\proclaim{Lemma 1} 
Let $h\in G_E\cap \LR G$ and suppose $\ups_F(\det h)$ is relatively prime to $n$. 
\roster 
\item 
The algebra $F[h]$ is a field, and the extension $F[h]/F$ is totally ramified of degree $n$. 
\item 
The field $F[h]$ contains $E$, and $E/F$ is the maximal tamely ramified sub-extension of $F[h]/F$. 
\endroster 
\endproclaim 
\demo{Proof}
Since $h\in \LR G$, the algebra $F[h]$ is a field, of degree $n$ over $F$. The condition on $\ups_F(\det h)$ ensures that $F[h]/F$ is totally ramified. Since $F[h]$ is a maximal subfield of $A$, it is its own centralizer in $A$. Since $h$ commutes with $E$, the field $F[h]$ must contain $E$. Therefore $[F[h]{:}E] = n/[E{:}F] = p^r$, whence (2) follows. \qed 
\enddemo 
We henceforward view the character $x\mapsto \roman{tr}\,\vL(x)$ of $\vL$ as a function on $G$, vanishing outside $\bk J$. For $g\in \LR G$, the Mackey formula of \cite{3, Appendix} then gives 
$$
\roman{tr}\,\pi(g) = \sum_{x\in G/\bk J} \roman{tr}\,\vL(x^{-1}gx). 
\tag 8.1.2
$$ 
There are only finitely many non-zero terms in the right hand side of this expansion \cite{14 (1.2 Lemma)}. Similar considerations apply to the representations $\pi_E^\gamma$ of $G_E$. Since we will only be concerned with term-by-term comparisons of these character expansions, we may re-arrange terms at will. 
\par 
We eliminate one class of elements $h$. Suppose that $h\in G_E\cap \LR G$ has no $G$-conjugate lying in $\bk J$. It then follows from (8.1.2) that $\roman{tr}\,\pi(h) = 0$. The element $h$ also has no $N_E$-conjugate in $\bk J_E$, whence both sides of (8.1.1) vanish. So, we need only consider those elements $h$ of $G_E\cap \LR G$ having a $G$-conjugate in $\bk J$. 
\proclaim{Lemma 2} 
\roster 
\item 
Let $h\in G_E\cap \LR G$, and suppose that $\ups_F(\det h)$ is relatively prime to $n$. Let $x\in G$, and suppose $x^{-1}hx\in \bk J$. There then exists $y\in N_E$ such that $x\bk J\cap N_E = y\bk J_E$. 
\item 
The map 
$$
\align 
N_E/\bk J_E &\longrightarrow G/\bk J, \\ 
y\bk J_E &\longmapsto y\bk J, 
\endalign 
$$
is injective. 
\endroster 
\endproclaim 
\demo{Proof}
In the present case, we have $\bk J = P^\times J^1$, whence $\bk J\cap N_E = \bk J_E$ (2.6 Proposition (2)) and (2) follows. In (1), let $h' = x^{-1}hx \in \bk J$, By the Conjugacy Lemma (2.6) and 5.5 Lemma, there exists $j\in J^1$ such that $h'' = j^{-1}h'j\in \bk J_E$. Surely $h''\in G_E\cap \LR G$, and we can apply Lemma 1: the field $F[h'']$ contains $E$ and $E/F$ is the maximal tamely ramified sub-extension of $F[h'']/F$. Conjugation by the element $y=xj$ gives an $F$-isomorphism $F[h]\to F[h'']$ which must carry $E$ to itself. That is, $y\in N_E$ and so the lemma is proved. \qed 
\enddemo 
We return to the expansion (8.1.2). Following Lemma 2, it reads 
$$
\roman{tr}\,\pi(h) = \sum_{x\in G/\bk J} \roman{tr}\,\vL(x^{-1}hx) 
= \sum_{y\in N_E/\bk J_E} \roman{tr}\,\vL(y^{-1}hy) ,  
$$ 
with $h\in G_E\cap \LR G$. For $y\in N_E$, we have $y^{-1}hy\in \bk J$ if and only if $y^{-1}hy\in \bk J_E$. In all cases, therefore, (5.3.3) yields 
$$ 
\roman{tr}\,\vL(y^{-1}hy) = \eps_{E/F}\,\roman{tr}\,\vL_E(y^{-1}hy). 
$$ 
Consequently, 
$$ 
\align 
\roman{tr}\,\pi(h) &= \eps_{E/F} \sum_{y\in N_E/\bk J_E} \roman{tr}\,\vL_E(y^{-1}hy) \\
&= \eps_{E/F} \sum_{\gamma\in\vG} \sum_{z\in G_E/\bk J_E} \roman{tr}\,\vL_E(z^{-1}h^\gamma z) \\ 
&= \eps_{E/F} \sum_{\gamma\in \vG} \roman{tr}\,\pi_E(h^\gamma), 
\endalign 
$$
as required. \qed 
\enddemo 
We record a technical consequence for later use. 
\proclaim{Corollary} 
Let $a$ be an integer, relatively prime to $n$. There exists $h\in \bk J_E\cap \LR G$ such that $\ups_F(\det h) = a$ and 
$$
\roman{tr}\,\pi(h) = \eps_{E/F} \sum_{\gamma\in \vG} \roman{tr}\,\pi_E^\gamma(h) \neq 0. 
$$
\endproclaim 
\demo{Proof} 
The representation $\pi$ is totally ramified, in that $t(\pi) = 1$ ({\it cf\.} 6.3 Lemma 1). The representations $\chi\pi$, $\chi\in X_0(F)_n$, are therefore distinct and the character functions $\roman{tr}\,\chi\pi$ are linearly independent on $\LR G$. It follows that there exists $h\in \LR G$ such that $\ups_F(\det h) = a$ and $\roman{tr}\,\pi(h)\neq 0$. The Mackey formula (8.1.2) shows we may as well take $h\in \bk J = P^\times J^1$. Applying the Conjugacy Lemma of 2.6 and 5.5 Lemma, there exists $j\in J^1$ such that $h' = j^{-1}hj\in \bk J_E$. The element $h'$ has all the desired properties. \qed 
\enddemo 
\subhead 
8.2 
\endsubhead 
Let $\rho\in \Gg1E\alpha$. We put 
$$
\nu_\rho = \upr L\rho,\quad \pi_\rho = \upr L(\Ind_{E/F}\,\rho). 
\tag 8.2.1 
$$
Thus $\nu_\rho\in \Aa1E{\vT_E}$ and $\pi_\rho\in \Aa1F\vT$ (6.4 Corollary). 
\par 
The representation $(\pi_\rho)_E$ lies in $\Aa1E{\vT_E}$ (5.3 Corollary 1) and $\Aa1E{\vT_E}$ is a principal homogeneous space over $X_1(E)$ (4.2 Proposition). It follows that there is a unique $\phi_\rho\in X_1(E)$ satisfying $(\pi_\rho)_E \cong \phi_\rho\nu_\rho$. We prove: 
\proclaim{Proposition} The function 
$$ 
\rho\longmapsto \phi_\rho, \quad \rho\in \Gg 1E\alpha, 
\tag 8.2.2 
$$ 
is constant. 
\endproclaim 
Before starting the proof of the proposition, we note that it implies the Comparison Theorem in this case: 
\proclaim{Corollary} 
There is a unique character $\mu = \mu_{1,\alpha}^F\in X_1(E)$ such that 
$$ 
\aligned 
\upr L(\Ind_{E/F}\,\sigma) &= \text{\it ind}_{E/F}\,\mu\cdot\upr L\sigma \\ 
&= \mu\odot_{\vT_E} \upr N\sigma, \quad \sigma\in \Gg1E\alpha.  
\endaligned 
\tag 8.2.3 
$$
\endproclaim 
\demo{Proof}
Continuing with the notation of the proposition, let $\phi$ denote the character $\phi_\rho$, defined by the condition $(\pi_\rho)_E = \phi\nu_\rho$, for all $\rho \in \Gg1E\alpha$. 
\par 
Let $\xi\in X_1(E)$ and consider the representation $\pi_{\xi\otimes\rho}$. According to the proposition, 
$$
(\pi_{\xi\otimes\rho})_E = \phi\nu_{\xi\otimes\rho} = \phi\xi\,\nu_\rho. 
$$
Invoking 5.3 Corollary 1, this relation is equivalent to $\pi_{\xi\otimes\rho} = \xi\odot_{\vT_E} \pi_\rho$. In other words, the map $\Gg1F{\scr O_F(\alpha)} \to \Aa1F\vT$, induced by the Langlands correspondence, is an isomorphism of $X_1(E)$-spaces. However, the composite map 
$$ 
\Gg1F{\scr O_F(\alpha)} \longrightarrow \Gg1E\alpha \longrightarrow \Aa1E{\vT_E} \longrightarrow \Aa1F{\scr O_F(\alpha)} 
$$
is also an isomorphism of $X_1(E)$ spaces. Since $\Gg1F{\scr O_F(\alpha)}$, $\Aa1F\vT$ are principal homogeneous spaces over $X_1(E)$ (1.5 Corollary, 4.2 Proposition), these two maps differ by a constant translation. This is precisely the assertion of the corollary. \qed 
\enddemo 
\subhead 
8.3 
\endsubhead 
We prove 8.2 Proposition by induction on $n = [E{:}F]$. 
\par 
We observe that $|\vG| = \roman{gcd}(e,q{-}1)$ where, we recall, $q = |\Bbbk_F|$ and $e = [E{:}F]$. We first treat the case where $|\vG| = 1$. Here, the norm map $\N EF:E^\times \to F^\times $ is surjective, and induces an isomorphism $E^\times/U^1_E \cong F^\times /U^1_F$. It therefore induces an isomorphism $X_1(F)\to X_1(E)$, denoted $\chi\mapsto \chi_E$. Let $\sigma\in \Gg1F{\scr O_F(\alpha)}$, $\phi\in X_1(E)$. Writing $\phi = \chi_E$, $\chi\in X_1(F)$, we get 
$$
\upr L(\phi\odot_\alpha\sigma) = \upr L(\chi\otimes \sigma) = \chi\cdot\upr L\sigma = \phi\odot_{\vT_E} \upr L\sigma.  
$$
The Langlands Correspondence thus induces an  isomorphism $\Gg1F{\scr O_F(\alpha)} \to \Aa1F\vT$ of principal homogeneous spaces over $X_1(E)$. The same applies to the na\"\i ve correspondence. The two correspondences therefore differ by a constant translation. The Comparison Theorem follows in this case, and with it the proposition of 8.2. 
\remark{Remark} 
In this case, the ramification index $e = e(E|F)$ is odd, so the discriminant character $d_{E/F}$ is unramified, of order $\le2$. It follows easily that the character $\mu = \mu^F_{1,\alpha}$ of the Comparison Theorem is unramified, of order dividing $2n$. 
\endremark 
\subhead 
8.4 
\endsubhead 
We therefore assume that $|\vG| \neq 1$, and we let $l$ be the largest prime divisor of $|\vG|$. Let $K/F$ be the unique sub-extension of $E/F$ of degree $l$. Thus $K/F$ is cyclic. We set $\vO = \Gal KF$. 
\par 
Let $\rho\in \Gg1E\alpha$ and define $\tau_\rho = \upr L(\Ind_{E/K}\,\rho) \in \Ao{n/l}K$. Taking $\pi_\rho = \upr L(\Ind_{E/F}\,\rho) \in \Ao nF$ as in 8.2, the representation $\pi_\rho$ is {\it automorphically induced\/} by $\tau_\rho$. We use the notation $\pi_\rho = \roman A_{K/F}\,\tau_\rho$.  
\par 
The representation $\pi_\rho$ lies in $\Aa1F\vT$, so $\pi_\rho \cong \cind_\bk J^G\,\vL$, for a unique element $\vL$ of $\scr T(\theta)$. In particular, $\pi_\rho$ contains the simple character $\theta\in \scr C(\frak a,\beta)$ of 8.1, and $E$ is identified with a subfield of $P = F[\beta]$. 
\par 
Following the conventions of \cite{13, \S1}, we choose a {\it transfer system\/} for $K/F$, in relative dimension $n/l$ and based on a character $\vk$ of $F^\times$. This gives us a transfer factor $\bs\delta = \bs\delta_{K/F}$. By definition, the kernel of $\vk$ is the norm group $\N KF(K^\times)$. Therefore $\vk$ is trivial on $\det \bk J \subset \N EF(E^\times)U^1_F$, and we may apply the Uniform Induction Theorem of \cite{13 (1.3)}. This gives the relation 
$$
\roman{tr}^\vk\,\pi_\rho(h) = c_\alpha^{K/F}\,\bs\delta(h) \sum_{\omega\in \vO} \roman{tr}\,\tau_\rho^\omega(h), 
\tag 8.4.1 
$$
valid for all $h\in G_K\cap \LR G$. In (8.4.1), the constant $c_\alpha^{K/F}$ depends on the choice of transfer system (which we regard as fixed for all time). Otherwise, it depends only on the endo-class $\text{\it cl\/}(\theta_K) = \Phi_K(\alpha) = \frak i_{E/K}\Phi_E(\alpha)$, and so only on $\alpha$. The transfer factor $\bs\delta$ is independent of all considerations of representations. The term $\roman{tr}^\vk\,\pi_\rho$ is the normalized $\vk$-twisted trace of $\pi_\rho$, formed relative to the $\theta$-normalized $\vk$-operator $\vF_{\pi_\rho}^\vk$ on $\pi_\rho$ (as in \cite{13 (1.3)}). It satisfies a Mackey formula 
$$
\roman{tr}^\vk\,\pi_\rho(h) = \sum_{x\in G/\bk J_F} \vk(\det x^{-1})\,\roman{tr}\,\vL(x^{-1}hx), \quad h\in G_K\cap \LR G, 
\tag 8.4.2 
$$ 
and this is the only property we shall use. 
\par 
We recall that if we treat $\tau_\rho$ as fixed, then the relation (8.4.1) determines $\pi_\rho$ uniquely. Conversely, if we view $\pi_\rho$ as given, then (8.4.1) determines the orbit $\{\tau_\rho^\omega:\omega\in \vO\}$. 
\par 
Let $h\in G_E\cap \LR G$, and suppose that $\ups_F(\det h)$ is relatively prime to $n$. Using (8.4.2), the argument of 8.1 Proposition applies unchanged to give 
$$
\roman{tr}^\vk\,\pi_\rho(h) = \eps_{E/F}\,\sum_{\gamma\in \vG} \vk(\det g^{-1}_\gamma)\, \roman{tr}\,\pi_{\rho,E}^\gamma(h). 
\tag 8.4.3 
$$
\proclaim{Lemma} 
Let $h\in G_E\cap \LR G$ and $\omega\in \vO$. Suppose that $\ups_F(\det h)$ is relatively prime to $n$ and $\roman{tr}\,\tau_\rho^\omega(h) \neq 0$. The automorphism $\omega$ then extends to an $F$-automorphism of $E$. 
\endproclaim 
\demo{Proof} 
Choose $f_\omega\in N_K = N_G(K^\times)$ such that $f_\omega^{-1} zf_\omega = z^\omega$, for all $z\in K^\times$. Since $\roman{tr}\,\tau_\rho^\omega(h) = \roman{tr}\,\tau_\rho(f_\omega h f_\omega^{-1}) \neq 0$, the Mackey formula (8.1.2) implies that $h' = f_\omega h f_\omega^{-1}$ has a $G_K$-conjugate $h'' = x^{-1}h'x$ lying in $\bk J_K$. By the Conjugacy Lemma 2.6, there exists $j\in J^1_K$ such that $h_1 = j^{-1}h''j\in \bk J_E$. The field $K[h_1]$ therefore contains $E$, and $E/K$ is the maximal tamely ramified sub-extension of $K[h_1]/K$. Likewise, $E/K$ is the maximal tamely ramified sub-extension of $K[h]/K$. Conjugation by the element $y = f_\omega^{-1}xj$ thus induces an $F$-isomorphism $K[h]\to K[h_1]$, stabilizing $E$ and extending the automorphism $\omega$ of $K$. Conjugation by $y$ thus induces an automorphism of $E/F$ with the desired property. \qed 
\enddemo 
\subhead 
8.5 
\endsubhead 
Before starting our analysis of the induction relation (8.4.1), we control the transfer factor $\bs\delta(h)$ for certain elements $h$. We use the notational conventions $\wt\Delta$, $\Delta^j$, etc., of \cite{13 (1.1.2)}. 
\proclaim{Transfer Lemma} 
Let $h_0\in \bk J_E$, and suppose that $\ups_F(\det h_0)$ is relatively prime to $n$. The function $u\mapsto \bs\delta(h_0u)$, $u\in J^1_E$, is constant. 
\endproclaim 
\demo{Proof} 
There is an element $x_0$ of $K$ such that $h_0^{n/l}\equiv x_0\pmod{J^1_E}$. For any $h\in h_0J^1_E$, we also have $h^{n/l}\equiv x_0 \pmod{J^1_E}$. Let $L/F$ be a finite Galois extension, containing $K$, all $F$-conjugates of $h$, $n$ distinct $n$-th roots of unity and an $(n/l)$-th root $x$ of $x_0$. An $F$-conjugate of $h$ is of the form $\zeta x u$, for some $\zeta\in \bs\mu_L$ and some $u\in U^1_L$. The quantity $\wt\Delta(h)$ is therefore a product of terms $(\zeta{-}\zeta')$, for various $\zeta\neq\zeta'\in \bs\mu_L$, a power of $x$ and a $1$-unit of $L$. The differences of roots of unity and the element $x$ depend only on the coset $h_0J^1_E$. Thus $\Delta^1$ and $\Delta^2$ are constant on this coset, so the same applies to $\bs\delta = \Delta^2/\Delta^1$. \qed 
\enddemo 
\subhead 
8.6  
\endsubhead 
It is now convenient to split into cases. For $\gamma\in \vG$, the element $g_\gamma$ is determined modulo $G_E$, and $G_E$ is contained in the kernel of the character $\vk\circ\det$ of $G$. Thus $\gamma\mapsto \vk(\det g_\gamma)$ is a character of $\vG$. 
\proclaim{Lemma} 
The character $\gamma\mapsto \vk(\det g_\gamma)$ of $\vG$ is trivial unless the following condition is satisfied: 
$$ 
l = |\vG| = 2 \quad \text{\rm and} \quad \text{\rm  $e/2$ is odd.} 
\tag 8.6.1 
$$ 
Suppose condition \rom{(8.6.1)} holds (and so $p\neq2$). The character is then trivial if $q\equiv 1\pmod 4$, of order $2$ if $q\equiv 3\pmod4$. 
\endproclaim 
\demo{Proof} 
By the Normal Basis Theorem, the automorphism $g_\gamma$ of the $E^\vG$-vector space $E$ is a permutation matrix. Its $E^\vG$-determinant is the signature of this permutation. Since $\vk$ has order $l$, the lemma follows immediately. \qed 
\enddemo 
\subhead 
8.7 
\endsubhead 
In this sub-section, we assume that condition (8.6.1) {\it fails.} Consequently, $\vk(\det g_\gamma) = 1$ for all $\gamma\in \vG$, and (8.4.3) is reduced to  
$$
\roman{tr}^\vk\,\pi_\rho(h) = \eps_{E/F} \sum_{\gamma\in \vG} \roman{tr}\,\pi_{\rho,E}^\gamma(h).  
\tag 8.7.1 
$$
We deduce from (8.4.1) that 
$$
\eps_{E/F}\,\sum_{\gamma\in \vG} \roman{tr}\,\pi_{\rho,E}^\gamma(h) = c_\alpha^{K/F}\, \bs\delta(h) \sum_{\omega\in\vO} \roman{tr}\,\tau_\rho^\omega(h). 
\tag 8.7.2 
$$ 
We apply 8.4 Lemma to simplify the relation (8.7.2). Writing $\vD = \roman{Aut}(E|K)$, we have an exact sequence 
$$
1\to \vD\longrightarrow \vG\longrightarrow \vO. 
$$
Let $\vO_0$ denote the image of $\vG$ in $\vO$: thus either $\vO_0 = \vO$ or $\vO_0$ is trivial. Either way, 8.4 Lemma implies 
$$
\eps_{E/F}\,\sum_{\gamma\in \vG} \roman{tr}\,\pi_{\rho,E}^\gamma(h) = c_\alpha^{K/F}\, \bs\delta(h) \sum_{\omega\in\vO_0} \roman{tr}\,\tau_\rho^\omega(h). 
$$  
Applying 8.1 Proposition, we expand the right hand side to get 
$$ 
\aligned 
\eps_{E/F}\,\sum_{\gamma\in \vG} \roman{tr}\,\pi_{\rho,E}^\gamma(h) &= \eps_{E/K}\, c_\alpha^{K/F}\, \bs\delta(h) \sum_{\omega\in\vO_0} \sum_{\delta\in\vD} \roman{tr}\,\tau_{\rho,E}^{\delta\omega}(h) \\  
&= \eps_{E/K}\, c_\alpha^{K/F}\, \bs\delta(h) \sum_{\gamma\in\vG} \roman{tr}\,\tau_{\rho,E}^\gamma(h), 
\endaligned 
\tag 8.7.3 
$$
the factor $\eps_{E/K}$ being a constant sign, given by (5.3.2), (5.3.3) relative to the base field $K$ in place of $F$. It depends only on the simple stratum $[\frak a_K,\beta]$. 
\par 
As in (8.2.1), (8.2.2), $\nu_\rho$ denotes $\upr L\rho$ and $\pi_{\rho,E} = \phi_\rho\nu_\rho$, for some $\phi_\rho \in X_1(E)$. We have to show that $\phi_\rho$ is independent of the choice of $\rho \in \Gg1E\alpha$. By inductive hypothesis, there is a character $\mu^K = \mu^K_{1,\alpha}\in X_1(E)$ such that $(\tau_\rho)_E = \mu^K\nu_\rho$, for all $\rho \in \Gg1E\alpha$.  
\par
We compare the central characters $\omega_{\pi_\rho}$, $\omega_{\tau_\rho}$ on $F^\times$. In (8.7.3), we replace $h$ by $zh$, with $z\in F^\times$. The left hand side changes by a factor $\omega_{\pi_\rho}(z)$ and the right hand side by $\omega_{\tau_\rho}(z)\bs\delta(zh)/\bs\delta(h)$. Since $\pi_\rho = \roman A_{K/F}\,\tau_\rho$, the quotient $\omega_{\pi_\rho}^{-1}\omega_{\tau_\rho}|_{F^\times}$ is $d_{K/F}^{n/l}$, where $d_{K/F}$ is the determinant of the regular representation $\Ind_{K/F}\,1_K$ of $\scr W_K\backslash \scr W_F$. It is therefore unramified except when $l=2$ and $n/l$ is odd. Since $l$ is the largest prime divisor of $|\vG|$, this case is excluded by 8.6 Lemma and the initial hypothesis of the sub-section. 
\par 
So, $\bs\delta(zh)/\bs\delta(h)$ depends only on $\ups_F(z)$. The characters $\phi_\rho^\gamma$, $(\mu^K)^\gamma$, thus all agree on $\bs\mu_F = \bs\mu_E$. We deduce that $\phi_\rho = \chi_\rho\mu^K$, where $\chi_\rho$ is unramified and, therefore, fixed by $\vG$. Thus (8.7.3) reduces to 
$$
\eps_{E/F}\, \chi_\rho({\det}_E\,h) \sum_{\gamma\in \vG} \roman{tr}\,\tau_{\rho,E}^\gamma(h) = \eps_{E/K}\, c_\alpha^{K/F}\, \bs\delta(h) \sum_{\gamma\in \vG} \roman{tr}\,\tau_{\rho,E}^\gamma(h). 
\tag 8.7.4
$$ 
Let $a$ be an integer relatively prime to $n$. By 8.1 Corollary, we may choose $h\in \bk J_E\cap \LR G$ with $\ups_F(\det h) = \ups_E(\det_Eh) = a$, such that $\sum_\gamma \roman{tr}\,\tau_{\rho,E}^\gamma(h) \neq 0$. We deduce that 
$$
\chi_\rho({\det}_E\, h) = \eps_{E/F}\,\eps_{E/K}\,c_\alpha^{K/F}\,\bs\delta(h), 
\tag 8.7.5 
$$
for all such $h\in \bk J_E\cap \LR G$. The Transfer Lemma of 8.5 implies that (8.7.5) holds for all $h\in \bk J_E$ with $\ups_E(\det_Eh)$ relatively prime to $n$. It follows that $\chi_\rho$ is independent of $\rho$, as desired. 
\subhead 
8.8 
\endsubhead 
We are reduced to the case where (8.6.1) {\it holds,} that is, where $|\vG| = 2$ and $[E^\vG{:}F]$ is odd. The canonical map $\vG \to \vO = \Gal KF$ is therefore an isomorphism. Let $\gamma$ be the non-trivial element of $\vG$. 
\par
The analogue of (8.7.2) reads 
$$
\multline 
\eps_{E/F}\, \big(\roman{tr}\,\pi_{\rho,E}(h) + \vk(-1)\roman{tr}\,\pi_{\rho,E}^\gamma(h) \big) \\ = \eps_{E/K}\,c_\alpha^{K/F}\,\bs\delta(h) \big( \roman{tr}\,\tau_{\rho,E}(h) + \roman{tr}\,\tau_{\rho,E}^\gamma(h)\big). \endmultline
\tag 8.8.1 
$$ 
As before, set $\nu_\rho = \upr L\rho$. There is a character $\phi_\rho\in X_1(E)$ such that $\pi_{\rho,E} = \phi_\rho\nu_\rho$. By inductive hypothesis, there is a character $\mu^K = \mu^K_{1,\alpha}\in X_1(E)$ such that $\tau_{\rho,E} = \mu^K\nu_\rho$ for all choices of $\rho\in \Aa1E\alpha$. We set $\chi_\rho = \phi_\rho/\mu^K$. 
\par 
Comparing the central characters of $\pi_\rho$ and $\tau_\rho$ on $F^\times$, we find: 
\proclaim{Lemma 1} 
\roster 
\item 
If $z\in F^\times$, then $\bs\delta(zh)/\bs\delta(h) = \vk(z)$ and 
\item 
$\chi_\rho|_{\bs\mu_F} = \vk|_{\bs\mu_F}$, this character having order $2$. 
\endroster 
\endproclaim 
We now invoke an elementary point. 
\proclaim{Lemma 2} 
Let $\xi\in X_1(E)$ and $y\in E^\times$. If $\xi|_{\bs\mu_F}$ has order $2$ and $\ups_E(y)$ is odd, then $\xi^\gamma(y) = \vk(-1)\xi(y)$. 
\endproclaim 
Taking the lemmas into account, the relation (8.8.1) reads 
$$ 
\multline 
\eps_{E/F}\, \chi_\rho({\det}_E\,h) \big(\roman{tr}\,\tau_{\rho,E}(h) + \roman{tr}\,\tau_{\rho,E}^\gamma(h) \big) 
\\ = c_\alpha^{K/F}\,\bs\delta(h) \big( \roman{tr}\,\tau_{\rho,E}(h) + \roman{tr}\,\tau_{\rho,E}^\gamma(h)\big). 
\endmultline 
$$
That is, 
$$
\chi_\rho({\det}_E\, h) = \eps_{E/K}\eps_{E/F}\,c_\alpha^{K/F}\,\bs\delta(h), 
\tag 8.8.2 
$$
for all $h\in \bk J_E\cap \LR G$ with $\ups_E(\det_E h)$ odd and $\roman{tr}\,\tau_{\rho,E}(h) + \roman{tr}\,\tau_{\rho,E}^\gamma(h) \neq 0$. Arguing as before, the character $\chi_\rho$ is independent of $\rho$. 
\par 
This completes the proof of 8.2 Proposition, and hence that of the Comparison Theorem in the present case. \qed 
\subhead 
8.9 
\endsubhead 
For ease of reference, we display an intermediate conclusion of the preceding argument. We use the notation of the Comparison Theorem in 7.3. 
\proclaim{Corollary} 
Suppose that $E/F$ is totally ramified and that $\vG = \roman{Aut}(E|F)$ is non-trivial. Let $l$ be the largest prime divisor of $|\vG|$ and let $K/F$ be the unique sub-extension of $E/F$ of order $l$. There is a unique character $\chi^K\in X_1(E)$ such that 
$$
\chi^K({\det}_E\,h) = \eps_{E/F}\eps_{E/K}\,c_\alpha^{K/F}\,\bs\delta_{K/F}(h), 
\tag 8.9.1 
$$
for all $h\in \bk J_E$ with $\ups_E(\det_Eh)$ relatively prime to $n$. The character $\mu^F_{1,\alpha}$ then satisfies 
$$
\mu^F_{1,\alpha} = \chi^K\cdot\mu^K_{1,\alpha}. 
\tag 8.9.2 
$$
\endproclaim 
We recall that the factors $\eps_{E/F}$, $\eps_{E/K}$ are symplectic signs: for example, $\eps_{E/F}$ is $t_{E^\times/F^\times U^1_E}(J^1_\theta/H^1_\theta) = \pm1$. 
In the essentially tame case, one may compute the character $\mu^F_{1,\alpha}$ iteratively, using (8.9.2). The first step, where $\vG = \roman{Aut}(E|F)$ is trivial, requires a base change argument as in \cite{10 (4.6)}. The steps coming from automorphic induction are computed in \cite{10 (4.5)} and \cite{11} {\it passim.} 
\par 
We comment briefly on the general case. For the first step, as in 8.3, a base change argument analogous to that of \cite{10} could not determine the desired character directly. In the cyclic case of the corollary, one may look at values $\chi^K(h_1h_2^{-1})$, for elements $h_i$ of $\bk J_E$ with $\ups_E(\det_E h_i)$ relatively prime to $n$. Suppose, for a simple example, that $l$ and $p$ are odd. The discriminant character $d_{K/F}$ is then trivial, so $\chi^K\circ\det_E$ is then trivial on $F^\times$. Since $p$ is odd, one may choose the $h_i$ so that $\ups_E(\det_E h_1h_2^{-1})$ is divisible by $e$ but relatively prime to $p$. The transfer factors $\bs\delta_{K/F}(h_i)$ are each the product of an $l$-th root of unity and a positive real number. One concludes that $\chi^K$ has order dividing $l$ and that $\eps_F\eps_Kc_\alpha>0$. An argument similar to that of \cite{10 (4.5 Lemma)} now implies that $\chi^K$ is trivial. 
\head\Rm 
9. Unramified automorphic induction 
\endhead 
We start the main part of the proof of the Comparison Theorem. We work in the general case, using the notation laid out at the start of 5.4. In particular, $E_0/F$ is a tame parameter field for $\theta$, $E/E_0$ is unramified of degree $m$ with $\vD = \Gal E{E_0}$, and $K/F$ is the maximal unramified sub-extension of $E/F$. In that situation, we set $\vT = \text{\it cl\/}(\theta)$ and use the notation $\vT_K$ etc., as in (5.8.1). As usual, we abbreviate $J^k = J^k_\theta$, $J^k_K = J^k_{\theta_K} = J^k\cap G_K$, and so on. In addition, we set $d=[K{:}F]$. 
\subhead 
9.1 
\endsubhead  
We recall a property of automorphic induction. 
\proclaim{Proposition} 
The operation of {\it automorphic induction\/} induces a bijection 
$$
\roman A_{K/F}: \vD\backslash \Aa1K{\vT_K}^{\vD\text{\rm -reg}} @>{\ \ \approx\ \ }>> \Aa mF\vT. 
\tag 9.1.1 
$$ 
\endproclaim 
\demo{Proof} 
The endo-classes $\vT_K^\gamma$, $\gamma \in \vD\backslash \vG$, are distinct, as follows from (2.3.5). Consequently, any representation $\tau \in \Aa1K{\vT_K}^{\vD \text{\rm -reg}}$ is $\vG$-regular. The automorphically induced representation $\pi = \roman A_{K/F}\,\tau$ is therefore cuspidal. The endo-class $\vt(\pi)$ is $\frak i_{K/F}\,\vT_K = \vT$ \cite{8 (7.1 Corollary)}, whence $\pi \in \Aa mF\vT$. 
\par
Conversely, let $\pi \in \Aa mF\vT$. The integer $t(\pi)$ is $d = [K{:}F]$ (6.3 Lemma 2). It follows that $\pi \cong \roman A_{K/F}\,\tau$, for some $\vG$-regular $\tau \in \Ao{n/d}K$. The representation $\tau$ is determined by $\pi$, up to $\vG$-conjugation. The endo-class $\vt(\tau)$ is a $K/F$-lift of $\vt(\pi) = \vT$. The set of $K/F$-lifts of $\vT$ form a single $\vG$-orbit, so we may choose $\tau$ to satisfy $\vt(\tau) = \vT_K$. Therefore $\tau\in \Aa1K{\vT_K}$. The representation $\tau$ is surely $\vD$-regular, and $\pi$ determines its $\vD$-orbit. \qed 
\enddemo 
We compare the bijection (9.1.1) with the bijection (5.8.2) given by the algebraic induction map $\text{\it ind}_{K/F}$. As before, we write $\chi_E = \chi\circ \N E{E_0}$, $\chi\in X_1(E_0)$, and recall that the map $\chi\mapsto \chi_E$ induces an isomorphism $X_1(E_0)/X_0(E_0)_m \to X_1(E)^\vD$. 
\proclaim{Unramified Induction Theorem} 
There is a unique character $\nu \in X_1(E)^\vD$ such that 
$$
\roman A_{K/F}\,\tau = \text{\it ind}_{K/F}\,\nu\odot_{\vT_E}\tau, 
\tag 9.1.2  
$$
for all $\tau\in \Aa1K{\vT_K}^{\vD\text{\rm -reg}}$. In particular, 
$$
\roman A_{K/F}(\chi_E\odot_{\vT_E}\tau) = \chi\odot_{\vT_{E_0}} \roman A_{K/F}\,\tau, 
\tag 9.1.3 
$$
for $\chi\in X_1(E_0)$ and $\tau \in \Aa1K{\vT_K}^{\vD\text{\rm -reg}}$. 
\endproclaim 
The proof will occupy us until the end of \S11. Before starting it, we remark that the second assertion of the theorem follows from the first and (5.8.3). The character $\nu$ is of the form $(\nu_0)_E$, for some $\nu_0\in X_1(E_0)$. In those terms, the relation (9.1.2) reads 
$$
\roman A_{K/F}\,\tau = \nu_0\odot_{\vT_{E_0}} \text{\it ind}_{K/F}\,\tau. 
\tag 9.1.4 
$$
While (9.1.2) determines $\nu$ uniquely, the character $\nu_0$ is only determined modulo $X_0(E_0)_m$. 
\subhead 
9.2 
\endsubhead 
We choose a prime element $\vp_F$ of $F$ and a representation $\kappa \in \scr H(\theta)$ satisfying the conditions of 5.6 Lemma 2. We define $\kappa_K = \ell_{K/F}(\kappa) \in \scr H(\theta_K)$, as in 5.6 Proposition. If $\xi\in X_1(E)^{\vD \text{\rm -reg}}$, we set 
$$ 
\tau_\xi = \cind_{\bk J_K}^{G_K}\, \xi \odot\kappa_K \in \Aa1K{\vT_K}^{\vD\text{\rm -reg}}. 
\tag 9.2.1 
$$
Lemma 1 of 5.7 implies that every element of $\Aa1K{\vT_K}^{\vD\text{\rm -reg}}$ is of this form, for some $\xi\in X_1(E)^{\vD\text{\rm -reg}}$. Likewise, we define 
$$
\pi_\xi = \cind_{\bk J}^G\, \lambda_\xi \ltimes \kappa = \text{\it ind}_{K/F}\,\tau_\xi \in \Aa mF\vT, 
\tag 9.2.2 
$$ 
in the notation of (3.6.2). 
\proclaim{Lemma} 
There exists $\chi = \chi_\xi\in X_1(E)^{\vD\text{\rm -reg}}$ such that 
$$
\roman A_{K/F}\,\tau_\xi = \pi_\chi. 
\tag 9.2.3  
$$ 
The correspondence $\xi\mapsto \chi_\xi$ induces a bijection of $\vD\backslash X_1(E)^{\vD\text{\rm -reg}}$ with itself. 
\endproclaim 
\demo{Proof} 
This follows from the bijectivity of the map (9.1.1) and 5.8 Proposition. \qed 
\enddemo 
In these terms, the Unramified Induction Theorem of 9.1 is equivalent to the following. 
\proclaim{Theorem} 
There is a unique character $\mu = \mu^{K/F}_\alpha\in X_1(E)^\vD$ such that 
$$
\roman A_{K/F}\,\tau_\xi = \pi_{\mu\xi} = \text{\it ind}_{K/F}\,\mu\odot_{\vT_E}\tau_\xi, 
$$
for all $\xi\in X_1(E)^{\vD \text{\rm -reg}}$. 
\endproclaim 
To prove the theorem, we use the automorphic induction equation, 
$$
\roman{tr}^\vk\,\pi_\chi(g) = c^{K/F}_\alpha\,\bs\delta(g) \sum_{\gamma\in \vG} \roman{tr}\,\tau_\xi^\gamma(g), \quad g\in G_K\cap \LR G,
\tag 9.2.4 
$$ 
as in \S8. In this case, $\vk$ is an unramified character of $F^\times$ of order $d = [K{:}F]$. It is trivial on the group $\det\bk J$, so we may again apply the Uniform Induction Theorem of \cite{13 (1.3)}. The $\vk$-twisted trace $\roman{tr}^\vk\,\pi_\chi$ is defined so as to satisfy the Mackey formula 
$$
\roman{tr}^\vk\,\pi(g) = \sum_{x\in G/\bk J} \vk(\det x^{-1})\,\roman{tr}\,\vL_\chi(x^{-1}gx), \quad g\in G_K \cap \LR G,
\tag 9.2.5 
$$
where $\vL_\chi = \lambda_\chi^\bk J \otimes \kappa =  \lambda_\chi \ltimes \kappa$. The constant $c_\alpha^{K/F}$ depends only on $\theta$ and $K$, and $\bs\delta = \bs\delta_{K/F}$ is a transfer factor. We recall that there are only finitely many terms in the expansion (9.2.5) \cite{14 (1.2)}. 
\subhead 
9.3 
\endsubhead 
We evaluate each side of the relation (9.2.4) at a special family of elements. 
\par 
We use the groups $\pre p{\bk J}(\vp_F) = \pre p{\bk J}$ and $\pre p{\bk J_K}$ defined in 5.5. Recalling that $p^r = [P{:}E]$, we say that $g\in G$ is {\it pro-unipotent in\/} $G/\ags{\vp_F}$ if $g^{p^{r+t}}\to 1$ in $G/\ags{\vp_F}$ as $t\to\infty$. In particular, any element of $\pre p{\bk J}$ is pro-unipotent in $G/\ags{\vp_F}$. 
\proclaim{Lemma 1} 
Let $g = \zeta h$, where $\zeta\in \bs\mu_K$ is $\vG$-regular and $h\in \pre p{\bk J_K}\cap (G_K\LR)$. The element $g$ then lies in $G_K\cap \LR G$. Let $x\in G$ satisfy $x^{-1}gx\in \bk J$. The coset $x\bk J$ then contains an element $y$ of $N_K = N_G(K^\times)$. For any such $y$, we have $h' = y^{-1}hy\in \pre p{\bk J_K}$. 
\endproclaim 
\demo{Proof} 
The first assertion is immediate. 
\par 
Of the element $g = \zeta h$, the factors $\zeta$, $h$ commute with each other; the first has finite order prime to $p$, and the second is pro-unipotent in $G/\ags{\vp_F}$. The factors $\zeta$, $h$ are thereby uniquely determined. So, any closed subgroup of $G$ containing $g$ and $\vp_F$ must contain both $\zeta$ and $h$. If $x\in G$ and if $g' = x^{-1}gx$ lies in $\bk J$, it follows that both $\zeta' = x^{-1}\zeta x$ and $h' = x^{-1}hx$ lie in $\bk J$. In particular, $\zeta'\in J^0$.  The algebra $F[\zeta']$ is a field in which $\zeta'$ is a root of unity of the same order as $\zeta$. It follows that $\zeta'$ is $J^0$-conjugate to an element of $\bs\mu_K$. Replacing $x$ by $y = xj$, for suitable $j\in J^0$, we therefore get $g'' = y^{-1}gy\in \bk J$ and $\zeta'' = y^{-1}\zeta y\in \bs\mu_K$. Thus $y\in N_K$. Again, $h'' = y^{-1}hy$ lies in $\bk J$ and commutes with $\zeta''$. It follows that $h''\in \bk J_K = P^\times J^1_K$. Since $h''$ is pro-unipotent in $G/\ags{\vp_F}$, it must lie in $\pre p{\bk J_K}$ as required.  \qed 
\enddemo 
We extract a remark made in the course of the proof. 
\proclaim{Lemma 2} 
Let $h\in \pre p{\bk J_K}$ and $y\in G_K$. If the element $y^{-1}hy$ lies in $\bk J_K$ then it lies in $\pre p{\bk J_K}$. 
\endproclaim 
The natural action of $N_K$ on $K^\times$ induces an isomorphism $N_K/G_K \to \vG$. For each $\gamma\in \vG$, we choose $g_\gamma\in N_K$ with image $\gamma$ in $\vG$. The groups $N_K\cap \bk J$, $N_K\cap  J^0$ both have image $\vD$ in $\vG$. The map $x\mapsto x\bk J$ thus induces a bijection 
$$
\bigcup_{\gamma \in \vG/\vD} g_\gamma G_K/\bk J_K \longrightarrow N_K\bk J/\bk J, 
$$ 
the union being disjoint. The character $\vk\circ\det$ of $G$ is trivial on $G_K$, while $\det g_\gamma \equiv \pm1 \pmod{\det G_K}$ (as in the discussion of 8.6). Since $\vk $ is unramified, this implies  $\vk(\det x) = 1$ for all $x\in N_K$. Thus (9.2.5) is reduced to 
$$
\roman{tr}^\vk\,\pi_\chi(\zeta h) = \sum_{\gamma\in \vG/\vD}\,\, \sum_{y\in G_K/\bk J_K}  \roman{tr}\, \vL_\chi(y^{-1}\zeta^\gamma h^\gamma y) = \roman{tr}\,\pi_\chi(\zeta h),  
$$
for any element $g = \zeta h$ as in Lemma 1. 
\par 
We recall that $\vL_\chi = \lambda_\chi^\bk J\otimes\kappa$, as in (9.2.5). We evaluate $\roman{tr}\, \kappa(y^{-1} \zeta^\gamma h^\gamma y)$, using (5.6.1). The character $\eps^1_{\bs\mu_K} = t^1_{\bs\mu_K}(J^1_\theta/H^1_\theta)$ of the cyclic group $\bs\mu_K$ has order at most $2$ and so is $\vG$-invariant. Therefore 
$$
\roman{tr}\,\kappa(y^{-1}\zeta^\gamma h^\gamma y) = \eps^0_{\bs\mu_K} \eps^1_{\bs\mu_K}(\zeta)\, \roman{tr}\,\kappa_K(y^{-1} h^\gamma y),  
$$ 
where $\kappa_K = \ell_{K/F}(\kappa)$ (as in the definition of $\tau_\xi$) and $\eps^0_{\bs\mu_K} = t^0_{\bs\mu_K}(J^1_\theta/H^1_\theta) = \pm1$. 
\par 
We consider the term $\roman{tr}\,\lambda_\chi^\bk J(\zeta^\gamma y^{-1}h^\gamma y)$. As in 3.1 (but with base field $K$), we have the character $\chi^{\bk J_K} \in X_1(\theta_K)$ satisfying 
$$ 
\chi^{\bk J_K}|_{P^\times} = \chi_P = \chi\circ \N PE. 
$$ 
The group $\pre p{\bk J_K}$ is contained in the subgroup $P_0^\times J^1_K$ of $\bk J_K = P^\times J^1_K$. The restriction of $\lambda_\chi^\bk J$ to $P_0^\times J^1_K$ is a multiple of the character $\chi^{\bk J_K}|_{P^\times_0 J^1_K}$. Let $X_0(\theta_K)$ denote the subgroup of $X_1(\theta_K)$ consisting of characters trivial on $J^0_K = \bs\mu_KJ^1_K$: thus $X_0(\theta_K)$ is the image of $X_0(E)$ under the map $\phi\mapsto \phi^{\bk J_K}$ of (3.1.1). 
\par 
We choose a character $\tilde a_\chi \in X_0(\theta_K)$ to agree with $\chi^{\bk J_K}$ on $\pre p{\bk J_K}$. This is of the form $\tilde a_\chi = a_\chi^{\bk J_K}$, for a character $a_\chi\in X_0(E)$, uniquely determined modulo $X_0(E)_e$, $e=e(E|F) = [E{:}K]$. We get 
$$ 
\multline 
\roman{tr}\,\lambda_\chi^\bk J(\zeta^\gamma y^{-1}h^\gamma y) \\ = \roman{tr}\,\lambda_\chi^\bk J(\zeta^\gamma)\,\tilde a_\chi(y^{-1} h^\gamma y) = (-1)^{m-1} \tilde a_\chi(y^{-1} h^\gamma y) \sum_{\delta\in\vD} \chi_P(\zeta^{\gamma\delta}), 
\endmultline 
$$ 
and so 
$$ 
\roman{tr}\,\vL_\chi(\zeta^\gamma y^{-1}h^\gamma y) = (-1)^{m-1} \eps^0_{\bs\mu_K}\eps^1_{\bs\mu_K} (\zeta) \, \roman{tr}(a_\chi\odot \kappa_K)(y^{-1}h^\gamma y) \sum_{\delta\in \vD} \chi_P(\zeta^{\gamma\delta}). 
$$ 
By definition, $\kappa_K \in \scr H(\theta_K)^\vD = \scr T(\theta_K)^\vD$, so 
$$ 
\rho = \cind_{\bk J_K}^{G_K}\,\kappa_K 
\tag 9.3.1 
$$ 
is an irreducible, cuspidal representation of $G_K$, lying in $\Aa1K{\vT_K}^\vD$. Assembling the preceding identities, and taking account of Lemma 2, we obtain 
$$ 
\aligned 
\roman{tr}^\vk\,\pi_\chi(\zeta h) 
&= (-1)^{m-1} \eps^0_{\bs\mu_K}\eps^1_{\bs\mu_K}(\zeta)\,\sum_{\gamma\in \vG/\vD} \bigg(\sum_{\delta\in \vD} \chi_P(\zeta^{\gamma\delta})\bigg) \,\roman{tr}(a_\chi\odot_{\vT_E} \rho)(h^\gamma) \\ 
&= (-1)^{m-1} \eps^0_{\bs\mu_K}\,\sum_{\gamma\in \vG/\vD} \bigg(\sum_{\delta\in \vD} \eps^1_{\bs\mu_K}\chi_P(\zeta^{\gamma\delta})\bigg)\, \roman{tr}(a_\chi\odot_{\vT_E} \rho)(h^\gamma). 
\endaligned 
\tag 9.3.2 
$$ 
\subhead 
9.4 
\endsubhead 
We look to the other side of the automorphic induction relation (9.2.4), using the same notation. We expand $\roman{tr}\,\tau_\xi(\zeta^\gamma h^\gamma)$. The representation $\kappa_K$ is trivial on $\bs\mu_K$ (remark following 5.6 Lemma 3). From 9.3 Lemma 2, we obtain 
$$ 
\align 
\roman{tr}\,\tau_\xi(\zeta ^\gamma h^\gamma) &= \xi_P(\zeta^\gamma) \sum_{y\in G_K/\bk J_K} \tilde a_\xi(y^{-1}h^\gamma y)\,\roman{tr}\, \kappa_K(y^{-1}h^\gamma y) \\ &= \xi_P(\zeta^\gamma)\, \roman{tr}(a_\xi \odot_{\vT_E} \rho)(h^\gamma). \endalign 
$$ 
Taking account of (9.3.2), the relation (9.2.4) now reads  
$$ 
\multline 
(-1)^{m-1} \eps^0_{\bs\mu_K}\,\sum_{\gamma\in \vG/\vD} \bigg(\sum_{\delta\in \vD} \eps^1_{\bs\mu_K}\chi_P(\zeta^{\gamma\delta})\bigg) \roman{tr}(a_\chi\odot_{\vT_E} \rho)(h^\gamma) \\ = 
c_\alpha^{K/F}\,\bs\delta(\zeta h) \,\sum_{\gamma\in \vG/\vD} \bigg(\sum_{\delta\in \vD} \xi_P(\zeta^{\gamma\delta})\bigg) \roman{tr}(a_\xi\odot_{\vT_E} \rho)(h^\gamma) . 
\endmultline \tag 9.4.1  
$$ 
\subhead 
9.5 
\endsubhead 
We simplify the relation (9.4.1), starting with a special case. Since $h\in \pre p{\bk J_K}$ and $\zeta\in \bs\mu_K$ is $\vG$-regular, the condition $\zeta h\in \LR G$ is equivalent to $h\in \LR{(G_K)}$. 
\proclaim{Transfer Lemma} 
If $h\in J^1_K\cap \LR{(G_K)}$ and if $\zeta\in \bs\mu_K$ is $\vG$-regular, then $\bs\delta(\zeta h) = 1$. 
\endproclaim 
\demo{Proof} 
Let $\beta,\gamma\in \vG$, $\beta\neq\gamma$. Let $y_\beta$, $y_\gamma$ be eigenvalues of $h^\beta$, $h^\gamma$ respectively, in some splitting field $L/K$. Thus $y_\beta, y_\gamma\in U^1_L$ and $\zeta^\beta y_\beta{-}\zeta^\gamma y_\gamma$ is a unit in $L$. We follow the standard notation laid out in \cite{13 (1.1.2)} and observe the convention of \cite{13, Remark 1}. In those terms, the quantity $\wt\Delta(\zeta h)$ is a unit and, since $\vk$ is unramified, the lemma follows. \qed 
\enddemo 
The identity (9.4.1) now implies 
$$ 
\multline 
(-1)^{m-1} \eps^0_{\bs\mu_K}\,\sum_{\gamma\in \vG/\vD} \bigg(\sum_{\delta\in \vD} \eps^1_{\bs\mu_K}\chi_P(\zeta^{\gamma\delta})\bigg)\, \roman{tr}(a_\chi\odot_{\vT_E} \rho)(h^\gamma) \\ = 
c_\alpha^{K/F}\,\sum_{\gamma\in \vG/\vD} \bigg(\sum_{\delta\in \vD} \xi_P(\zeta^{\gamma\delta})\bigg)\, \roman{tr}(a_\xi\odot_{\vT_E} \rho)(h^\gamma) , 
\endmultline \tag 9.5.1  
$$ 
for $h\in J^1_K\cap \LR{(G_K)}$. 
The characters $a_\chi$ and $a_\xi$ are unramified, so we get 
$$ 
\roman{tr}(a_\chi\odot_{\vT_E} \rho)(h^\gamma) = \roman{tr}\,\rho(h^\gamma) = \roman{tr}(a_\xi\odot_{\vT_E} \rho)(h^\gamma), 
$$ 
and the relation (9.5.1) becomes  
$$ 
\multline 
(-1)^{m-1} \eps^0_{\bs\mu_K}\,\sum_{\gamma\in \vG/\vD} \bigg(\sum_{\delta\in \vD} \eps^1_{\bs\mu_K}\chi_P(\zeta^{\gamma\delta})\bigg) \roman{tr}\,\rho(h^\gamma) \\ = 
c_\alpha^{K/F}\,\sum_{\gamma\in \vG/\vD} \bigg(\sum_{\delta\in \vD} \xi_P(\zeta^{\gamma\delta})\bigg) \roman{tr}\,\rho(h^\gamma), \quad h\in J^1_K\cap \LR{(G_K)}. 
\endmultline \tag 9.5.2 
$$ 
To proceed further, we need: 
\proclaim{Linear Independence Lemma} 
Let $h\in \bk J_K$, and let $u$ range over the set of elements of $J^1_K$ such that $hu\in \LR{(G_K)}$. The set of functions 
$$
u\longmapsto \roman{tr}\,\rho^\gamma(hu),\quad \gamma\in \vD\backslash\vG, 
$$
is linearly independent. 
\endproclaim 
We shall prove this lemma in 9.7 below. It implies 
$$ 
(-1)^{m-1} \eps^0_{\bs\mu_K}\,\sum_{\delta\in \vD} \eps^1_{\bs\mu_K}\chi_P(\zeta^\delta) = 
c_\alpha^{K/F}\,\sum_{\delta\in \vD} \xi_P(\zeta^\delta) . 
$$ 
We recall that the constant $c_\alpha^{K/F}$ does not depend on $\xi$. The lemma of 9.2 therefore enables us to apply \cite{13 (2.3 Corollary 1)} and so obtain: 
\proclaim{Proposition} 
The characters $\xi$ and $\eps^1_{\bs\mu_K}\chi$ of $\bs\mu_K$ lie in the same $\vD$-orbit, and 
$$ 
c_\alpha^{K/F} = (-1)^{m-1}\,\eps^0_{\bs\mu_K}. 
$$ 
\endproclaim 
The first conclusion of the proposition may be strengthened. 
\proclaim{Corollary} 
Let $\xi\in X_1(E)^{\vD \text{\rm -reg}}$. There exists a unique character 
$$ 
\chi = \chi_\xi \in X_1(E)^{\vD\text{\rm -reg}} 
$$ 
such that 
$$
\chi|_{U_E} = \eps^1_{\bs\mu_K}\xi|_{U_E} \quad \text{\rm and} \quad 
\roman A_{K/F}\,\tau_\xi = \pi_\chi. 
$$ 
The character $\xi^{-1}\chi_\xi$ has finite order. 
\endproclaim 
\demo{Proof} 
Given $\xi$, there surely exists $\chi\in X_1(E)^{\vD\text{\rm -reg}}$ satisfying the two requirements. The first of them determines $\chi$ modulo $X_0(E)$. We therefore take $\phi \in X_0(E)$, $\phi\neq 1$, and show that $\pi_{\phi\chi}\neq \pi_\chi$. 
\par
The character $\phi$ is of the form $\psi\circ\N EF$, for some $\psi\in X_0(F)$. We then have $\pi_{\phi\chi} = \psi\pi_\chi$. The relation $\psi\pi_\chi = \pi_\chi$ would imply $\psi\circ\N KF = 1$, which is false since $\psi\circ \N EF \neq 1$. Thus $\phi=1$ as required. 
\par
For the final assertion, it is enough to show that $\chi^{-1}\xi(\vp_F)$ has finite order. Let $\psi_2$ temporarily denote the unramified character of $F^\times$ of order $2$. By construction, both $\kappa$ and $\kappa_K$ are trivial on $\vp_F$. The central character $\omega_{\pi_\chi}$ of $\pi_\chi$ thus satisfies $\omega_{\pi_\chi}(\vp_F) = \chi_P(\vp_F)$. On the other hand, the central character of $\roman A_{K/F}\,\tau_\xi$ is $\psi_2^{(n-d)}\,\omega_{\tau_\xi}|_{F^\times}$. Its value at $\vp_F$ is therefore $(-1)^{(n-d)}\xi_P(\vp_F)$, whence $\chi^{-1}\xi(\vp_F)^{p^r} = (-1)^{(n-d)}$. \qed 
\enddemo 
In the situation of the corollary, we define 
$$
\mu_\xi = \xi^{-1}\chi_\xi, \quad \xi\in X^1(E)^{\vD \text{\rm -reg}}. 
\tag 9.5.3 
$$ 
We show that $\mu_\xi$ does not depend on $\xi$: this will prove 9.2 Theorem, with $\mu^{K/F}_\alpha$ being the common value of the $\mu_\xi$. 
\subhead 
9.6 
\endsubhead 
We return to the identity (9.4.1), relative to an arbitrary element $h$ of $\pre p{\bk J_K} \cap \LR{(G_K)}$. The character $a_\chi$, by definition, lies in $X_0(E)$. So, there exists an unramified character $b_\chi$ of $K^\times$ such that $a_\chi = b_\chi\circ\N EK$. This gives 
$$
\roman{tr}(a_\chi\odot_{\vT_E}\rho)(h^\gamma) = b_\chi({\det}_K\,h^\gamma)\,\roman{tr}\,\rho(h^\gamma) = b_\chi({\det}_K\,h)\,\roman{tr}\,\rho(h^\gamma). 
$$ 
Looking back at the definition of $a_\chi$ in 9.3, we see that $b_\chi(\det_K h) = \chi^{\bk J_K}(h)$. Similar considerations apply to the character $a_\xi$. Taking account of 9.5 Proposition, the relation (9.4.1) is therefore reduced to 
$$
\chi^{\bk J_K}(h) \sum_{\gamma\in \vG/\vD} \roman{tr}\,\rho(h^\gamma) = \bs\delta(\zeta h)\, \xi^{\bk J_K}(h) \sum_{\gamma\in \vG/\vD} \roman{tr}\,\rho(h^\gamma). 
\tag 9.6.1 
$$ 
This relation is valid for all $h\in \pre p{\bk J_K}\cap \LR{(G_K)}$. The characters $\chi$, $\xi$ therefore satisfy 
$$
\chi^{\bk J_K}(h) = \bs\delta(\zeta h)\,\xi^{\bk J_K}(h), 
\tag 9.6.2 
$$
for any $h\in \pre p{\bk J_K}$ such that $\sum_\gamma\roman{tr}\,\rho(h^\gamma) \neq 0$. For any such $h$, 9.5 Corollary asserts that $\bs\delta(\zeta h)$ is a root of unity. 
\proclaim{Transfer Lemma} 
Let $h\in \pre p{\bk J_K}\cap \LR{(G_K)}$ and let $\zeta \in \bs\mu_K$ be $\vG$-regular. Let 
$$ 
\roman v_K(h) = \ups_K({\det}_K\,h) = \ups_F(\det h)/d. 
$$
Let $\vp_F$ be a prime element of $F$, and write $s = [P{:}K] = ep^r$. If $\bs\delta(\zeta h)$ is a root of unity, then 
$$
\bs\delta(\zeta h) = \vk(\det h)^{s(d-1)/2} = \vk(\vp_F)^{\roman v_K(h)n(d-1)/2}. 
$$
\endproclaim 
\demo{Proof} 
We use the standard notation for transfer factors, as summarized in \cite{13, \S1}. The relevant relations are 
$$
\aligned 
\Delta^2(\zeta h) &= \vk(e_s\wt\Delta(\zeta h)), \\ 
\Delta^1(\zeta h) &= \|\wt\Delta(\zeta h)^2\|_F^{1/2}\,\|\det\zeta h\|_F^{s(1-d)/2}, \\ 
\bs\delta(\zeta h) &= \Delta^2(\zeta h)/\Delta^1(\zeta h). 
\endaligned 
\tag 9.6.3 
$$
Here, $e_s$ is a unit in $K$ such that $e_s^\gamma = (-1)^{s(d-1)}e_s$, for a generator $\gamma$ of $\vG$. We do not need to give the definition of the function $\wt\Delta$: we recall only that $\wt\Delta(\zeta h) \in K^\times$ satisfies $\wt\Delta(\zeta h)^2\in F^\times$ and $e_s\wt\Delta(\zeta h) \in F^\times$. 
\par
If $\bs\delta(\zeta h)$ is a root of unity, then $\Delta^1(\zeta h) = 1$ and so  
$$
\ups_F(e_s\wt\Delta(\zeta h)) = \ups_K(\wt\Delta(\zeta h)) = s(d{-}1)\,\ups_F(\det h)/2, 
$$
whence the lemma follows. \qed 
\enddemo 
\proclaim{Proposition} 
The character $\mu_\xi$ of \rom{(9.5.3)} satisfies 
$$
\mu_\xi^{\bk J_K}(h) = \vk(\det h^{ep^r(d-1)/2}) = \vk(\vp_F)^{\roman v_K(h)n(d-1)/2}.  
\tag 9.6.4
$$
\endproclaim  
\demo{Proof}
Let $h\in \pre p{\bk J_K}$. The Linear Independence Lemma of 9.5 gives $u\in J^1_K$ such that $hu\in \LR{(G_K)}$ and $\sum_\gamma \roman{tr}\,\rho((hu)^\gamma) \neq 0$. The relation (9.6.1) therefore implies $\mu_\xi^{\bk J_K}(hu) = \bs\delta(\zeta hu)$, and this quantity is a root of unity, by 9.5 Corollary. The result now follows from the Transfer Lemma. \qed 
\enddemo 
The proposition determines $\mu_\xi^{\bk J_K}$ on $\pre p{\bk J_K}$. It therefore determines $\mu_\xi$ on the group $\det_E(G_E\cap \pre p{\bk J_K})$, which contains $\vp_F$. We deduce: 
\proclaim{Corollary} 
Let $\xi\in X_1(E)^{\vD \text{\rm -reg}}$. The character $\mu_\xi = \xi^{-1}\chi_\xi$ of \rom{(9.5.3)} satisfies 
$$
\aligned 
\mu_\xi|_{U^1_E} &= 1, \\ 
\mu_\xi|_{\bs\mu_E} &= \eps^1_{\bs\mu_E}, \\ 
\mu_\xi(\vp_F) &= \vk(\vp_F)^{n(d-1)/2} = \pm1, 
\endaligned 
\tag 9.6.5 
$$ 
for a prime element $\vp_F$ of $F$. In particular, if $\xi'\in X_1(E)^{\vD \text{\rm -reg}}$, then 
$$ 
\mu_{\xi'}(z) = \mu_\xi(z), \quad z\in F^\times U_E. 
$$
In particular, if $E_0/F$ is unramified, then $\mu_\xi = \mu_{\xi'}$. 
\endproclaim 
\demo{Proof} 
Only the third property requires comment. From its definition, the character $\eps^1_{\bs\mu_E}$ is trivial on $\bs\mu_F$, so the asserted property is independent of the choice of $\vp_F$. We therefore assume $\vp_F$ is the prime used to define $\pre p{\bk J}$, that is, $\pre p{\bk J} = \pre p{\bk J}(\vp_F)$. The infinite cyclic group $\pre p{\bk J_K}/J^1_K$ is therefore generated by an element $h_0$ of  $P^\times$ such that $h_0^{p^r} \equiv \vp_F \pmod{U^1_P}$. In the notation of the proposition, $\roman v_K(h_0) = e$. On the other hand, $\det_E(h_0) = \N PE(h_0) \equiv h_0^{p^r} \pmod{U^1_P}$, since the field extension $P/E$ is totally wildly ramified. This yields 
$$
\mu_\xi(\vp_F) = \mu_\xi^{\bk J_K}(h_0) = \vk(\vp_F)^{e^2p^rd(d-1)/2} = \vk(\vp_F)^{n(d-1)/2}, 
$$ 
as required. \qed 
\enddemo 
The corollary completes the proof of the Unramified Induction Theorem in the case where $E_0/F$ is unramified. 
\par 
Returning to the general case, the character $\eps_{\bs\mu_K}^1$ is trivial when $p=2$, of order $\le 2$ otherwise. Consequently, there is no need to distinguish between $\eps^1_{\bs\mu_K}$ and $\eps^1_{\bs\mu_K,P}$ as characters of $\bs\mu_K = \bs\mu_E = \bs\mu_P$. 
\remark{Remark} 
In (9.6.5), the character $\eps^1_F(\bs\mu_E)$ is invariably trivial on $\bs\mu_F$. Consequently, the formula for $\vp_F$ holds for all prime elements $\vp_F$ of $F$, not just the one on which we have chosen to base our constructions. Observe also that $\mu_\xi(\vp_F) = \pm1$ and, indeed, takes only the value $+1$ when $d$ is odd. 
\endremark 
\subhead 
9.7 
\endsubhead 
We prove the Linear Independence Lemma of 9.5. To summarize the key hypotheses, the representation $\rho$ of $G_K$ is irreducible, cuspidal and totally ramified, the endo-class $\vt(\rho)\in \scr E(K)$ is $\vT_K$, and $E/K$ is a tame parameter field for $\vT_K$. 
\proclaim{Lemma 1} 
The endo-classes $\vt(\rho^\gamma)$, $\gamma\in \vD\backslash\vG$, are distinct. 
\endproclaim 
\demo{Proof} 
Recall that $K\cap E_0 = K_0/F$ is the maximal unramified sub-extension of $E_0/F$. We think of $\vD\backslash\vG$ as $\Gal {K_0}F$. We set $\vT_0 = \frak i_{K/K_0}(\vT_K)$, so that $\vT_0$ is a $K_0/F$-lift of $\vT$. The set of $K_0/F$-lifts of $\vT$ is therefore the orbit $\{\vT_0^\gamma:\gamma\in \vD\backslash\vG\}$, the conjugates $\vT_0^\gamma$, $\gamma\in \vD\backslash\vG$, being distinct. A tame parameter field for $\vT_0$ is provided by $E_0/K_0$. The field extensions $K/K_0$, $E_0/K_0$ are linearly disjoint, whence $\vT_0^\gamma$ has a unique $K/K_0$-lift, namely $\vt(\rho^\gamma)$. \qed 
\enddemo 
We temporarily set $s=[P{:}K] = [P_0{:}K_0]$ and $t = [K_0{:}F] = |\vD\backslash\vG|$. We fix a prime element $\vp$ of $P_0$ and let $\frak S_a = \vp^aJ^1_K\cap \LR{(G_K)}$. We recall that $\rho = \cind_{\bk J_K}^{G_K}\,\kappa_K$. In particular, the central character of $\rho$ is trivial on $\bs\mu_K$ and a pre-determined prime element $\vp_F$ of $F$. Consequently: 
\proclaim{Lemma 2} 
Let $g\in \LR{(G_K)}$ and suppose $\roman{tr}\,\rho(g)\neq 0$. The element $g$ then has a $G_K$-conjugate lying in $\mu\frak S_a$, where $a = \ups_K(\det g)$ and $\mu$ is some element of $\bs\mu_P = \bs\mu_K$. 
\endproclaim 
Consider the set of $st$ functions $\roman{tr}\,\chi\rho^\gamma$ on $\LR{(G_K)}$, where $\gamma$ ranges over $\vD\backslash\vG$ and $\chi$ over $X_0(K)_s$. This set is linearly independent. Let $\scr D$ denote the space it spans. 
\par 
Each of the representations $\chi\rho^\gamma$ has central character trivial on $\ags{\bs\mu_K,\vp_F}$. Consequently, Lemma 2 implies that their characters form a linearly independent set of functions on the space $\bigcup_{0\le a\le s{-}1} \frak S_a$. Let $\delta_a$ denote the dimension of the space $\{f|_{\frak S_a}: f\in \scr D\}$. For fixed $\gamma$ and $a$, the functions $\roman{tr}\,\chi\rho^\gamma|_{\frak S_a}$, $\chi\in X_0(K)_s$, are constant multiples of each other. We conclude that $\delta_a\le t$, for all $a$. However, 
$$
st = \dim{\scr D} \le \sum_{0\le a\le s{-}1} \delta_a \le st. 
$$
It follows that $\delta_a=t$ for all $a$, whence the functions $\roman{tr}\,\rho^\gamma$ are linearly independent on each set $\frak S_a$. This proves the Linear Independence Lemma. \qed  
\head\Rm 
10. Discrepancy at a prime element 
\endhead 
We continue our investigation of the character $\mu_\xi$ of (9.5.3), attached to a character $\xi\in X_1(E)^{\vD \text{\rm -reg}}$. We take a particular prime element $\vp$ of $E_0$ and show that the value $\mu_\xi(\vp)$ is independent of $\xi$. It will follow that $\mu_\xi$ is independent of $\xi$, as required to finish the proof of 9.2 Theorem. We shall rely on the fact that the proof of that result is already complete when $E_0/F$ is unramified (9.6 Corollary). 
\subhead 
10.1 
\endsubhead  
We have chosen in \S5 a prime element $\vp_F$ of $F$; we take $\vp\in C_{E_0}(\vp_F)$ (notation of 5.2). We set $L = F[\vp]$, and note that $E/L$ is {\it unramified.} We put $\vU = \roman{Aut}(E|F)$ and $\vU_L = \Gal EL \subset \vU$, so that $\vD\subset \vU_L\subset \vU$. Other notation is as in \S9, particularly $d = [K{:}F]$. 
\par 
We return to the automorphic induction equation as in (9.2.4), but incorporating the conclusion of 9.5 Proposition. It reads   
$$
\roman{tr}^\vk\,\pi_\chi(g) = (-1)^{m-1}\eps^0_F(\bs\mu_K)\,\bs\delta_{K/F}(g) \sum_{\gamma\in \vG} \roman{tr}\,\tau_\xi^\gamma(g), \quad  g\in G_K\cap \LR G. 
\tag 10.1.1
$$
(Since we will vary the base field, it is now necessary to specify the field extension to which the transfer factor $\bs\delta$ is attached.) We evaluate each side of (10.1.1) at an element $g = \vp h$, where $h$ satisfies the following conditions: 
\proclaim{(10.1.2) Hypotheses} 
\roster 
\item $h\in \pre p{\bk J_E} = \pre p{\bk J_E}(\vp_F)$; 
\item the group $\pre p{\bk J_E}$ is generated by $h$ and $J^1_E$; 
\item the element $\vp h$ lies in $G_E\cap \LR G$. 
\endroster 
\endproclaim 
We recall that $\pi_\chi$ contains the extended maximal simple type $\vL_\chi = \lambda^\bk J_\chi\otimes\kappa\in \scr T(\theta)$. 
\proclaim{Proposition} 
Let $g = \vp h$, where $h$ satisfies the conditions \rom{(10.1.2).} Let $x\in G$ satisfy $x^{-1}gx\in \bk J$, and suppose that $\roman{tr}\, \vL_\chi (x^{-1}gx) \neq 0$. The coset $x\bk J$ then contains an element $y$ for which $y^{-1}\vp y \in E^\times$. 
\endproclaim 
\demo{Proof} 
Write $g' = x^{-1}gx$, $g'_1 = x^{-1}\vp x$ and $g'_2 = x^{-1}hx$. In the group $G/\ags{\vp_F}$, the element $g'_1$ has finite order relatively prime to $p$, while $g'_2$ is pro-unipotent. Moreover, the elements $g'_i$ commute with each other. It follows that any closed subgroup of $G$ containing $\vp_F$ and $g'$ must also contain both $g'_i$. 
\par 
Since $\vp$ is central in $\bk J/J^1$, we may write $g'_1 = \vp h_1$, where $h_1\in J^0$ has finite order, relatively prime to $p$, in $J^0/J^1$. There exists $\vp'\in P_0^\times\cap \pre p{\bk J}$ such that $g'_2 = \vp' h_2$, where $h_2\in J^0$. This element $\vp'$ is also central in $\bk J/J^1$. Thus $h_2$ is unipotent in $J^0/J^1$, Moreover, the elements $h_i$ commute with each other in $J^0/J^1$. 
\par 
We have $g' \equiv h_1h_2z \pmod{J^1}$, for some $z$ which is central in $\bk J/J^1$. As $\lambda^\bk J_\chi$ is an irreducible representation of $\bk J$, trivial on $J^1$, we deduce that $\roman{tr}\,\lambda^\bk J_\chi(g') \neq 0$ if and only if $\roman{tr}\,\lambda^\bk J_\chi(h_1h_2) \neq 0$. 
\par 
As in the proof of 6.1 Proposition 12 of \cite{13}, the condition $\roman{tr}\,\lambda^\bk J_\chi (h_1h_2)\neq 0$ implies that $h_1$ is conjugate in $J^0/J^1$ to an element of $\bs\mu_K$. We may therefore adjust $x$ within the coset $x\bk J$ in order to assume that $h_1 = \alpha_1j$, where $\alpha_1\in \bs\mu_K$ and $j\in J^1$. The element $g'_1 = \vp h_1$ acts on $J^1$ (by conjugation) as an automorphism of finite order, relatively prime to $p$. Applying the Conjugacy Lemma of 2.6, we may further adjust $x$ by an element of $J^1$ to assume that $j$ commutes with $h_1\vp$. However, $\alpha_1\vp j = g_1'$ has finite $p$-prime order in $G/\ags{\vp_F}$. It follows that $j=1$, whence $g_1'\in E^\times$. \qed 
\enddemo 
The extension $E/L$ is unramified. Consequently, any $F$-embedding $L\to E$ extends to an $F$-automorphism of $E$. The set of $x\in G$ for which $x^{-1}\vp x\in E$ is therefore $G_LN_E$, where $N_E$ denotes the $G$-normalizer of $E^\times$. Thus $N_E/G_E$ is canonically isomorphic to $\vU = \roman{Aut}(E|F)$. 
\par 
We choose a section $\alpha \mapsto g_\alpha$ of the canonical map $N_E/G_E \to \vU$, such that $g_\alpha\in J^0$ when $\alpha\in \vD$. The coset space $G_LN_E\bk J/\bk J$ can then be decomposed as a disjoint union 
$$
G_LN_E\bk J/\bk J = \bigcup_\alpha G_Lg_\alpha\bk J/\bk J = \bigcup_\alpha g_\alpha G_{L^\alpha} \bk J/\bk J = \bigcup_\alpha g_\alpha G_{L^\alpha} /\bk J_{L^\alpha}, 
$$ 
where $\alpha$ ranges over $\vU_L\backslash \vU/\vD$. Here, we view $\vD$ as $\Gal E{E_0}$ and note that $\bk J_{L^\alpha}$ is not necessarily the same as $(\bk J_L)^\alpha$. We accordingly decompose the $\vk$-trace 
$$
\roman{tr}^\vk\,\pi_\chi(g) = \sum_{\alpha\in \vU_L\backslash \vU/\vD} \scr L_\chi(g;\alpha), 
\tag 10.1.3 
$$
where 
$$
\scr L_\chi(g;\alpha) = \sum_{y\in G_L^\alpha/\bk J_{L^\alpha}} \vk(\det y^{-1})\,\roman{tr}\, \vL_\chi (y^{-1}g^\alpha y). 
\tag 10.1.4 
$$ 
\remark{Remark} 
For $\alpha\in \vU$, the field $E_0^\alpha\subset E$ is isomorphic to a tame parameter field for $\theta$. However, the decomposition of (10.1.3) and 2.6 Proposition show that it will not be a tame parameter field unless it equals $E_0$. 
\endremark 
\subhead 
10.2 
\endsubhead 
We return to (10.1.1) and write 
$$
\scr R(g) = (-1)^{m-1}\eps^0_F(\bs\mu_K)\bs\delta_{K/F}(g) \sum_{\gamma\in \vG} \roman{tr}\,\tau_\xi^\gamma(g), \quad  g\in G_K\cap \LR G. 
$$ 
We analyze this expression at the element $g = \vp h$, where $h\in G_E$ satisfies (10.1.2). 
\proclaim{Lemma 1} 
Let $\gamma\in \vG$ and suppose that $\roman{tr}\,\tau^\gamma(\vp h) \neq 0$. The automorphism $\gamma$ of $K$ then extends to an $F$-automorphism of $E$. 
\endproclaim 
\demo{Proof } 
By hypothesis, $\vp h\in G_E\cap \LR G$. The algebra $F[\vp h]$ is therefore a field, of degree $n$ over $F$. It commutes with, and therefore contains, $K$. Condition (10.1.2)(2) implies $F[\vp h]/K$ to be totally ramified, of degree $n/d$. The argument is now identical to that of 8.4 Lemma. \qed 
\enddemo 
Lemma 1 implies 
$$
\scr R(\vp h) = (-1)^{m-1\,} \eps^0_F(\bs\mu_K)\,\bs\delta_{K/F}(\vp h) \sum_{\gamma \in \vU/\vU_K} \roman{tr}\, \tau_\xi^\gamma(\vp h), 
$$ 
where $\vU_K$ is the subgroup $\roman{Aut}(E|K)$ of $\vU$. We take $\kappa\in \scr H(\theta)$, as in 9.2. We form $\kappa_K = \ell_{K/F}(\kappa) \in \scr H(\theta_K)^\vD$, as in 5.6 Proposition, and then $\kappa_E = (\kappa_K)_E\in \scr H(\theta_E) = \scr T(\theta_E)$, as in 5.3 Proposition. We set 
$$ 
\rho = \cind_{\bk J_E}^{G_E}\,\kappa_E. 
\tag 10.2.1 
$$ 
Thus $\tau_\xi \cong \text{\it ind}_{E/K}\,\xi\rho$ and, equivalently, $(\tau_\xi)_E = \xi\rho$. Applying 8.1 Proposition, relative to the base field $K$, we get
$$
\roman{tr}\,\tau_\xi(\vp h) = \eps_K(\vp)\sum_{\beta\in \vU_K} \roman{tr}\,(\xi\rho)^\beta(\vp h), 
$$ 
where $\eps_K(\vp)$ is the symplectic sign $t_{\ags\vp}(J^1_K/H^1_K)$. In all, 
$$
\scr R(\vp h) = (-1)^{m-1\,} \eps^0_F(\bs\mu_K)\,\bs\delta_{K/F}(\vp h) \sum_{\gamma\in \vU} \eps_K(\vp^\gamma)\,\roman{tr}(\xi \rho)^\gamma(\vp h). 
$$ 
The element $\vp^\gamma/\vp$ lies in $\bs\mu_K$, and so acts trivially on $J^1_K/H^1_K$. Therefore $\eps_K(\vp^\gamma) = \eps_K(\vp)$ and so 
$$
\scr R(\vp h) = (-1)^{m-1\,} \eps^0_F(\bs\mu_K)\,\eps_K(\vp)\,\bs\delta_{K/F}(\vp h) \sum_{\gamma\in \vU} \roman{tr}(\xi \rho)^\gamma(\vp h). 
\tag 10.2.2 
$$  
We write the inner sum here in the form $\sum_{\psi\in\vU} \roman{tr}\,\xi\rho(\vp^\psi h^\psi)$, and consider an individual term $\roman{tr}\,\xi\rho(\vp^\psi h^\psi)$. By construction, the central character of $\rho$ is trivial on $C_E(\vp_F)$, so 
$$
\roman{tr}\,\xi\rho(\vp^\psi h^\psi) = \xi({\det}_E(\vp^\psi h^\psi))\,\roman{tr}\,\rho(h^\psi).  
$$
The element $h\in \pre p{\bk J_E}$ satisfies $h^{p^r} \equiv \vp_F^a\pmod{J^1_E}$, for an integer $a$ relatively prime to $p^re$. It follows that $\det_E h\equiv \vp_F^a \pmod{U^1_E}$, whence $\xi({\det}_E\,h^\psi) = \xi(\vp_F)^a$. For the other factor, there is a root of unity $\zeta_\psi\in \bs\mu_E$ such that $\vp^\psi = \zeta_\psi\vp$, giving $\xi({\det}_E\,\vp^\psi) = \xi(\zeta_\psi)^{p^r}\,\xi(\vp)^{p^r}$. We therefore have 
$$
\multline 
\scr R(\vp h) \\ 
= (-1)^{m-1} \eps^0_F(\bs\mu_E)\,\eps_K(\vp)\,\bs\delta_{K/F}(\vp h)\, \xi(\vp_F^a\vp^{p^r}) \\ 
\sum_{\psi\in \vU} \xi(\zeta_\psi)^{p^r}\,\roman{tr}\,\rho(h^\psi) . 
\endmultline  \tag 10.2.3 
$$ 
\indent 
Before passing on, we record some properties of the root of unity $\zeta_\psi$. 
\proclaim{Lemma 2} 
\roster 
\item 
If $\gamma\in \vU_L$, then $\zeta_\gamma = 1$. In particular, $\zeta_\delta = 1$ if $\delta\in \vD$. 
\item 
If $\psi = \gamma\alpha\delta$, $\gamma\in \vU_L$, $\delta \in \vD$, then $\zeta_\psi = \zeta_\alpha^\delta$. 
\endroster 
\endproclaim 
\demo{Proof} 
If $\beta\in \vU_L$, then $\vp^\beta = \vp$ by definition, and this proves (1). The map $\psi\mapsto \zeta_\psi$ satisfies $\zeta_{\psi\phi} = \zeta_\psi^\phi\zeta_\phi$, $\psi,\phi\in \vU$, whence (2) follows. \qed 
\enddemo 
The representation $\rho$ is fixed by $\vD$. Indeed, $\vD$ is the $\vU$-isotropy group of $\rho$. Using Lemma 2, we may now re-write (10.2.3) in the form 
$$ 
\multline 
\scr R(\vp h) \\ 
= (-1)^{m-1} \eps^0_F(\bs\mu_E)\,\eps_K(\vp)\,\bs\delta_{K/F}(\vp h)\, \xi(\vp_F^a\vp^{p^r}) \\ 
\sum_{\psi\in \vU/\vD}\roman{tr}\,\rho(h^\psi) \sum_{\delta\in \vD} \xi(\zeta_\psi^{p^r})^\delta. 
\endmultline 
\tag 10.2.4 
$$ 
\subhead 
10.3 
\endsubhead 
We return to the expression (10.1.4), 
$$
\scr L_\chi(\vp h;\alpha) = \sum_{y\in G_L^\alpha/\bk J_{L^\alpha}} \vk(\det y^{-1})\,\roman{tr}\, \vL_\chi (y^{-1}\vp^\alpha h^\alpha y), 
$$
for $\alpha\in \vU_L\backslash \vU/\vD$. Recalling that $\vL_\chi = \lambda_\chi^{\bk J} \otimes \kappa$, we first have 
$$
\roman{tr}\,\lambda_\chi^\bk J(\vp^\alpha y^{-1}h^\alpha y) = \chi(\vp)^{p^r}\,\roman{tr}\,\lambda_\chi^\bk J(\zeta_\alpha y^{-1}h^\alpha y). 
\tag 10.3.1 
$$ 
To evaluate the expression (10.3.1), we first use the method of \cite{13 (6.10)}. We note that $E_0L^\alpha/L^\alpha$ is a tame parameter field for $\theta_{L^\alpha}$ and $\Gal E{E_0L^\alpha} = \vD\cap \vU_L^\alpha$. We put $m_\alpha = [E{:}E_0L^\alpha]$. Lemma 6.10 of \cite{13} then yields 
$$
\roman{tr}\,\lambda_\chi^\bk J(\zeta_\alpha y^{-1}h^\alpha y) = (-1)^{m-m_\alpha} \sum_{\beta\in \vD\cap \vU_L^\alpha\backslash \vD} \roman{tr}\,\lambda_{\chi^\beta}^{\bk J_{L^\alpha}} (\zeta_\alpha y^{-1}h^\alpha y). 
\tag 10.3.2 
$$ 
Second, we use the Glauberman correspondence relative to the action of the element $\vp^\alpha$ of $C_E(\vp_F)$. We obtain a representation $\kappa_{L^\alpha} \in \scr H(\theta_{L^\alpha})$ and a constant $\eps_F(\vp^\alpha) = t_{\ags{\vp^\alpha}}(J^1/H^1) = \pm1$ such that 
$$
\roman{tr}\,\kappa(\vp^\alpha y^{-1}h^\alpha y) = \eps_F(\vp^\alpha)\, \roman{tr}\,\kappa_{L^\alpha} (y^{-1}h^\alpha y). 
\tag 10.3.3 
$$ 
For $\beta \in \vD$, we set 
$$
\sigma_{\alpha,\beta} = \cind_{\bk J_{L^\alpha}}^{G_{L^\alpha}}\,\lambda^{\bk J_{L^\alpha}}_{\chi^\beta} \otimes \kappa_{L^\alpha} \in \Aa{m_\alpha}{L^\alpha}{\vT_{L^\alpha}} , 
$$
where $\vT_{L^\alpha} = \text{\it cl\/}(\theta_{L^\alpha})$. Setting $\vk_{L^\alpha} = \vk_L^\alpha = \vk \circ \N{L^\alpha}F$, the relations (10.3.2), (10.3.3) yield 
$$
\scr L_\chi(\vp h;\alpha) = (-1)^{m-m_\alpha}\,\eps_F(\vp^\alpha) \sum_{\beta\in \vD\cap \vU_L^\alpha \backslash \vD} \roman{tr}^{\vk_L^\alpha} \sigma_{\alpha,\beta}(\vp^\alpha h^\alpha). 
$$
We may now form the representation $\ell_{E/L^\alpha}(\kappa_{L^\alpha})\in \scr H(\theta_E)$. Transitivity of the Glauberman correspondence (5.1.1) implies that $\ell_{E/L^\alpha}(\kappa_{L^\alpha})$ is equivalent to the representation $\kappa_E$ defined in 10.2. As there, we set $\rho = \cind_{\bk J_E}^{G_E}\,\kappa_E$. It follows that $\sigma_{\alpha,\beta}$ is automorphically induced from a representation $\chi^\beta\eta_{\alpha,\beta}\cdot\rho$ of $G_E$, where $\eta_{\alpha,\beta}\in X_1(E)$ is given by 9.6 Corollary. The explicit formula (9.6.5) shows $\eta_{\alpha,\beta}$ depends only on the subfield $L^\alpha$ of $P$, so we write $\eta_{\alpha,\beta} = \eta_\alpha$. Therefore 
$$ 
\multline 
\roman{tr}^{\vk_L^\alpha} \sigma_{\alpha,\beta}(\vp^\alpha h^\alpha) \\ 
= (-1)^{m_\alpha-1}\,\eps^0_{L^\alpha}(\bs\mu_E)\, \bs\delta_{E/L^\alpha}(\vp^\alpha h^\alpha) \sum_{\delta\in \vU_L^\alpha} \roman{tr}(\chi^\beta \eta_\alpha\cdot\rho)^\delta (\vp^\alpha h^\alpha).  
\endmultline 
$$
The definition gives $\bs\delta_{E/L^\alpha}(\vp^\alpha h^\alpha) = \bs\delta_{E/L}(\vp h)$. The representations $\eta_\alpha$, $\rho$ are fixed by $\vD$ and, in particular, by the element $\beta$. So, adding up, we get 
$$ 
\multline 
\scr L_\chi(\vp h;\alpha) \\ 
= (-1)^{m-1}\, \eps_F(\vp^\alpha)\, \eps^0_{L^\alpha}(\bs\mu_E)\, \bs\delta_{E/L}(\vp h) 
\sum_{\psi\in \vU_L \alpha\vD} \roman{tr}\,(\chi\eta_\alpha\cdot\rho)((\vp h)^\psi). 
\endmultline 
\tag 10.3.4 
$$ 
\subhead 
10.4 
\endsubhead  
In (10.3.4), we consider the term 
$$ 
\roman{tr}\,(\chi\eta_\alpha\cdot\rho)((\vp h)^\psi) = \chi({\det}_E(\vp h)^\psi)\,\eta_\alpha({\det}_E(\vp h)^\psi)\, \roman{tr}\,\rho(h^\psi), 
$$
for $\psi \in \vU_L \alpha\vD$. The element $h$ lies in $\pre p{\bk J}(\vp_F)$, and satisfies $h^{p^r} \equiv \vp_F^a \pmod{J^1}$. This implies $\det_E h\equiv \vp_F^a \pmod{U^1_E}$, and therefore 
$$
\chi({\det}_E\,h^\psi) = \chi(\vp_F^a).  
\tag 10.4.1 
$$
Also, 
$$
\chi({\det}_E\, \vp^\psi) = \chi(\vp^\psi)^{p^r} = \chi(\zeta_\psi^{p^r})\,\chi(\vp)^{p^r},  
\tag 10.4.2  
$$ 
where $\zeta_\psi = \vp^\psi/\vp \in \bs\mu_K$, as in 10.2. 
\par 
We next consider the contributions $\eta_\alpha(\det_E(\vp h)^\psi)$. We write $\psi = \gamma\alpha\delta$, with $\gamma\in \vU_L$ and $\delta\in \vD$. It will be useful to have the notation 
$$ 
d_L = [E{:}L] = [E{:}L^\alpha], 
$$ 
and to recall that $e = e(E|F)$. By (9.6.5), 
$$ 
\eta_\alpha(\vp^\alpha) = \vk_L^\alpha(\vp^\alpha)^{[P{:}L](d_L-1)/2} = \vk(\vp_F)^{n(d_L-1)/2e}. 
$$ 
Therefore  
$$
\eta_\alpha({\det}_E(\vp^\alpha)) = \vk(\vp_F)^{p^rn(d_L-1)/2e}. 
$$ 
Since $(\vp^\alpha)^\delta = \zeta_{\alpha\delta}\zeta_\alpha^{-1}\vp^\alpha = \zeta_\alpha^\delta \zeta_\alpha^{-1}\vp^\alpha$, we get  
$$ 
\align 
\eta_\alpha({\det}_E(\vp^\psi)) &= \eta_\alpha({\det}_E(\vp^{\alpha\delta})) \\ 
&= \eps^1_{L^\alpha}(\bs\mu_E, \zeta_\alpha^\delta\zeta_\alpha^{-1})\,\vk(\vp_F)^{p^rn(d_L-1)/2e}. 
\endalign 
$$ 
All $\eps^1$-characters of $\bs\mu_K$ have order $\le2$, and so are fixed by all automorphisms of $\bs\mu_K$. Therefore  
$$ 
\eta_\alpha({\det}_E(\vp^\psi)) = \vk(\vp_F)^{p^rn(d_L-1)/2e}. 
\tag 10.4.3 
$$ 
The next term to consider is $\eta_\alpha(\det_E h^\psi) = \eta_\alpha(\vp_F^a)$. We have $\vp^e = \nu\vp_F$, for some $\nu\in \bs\mu_L$, whence $(\vp^\alpha)^e = \nu^\alpha\vp_F$ and $\nu^\alpha\in \bs\mu_{L^\alpha}$. In particular, $\eps^1_{L^\alpha}(\bs\mu_E;\nu^\alpha) = 1$. This gives us 
$$
\eta_\alpha(\vp_F) = \eta_\alpha(\vp^\alpha)^e = \kappa(\vp_F)^{n(d_L-1)/2} ,  
$$ 
and so 
$$
\eta_\alpha({\det}_E\,h^\psi) = \eta_\alpha(\vp_F^a) = \vk(\vp_F)^{an(d_L-1)/2}. 
\tag 10.4.4 
$$
Incorporating (10.4.1--4) into (10.3.4), it becomes 
$$ 
\multline 
\scr L_\chi(\vp h;\alpha)  \\ 
= (-1)^{m-1}\, \bs\delta_{E/L}(\vp h)\, \chi(\vp_F^a\vp^{p^r}) \, 
\vk(\vp_F)^{(p^rn/e + an)(d_L-1)/2} \\
\sum_{\psi\in \vU_L \alpha\vD}\eps_F(\vp^\alpha)\, \eps^0_{L^\alpha}(\bs\mu_E)\, \chi(\zeta_\psi^{p^r})\, \roman{tr}\,\rho(h^\psi).
\endmultline 
\tag 10.4.5 
$$ 
The factor $\chi(\zeta_\psi)$ is $\xi(\zeta_\psi) \eps^1_F(\bs\mu_E; \zeta_\psi)$ (9.6 Corollary). If $\psi = \gamma\alpha\delta$, we have $\zeta_\psi = \zeta_\alpha^\delta$, and so 
$$
\chi(\zeta_\psi)^{p^r} = \xi(\zeta_\psi)^{p^r}\,\eps^1_F(\bs\mu_E;\zeta_\alpha). 
$$ 
We introduce a new function 
$$ 
\eurm k(h) = \vk(\vp_F)^{(p^rn/e + an)(d_L-1)/2} = \pm1, 
$$
where $a = a(h)$ is the integer defined by $h^{p^r}\equiv \vp_F^a\pmod{J^1}$. 
\remark{Remark} 
The quantity $\eurm k(h)$ is the product of the constant $\vk(\vp_F)^{p^rn(d_L-1)/2e} = \pm1$ and the character $\pre p{\bk J_E}\to \{\pm1\}$ given by $h\mapsto \vk(\vp_F)^{an(d_L-1)/2}$, where $a = a(h)$ is defined by $h^{p^r} \equiv \det_E h \equiv \vp_F^a \pmod{J^1}$. This character is non-trivial exactly when $d_L$ is even and $ep^r = n/d$ is odd. 
\endremark 
We return to the automorphic induction relation (10.1.1) in the form 
$$
\scr R(g) = \sum_{\alpha\in \vU_L\backslash \vU/\vD} \scr L_\chi(g;\alpha). 
$$ 
Following (10.2.3) and (10.4.5), it becomes 
$$
\multline 
 \eps_F^0(\bs\mu_E)\eps_K(\vp)\, \xi(\vp_F^a\vp^{p^r}) \,\bs\delta_{K/F}(\vp h) 
\sum_{\psi\in \vU} \xi(\zeta_\psi)^{p^r}\,\roman{tr}\,\rho(h^\psi) \\ 
= 
\eurm k(h)\,\chi(\vp_F^a\vp^{p^r})\,\bs\delta_{E/L}(\vp h) \sum_{\alpha\in \vU_L\backslash \vU/\vD} \eps_F(\vp^\alpha)\, \eps^0_{L^\alpha}(\bs\mu_E)\, \eps^1_F(\bs\mu_E;\zeta_\alpha) \\ 
\sum_{\psi\in \vU_L\alpha \vD} \xi(\zeta_\psi)^{p^r}\,\roman{tr}\,\rho(h^\psi). 
\endmultline 
\tag 10.4.6 
$$ 
By (9.6.5), we have 
$$
\chi(\vp_F^a) = \xi(\vp_F^a)\,\vk(\vp_F)^{an(d-1)/2}. 
\tag 10.4.7 
$$ 
To simplify further, we need to control the coefficients in the sum. This is achieved via: 
\proclaim{Symplectic Signs Lemma} 
The function $\vU_L\backslash \vU/\vD \to \{\pm1\}$, given by 
$$
\vU_L\alpha \vD \longmapsto \eps_F(\vp^\alpha)\, \eps^0_{L^\alpha}(\bs\mu_E)\, \eps^1_F(\bs\mu_E;\zeta_\alpha), 
$$
is constant. 
\endproclaim 
We prove this lemma in \S11 below. In particular, it implies 
$$ 
\eps_F(\vp^\alpha)\, \eps^0_{L^\alpha}(\bs\mu_E)\, \eps^1_F(\bs\mu_E;\zeta_\alpha) = \eps_F(\vp)\, \eps^0_L (\bs\mu_E), \quad \alpha\in \vU. 
\tag 10.4.8 
$$ 
We accordingly set 
$$ 
\aligned 
\eurm k'(h) &= \eps_K(\vp)\,\eps_F(\vp)\, \eps^0_L (\bs\mu_E)\,\eps^0_F(\bs\mu_E)\, \vk(\vp_F)^{an(d-1)/2}\, \eurm k(h) \\ 
&= \eps_K(\vp)\,\eps_F(\vp)\, \eps^0_L (\bs\mu_E)\,\eps^0_F(\bs\mu_E)\, \vk(\vp_F)^{p^rn(d_L-1)/2e + an(d-d_L)/2} \\ &= \pm1. 
\endaligned 
\tag 10.4.9 
$$
\remark{Remark} 
As for the function $\eurm k$, the quantity $\eurm k'(h)$ is the product of a constant, with value $\pm1$, and the character $\pre p{\bk J_E}\to \{\pm1\}$ given by $h\mapsto \vk(\vp_F)^{an(d-d_L)/2}$. This character is non-trivial exactly when $d\not\equiv d_L\equiv ep^r \equiv 1\pmod2$. 
\endremark 
The fundamental relation (10.4.6) so becomes 
$$
\multline 
 \xi(\vp^{p^r}) \,\bs\delta_{K/F}(\vp h) \sum_{\psi\in \vU/\vD} \roman{tr}\,\rho(h^\psi) \sum_{\delta\in \vD} \xi(\zeta^\delta_\psi)^{p^r} \\ 
= \eurm k'(h)\,\chi(\vp^{p^r})\,\bs\delta_{E/L}(\vp h) \sum_{\psi\in \vU/\vD} \roman{tr}\,\rho(h^\psi) \sum_{\delta\in \vD} \xi(\zeta^\delta_\psi)^{p^r}  
\endmultline 
\tag 10.4.10 
$$ 
({\it cf\.} (10.2.4)). 
\subhead 
10.5 
\endsubhead 
To get any further with (10.4.10), we have to cancel the double sum occurring in each side. 
\proclaim{Linear Independence Lemma} 
Let $h$ satisfy \rom{(10.1.2),} and let $S_h$ be the set of $u\in J^1_E$ for which $hu$ satisfes the same conditions. The set of functions $u\mapsto \roman{tr}\,\rho^\psi(hu)$, $\psi\in \vD\backslash\vU$, is then linearly independent on $S_h$. 
\endproclaim 
\demo{Proof} 
Let $S'_h$ be the set of $u\in J^1_E$ for which $\vp hu\in \LR{(G_E)}$. The set of functions $u\mapsto \roman{tr}\,\rho^\psi(hu)$, $\psi\in \vD\backslash \vU$, is then linearly independent on $S'_u$: the argument is the same as the proof of the Linear Independence Lemma of 9.5. However, these functions are all locally constant on $S'_h$, and $S_h$ is dense in $S'_h$. \qed 
\enddemo 
\proclaim{Non-vanishing Lemma} 
Let $h$ satisfy \rom{(10.1.2),} and define $S_h$ as in the Linear Independence Lemma. Let $\xi\in X_1(E)^{\vD \text{\rm -reg}}$. There exists $u\in S_h$ such that 
$$
\sum_{\psi\in \vU/\vD} \roman{tr}\,\rho((hu)^\psi) \sum_{\delta\in \vD} \xi(\zeta^\delta_\psi)^{p^r} \neq 0. 
$$
\endproclaim 
\demo{Proof} 
In light of the Linear Independence Lemma, it is enough to show that one of the coefficients $\sum_\delta \xi(\zeta^\delta_\psi)^{p^r}$, $\psi\in \vU/\vD$, is non-zero. However, if we take $\psi = 1$, this coefficient reduces to $|\vD| = m$. \qed 
\enddemo 
The character $\xi^{-1}\chi$ has finite order (9.5 Corollary), so we conclude: 
\proclaim{Proposition} 
Let $\xi\in X_1(E)^{\vD \text{\rm -reg}}$, Let $\chi = \chi_\xi$ and let $h$ satisfy \rom{(10.1.2).} There exists $u\in S_h$ such that 
$$
\xi^{-1}\chi(\vp^{p^r}) = \eurm{k'}(hu)\,\bs\delta_{K/F}(\vp hu)/\bs\delta_{E/L}(\vp hu). 
\tag 10.5.1
$$
In particular, $\eurm k'(hu) = \eurm k'(h) = \pm1$, and $\bs\delta_{K/F}(\vp hu)/\bs\delta_{E/L}(\vp hu)$ is a root of unity. 
\endproclaim 
\subhead 
10.6 
\endsubhead 
To reduce the notational load, we assume $h$ was chosen so that 10.5 Proposition holds with $u=1$: this is certainly achievable but the choice of $h$ may still depend on $\xi$. 
\proclaim{Transfer Lemma} 
If the quantity $\bs\delta_{K/F}(\vp h)/\bs\delta_{E/L}(\vp h)$ is a root of unity, then 
$$
\bs\delta_{K/F}(\vp h)/\bs\delta_{E/L}(\vp h) = \eurm d\,\vk(\vp_F)^{an(d-d_L)/2},  
$$ 
where $a$ is the integer defined by $h^{p^r}\equiv \vp_F^a\pmod{J^1}$, and $\eurm d$ is a sign depending only on the ramification indices and inertial degrees in the tower of fields $F\subset L\subset P$. 
\endproclaim 
\demo{Proof} 
We set $t= \vp h$ and return to the defining formulas 
$$
\gathered
\aligned 
\Delta^1_{K/F}(t) &= \|\wt\Delta_{K/F}(t)^2\|_F^{1/2}\,\|\det t\|_F^{ep^r(1-d)/2}, \\ \Delta^1_{E/L}(t) &= \|\wt \Delta_{E/L}(t)^2\|_L^{1/2}\,\|{\det}_L\,t\|_L^{p^r(1-d_L)/2}, \endaligned \\ 
\alignedat3 
\Delta^2_{K/F}(t) &= \vk(e_F\wt\Delta_{K/F}(t)), &\quad \Delta^2_{E/L}(t) &=\vk_L(e_L\wt\Delta_{E/L}(t)), \\ 
\bs\delta_{K/F}(t) &= \Delta^2_{K/F}(t)/\Delta^1_{K/F}(t), &\quad \bs\delta_{E/L}(t) &= \Delta^2_{E/L}(t)/ \Delta^1_{E/L}(t). 
\endalignedat 
\endgathered 
$$ 
The elements $e_F$, $e_L$ are units of $K$, $E$ respectively, such that $e_F\wt\Delta_{K/F}(t)\in F^\times$ and $\wt\Delta_{E/L}(t)\in L^\times$. Both are independent of $t$. 
\par 
As $\bs\delta_{K/F}(t)/ \bs\delta_{E/L}(t)$ is a root of unity, it follows that $\Delta^1_{K/F}(t) = \Delta^1_{E/L}(t)$, or that 
$$
\frac{\|\wt\Delta_{K/F}(t)^2\|_F}{\|\det t\|_F^{ep^r(d-1)}} = \frac{\|\wt \Delta_{E/L}(t)^2\|_L}{\|{\det}_L\,t\|_L^{p^r(d_L-1)}} = \frac{\|\N LF\, \wt\Delta_{E/L}(t)^2\|_F}{\|\det t\|_F^{p^r(d_L-1)}}. 
\tag 10.6.1 
$$ 
Let $\roman v_{K/F}$ denote $\ups_K(\wt\Delta_{K/F}(t))$, so that $\Delta^2_{K/F}(t) = \vk(\vp_F)^ {\roman v_{K/F}(t)}$. Let $\roman v_{E/L} = \ups_E(\wt\Delta_{E/L}(t))$. Thus, if $\vp_L$ is a prime element of $L$, then  $$
\Delta^2_{E/L}(t) = \vk_L(\vp_L)^{\roman v_{E/L}} = \vk(\vp_F)^{d\roman v_{E/L}/d_L}. 
$$
Moreover, 
$$
\bs\delta_{K/F}(t)/\bs\delta_{E/L}(t) =\Delta^2_{K/F}(t)/\Delta^2_{E/L}(t) =  \vk(\vp_F)^ {\roman v_{K/F}-d \roman v_{E/L}/d_L}. 
\tag 10.6.2 
$$ 
If $\roman v_0 = \ups_F(\det t)$, the relation (10.6.1) implies 
$$ 
\roman v_{K/F}-d\roman v_{E/L}/d_L = \roman v_0\big(ep^r(d{-}1)/2 - p^r(d_L{-}1)/2\big). 
$$ 
However, $\roman v_0 = ade+p^rd$, giving 
$$
\multline 
\aligned 
\roman v_{K/F}-d\roman v_{E/L}/d_L &= \big(ade{+}p^rd/d_L\big)\big(ep^r(d{-}1)/2 - p^r(d_L{-}1)/2\big) \\ 
&= p^r\bigg(ade^2(d{-}1)/2 - ade(d_L{-}1)/2 
\endaligned \\ 
+ ep^rd(d{-}1)/2 - p^rd(d_L{-}1)/2\bigg). 
\endmultline  
$$ 
The third and fourth terms contribute a constant integer 
$$ 
w = p^{2r}d\big(e(d{-}1)-(d_L{-}1)\big)/2, 
$$ 
depending only on the fields $F\subset L\subset P$. Since $\vk(\vp_F)$ is a primitive $d$-th root of unity and $ep^rd = n$, the relation (10.6.2) reduces to 
$$ 
\bs\delta_{K/F}(t)/\bs\delta_{E/L}(t) = \eurm d\,\vk(\vp_F)^{an(d-d_L)/2}, 
$$ 
where $\eurm d = \vk(\vp_F)^w = \pm1$. \qed 
\enddemo 
We combine the conclusion of the Transfer Lemma with (10.4.9) and (10.5.1) to get 
$$
\xi^{-1}\chi(\vp)^{p^r} = \eurm d'\, \eps_K(\vp)\,\eps_F(\vp)\, \eps^0_L (\bs\mu_E)\,\eps^0_F(\bs\mu_E), 
\tag 10.6.3 
$$ 
where $\eurm d'$ is a sign depending only on the ramification indices and inertial degrees in the tower of fields $F\subset L\subset P$. By 9.6 Corollary, $\xi^{-1}\chi(\vp)$ is an $e$-th root of $\pm1$. Since $e$ is relatively prime to $p$, there is a unique character $\mu_\xi$ of $E^\times$ satisfying (9.6.5) and 
$$
\mu_\xi(\vp)^{p^r} = \eurm d'\, \eps_K(\vp)\,\eps_F(\vp)\, \eps^0_L (\bs\mu_E)\,\eps^0_F(\bs\mu_E) = \pm1. 
\tag 10.6.4
$$ 
This is visibly independent of $\xi\in X_1(E)^{\vD \roman{-reg}}$ and is fixed by $\vD$. 
\par
We have completed the proof of 9.2 Theorem and, with it, the proof of the Unramified Induction Theorem of 9.1. \qed 
\subhead 
10.7 
\endsubhead 
For convenience, we summarize our conclusions concerning the character $\mu = \mu^{K/F}_\alpha$ of 9.2 Theorem. 
\proclaim{Corollary} 
\roster 
\item 
The character $\mu = \mu^{K/F}_\alpha \in X_1(E)^\vD$ satisfies: 
$$ 
\aligned 
\mu|_{U^1_E} &= 1, \\ 
\mu|_{\bs\mu_E} &= \eps^1_{\bs\mu_E} = t^1_{\bs\mu_E}(J^1_\theta/H^1_\theta), \\ 
\mu(\vp_F) &= \vk(\vp_F)^{n(d-1)/2} = \pm1, 
\endaligned 
$$ 
where $\vp_F$ is a prime element of $F$. 
\item 
If $\vp$ is a prime element of $E_0$ lying in $C_{E_0}(\vp_F)$, then 
$$
\mu(\vp)^{p^r} = \eurm d'\,\eps_K(\vp)\,\eps_F(\vp)\,\eps_L^0(\bs\mu_E)\,\eps^0_F(\bs\mu_E),
$$ 
where $\eurm d'$ is a sign depending only on the ramification indices and inertial degrees in the tower of fields $F\subset L = F[\vp]\subset P$. 
\item These conditions determine $\mu$ uniquely, and $\mu$ takes its values in $\{\pm1\}$. 
\endroster 
\endproclaim 
\demo{Proof} 
Parts (1) and (2) have already been done. Setting $e = e(E|F)$, we have $\vp^e = \zeta\vp_F$, for some $\zeta\in \bs\mu_{E_0}$. We choose integers $a$ and $b$ so that $ae+bp^r=1$. We get 
$$
\mu(\vp) = \big(\eps^1_{\bs\mu_E}(\zeta)\mu(\vp_F)\big)^a \cdot \mu(\vp)^{bp^r} = \pm1, 
$$
whence all outstanding assertions follow. \qed 
\enddemo 
\head\Rm 
11. Symplectic signs 
\endhead 
We prove the Symplectic Signs Lemma of 10.4. This concerns the sign invariants attached to symplectic representations of various {\it cyclic\/} groups over the field $\Bbbk = \Bbb F_p$ of $p$ elements. The theory of such representations is developed in \cite{2} and summarized, with some additions and modifications, in \S3 of \cite{13}. To prove the desired result, it is necessary to expand that framework a little, to cover symplectic $\Bbbk$-representations of finite {\it abelian\/} groups of order prime to $p$.  
\subhead 
11.1 
\endsubhead 
 $C$ be a finite abelian group of order not divisible by $p$. A {\it symplectic $\Bbbk C$-module\/} is a pair $(M,h)$, consisting of a finite $\Bbbk C$-module $M$ and a nondegenerate, alternating, bilinear form $h:M\times M\to \Bbbk$ which is $C$-invariant, 
$$
h(cm_1,cm_2) = h(m_1,m_2), \quad m_i\in M,\ c\in C. 
$$ 
\remark{\bf Example} 
Let $V$ be a finite $\Bbbk C$-module. Let $V^* = \Hom{}V\Bbbk$ and let $\ags{\,,\,}$ be the canonical pairing $V\times V^*\to \Bbbk$. We impose on $V^*$ the $C$-action which makes this pairing $C$-invariant. The module $V\oplus V^*$ carries the bilinear pairing 
$$
h((v_1,v_1^*),(v_2,v_2^*)) = \ags{v_1,v_2^*} - \ags{v_2,v_1^*}, 
$$
which is nondegenerate, alternating and $C$-invariant. We refer to the pair $H(V) = (V\oplus V^*,h)$ as the {\it $C$-hyperbolic space\/} on $V$. 
\endremark 
Returning to the general situation, the obvious notions of {\it orthogonal sum\/} and {\it $C$-isometry\/} apply. A symplectic $\Bbbk C$-module $(M,h)$ is said to be {\it irreducible\/} if it cannot be written as an orthogonal sum $(M_1,h_1)\perp (M_2,h_2)$, for non-zero symplectic modules $(M_i,h_i)$. As an instance of the general theory of Hermitian forms over semisimple rings with involution \cite{19}, \cite{20}, we have the following standard properties. 
\proclaim{Lemma} 
\roster 
\item 
If $(M,h)$ is a symplectic $\Bbbk C$-module, then 
$$
(M,h) \perp (M,-h) \cong H(M). 
$$ 
\item 
Let $(M,h)$ be a symplectic $\Bbbk C$-module. There then exist irreducible symplectic $\Bbbk C$-modules $(M_i,h_i)$, $1\le i\le r$, such that 
$$
(M,h) \cong (M_1,h_1) \perp (M_2,h_2) \perp \dots \perp (M_r,h_r). 
$$
The factors $(M_i,h_i)$ are uniquely determined, up to $C$-isometry and permutation. 
\item 
Let $(M,h)$ be an irreducible symplectic $\Bbbk C$-module. If $M$ is not irreducible as $\Bbbk C$-module, there exists an irreducible $\Bbbk C$-module $V$ such that $(M,h) \cong H(V)$. In this case, $(M,h) \cong H(U)$ if and only if either $U\cong V$ or $U\cong V^*$. 
\endroster 
\endproclaim 
\subhead 
11.2 
\endsubhead 
We consider briefly the structure of irreducible $\Bbbk C$-modules. 
\par 
Let $\bar\Bbbk/\Bbbk$ be an algebraic closure, and set $\vO = \Gal{\bar\Bbbk}\Bbbk$. Let $\chi:C\to \bar\Bbbk^\times$ be a character. In particular, the finite subgroup $\chi(C)$ of $\bar\Bbbk^\times$ is {\it cyclic.} Let $\Bbbk[\chi]$ denote the field extension of $\Bbbk$ generated by the values $\chi(c)$, $c\in C$. The group $C$ acts on $\Bbbk[\chi]$ by $c:x\mapsto \chi(c)x$, and so $\Bbbk[\chi]$ provides a finite $\Bbbk C$-module which we denote $V_\chi$. 
\proclaim{Lemma} 
The $\Bbbk C$-module $V_\chi$ is irreducible. The map $\chi\mapsto V_\chi$ induces a bijection between the set of $\vO$-orbits in $\roman{Hom}(C,\bar\Bbbk^\times)$ and the set of isomorphism classes of irreducible $\Bbbk C$-modules. 
\endproclaim 
\demo{Proof} 
The group algebra $\Bbbk C$ is commutative and semisimple. It is therefore isomorphic to a direct product of finite field extensions of $\Bbbk$. Indeed, $\Bbbk C \cong \prod_\chi \Bbbk[\chi]$, where $\chi$ ranges over $\vO \backslash \roman{Hom}(C,\bar\Bbbk^\times)$. \qed 
\enddemo 
Observe that, in this formulation, $(V_\chi)^* \cong V_{\chi^{-1}}$. We may now list the irreducible symplectic $\Bbbk C$-modules. 
\proclaim{Proposition} 
\roster 
\item 
Let $\chi\in \Hom{}C{\bar\Bbbk^\times}$. The symplectic module $H(V_\chi)$ is irreducible except when $\chi\neq \chi^{-1}$ but $\chi^{-1} = \chi^\omega$, for some $\omega \in \vO$. 
\item 
Let $\chi,\psi\in \Hom{}C{\bar\Bbbk^\times}$. The spaces $H(V_\chi)$, $H(V_\psi)$ are $C$-isometric if and only if $\chi$ is $\vO$-conjugate to $\psi$ or $\psi^{-1}$. 
\item 
Suppose $\chi\in \Hom{}C{\bar\Bbbk^\times}$ satisfies $\chi^{-1} = \chi^\omega \neq \chi$, for some $\omega\in \vO$. The module $V_\chi$ then admits a $C$-invariant, nondegenerate alternating form. This form is uniquely determined up to $C$-isometry. 
\endroster 
\endproclaim 
\demo{Proof} 
This is effectively the same as the case of $C$ cyclic, treated in \cite{2 (8.2.2), (8.2.3)}. \qed 
\enddemo 
\proclaim{Corollary} 
Let $(M,h)$ be an irreducible symplectic $\Bbbk C$-module. 
\roster 
\item 
If $h':M\times M\to \Bbbk$ is a $C$-invariant, nondegenerate, alternating form, then $h'$ is $C$-isometric to $h$. 
\item 
The image of $C$ in $\roman{Aut}\,M$ is cyclic. 
\endroster 
\endproclaim 
\demo{Proof} 
If $(M,h)$ is hyperbolic, assertion (1) is straightforward. Otherwise, $M\cong V_\chi\cong \Bbbk[\chi]$, for some $\chi$. In particular, $C$ acts on $V_\chi$ via the cyclic quotient $\chi(C)$ and we are reduced to the case of the proposition. 
\par
Assertion (2) is therefore clear in the case where $(M,h)$ is not hyperbolic. If it is hyperbolic, then $M \cong V_\chi \oplus (V_\chi)^*$, for some $\chi$. The group $C$ acts on both factors via the same cyclic quotient $\chi(C) = \chi^{-1}(C)$. \qed 
\enddemo 
So, if $(M,h)$ is any symplectic $\Bbbk C$-module, the $C$-isometry class of $h$ is determined by the $\Bbbk C$-isomorphism class of $M$. We therefore now tend to suppress explicit mention of the form $h$. 
\subhead 
11.3 
\endsubhead 
Let $M$ be an irreducible symplectic $\Bbbk C$-module. Suppose, for the moment, that $C$ is {\it cyclic.} Following \cite{13 (3.4 Definition 3)}, which we now recall, we define a constant $t_C^0(M) = \pm1$ and a character $t^1_C(M):C\to \{\pm1\}$. If $M$ is hyperbolic, say $M\cong H(V)$, then $t^0_C(M) = +1$ and $t^1_C(M;c)$, $c\in C$, is the signature of the permutation of $V$ induced by the action of $c$. Otherwise, $t^0_C(M) = -1$, and we may think of $M$ as the field $\Bbb F_{p^{2r}}$ of $p^{2r}$ elements, for an integer $r\ge1$. The group $C$ acts via a homomorphism $\chi_M:C \to \Gamma$, where $\Gamma$ is the subgroup of $\Bbb F_{p^{2r}}^\times$ of order $1{+}p^r$. If $p=2$, the character $t^1_C(M)$ is trivial, otherwise it is the composite of $\chi_M$ with the unique character of $\Gamma$ of order $2$. 
\par 
We further set $t_C(M) = t^0_C(M)\,t^1_C(M;c_0)$, for any generator $c_0$ of the cyclic group $C$. It follows from the definitions that these quantities depend only on the image of $C$ in $\roman{Aut}\,M$. 
\par
We return to the general case of a finite {\it abelian\/} group $C$, of order prime to $p$. Let $M$ be an irreducible symplectic $\Bbbk C$-module. The image $\overline C$ of $C$ in $\roman{Aut}\,M$ is cyclic, by 11.2 Corollary (2). We accordingly define 
$$
\align 
t^0_C(M) &= t^0_{\overline C}(M), \\ 
t_C(M) &= t_{\overline C}(M), \\ 
t^1_C(M) &= \text{the inflation to $C$ of $t^1_{\overline C}(M)$.} 
\intertext{The definitions then yield} 
t_C(M) &= t^0_C(M)\,t^1_C(M;c), 
\endalign 
$$
for any $c\in C$ whose image in $\roman{Aut}\,M$ generates $\overline C$. 
\par
We extend the definitions in the obvious way: if $M$ is a symplectic $\Bbbk C$-module, we write it as an orthogonal sum $M = M_1\perp M_2\perp \dots\perp M_r$, where the $M_i$ are irreducible symplectic $\Bbbk C$-modules, and set 
$$
\aligned 
t^k_C(M) &=\prod_{i=1}^r t^k_C(M_i), \quad k=0,1, \\
t_C(M) &=\prod_{i=1}^r t_C(M_i). 
\endaligned 
\tag 11.3.1 
$$ 
\proclaim{Proposition} 
Let $M$ be an irreducible symplectic $\Bbbk C$-module. If $D$ is a subgroup of $C$, the space $M^D$ of $D$-fixed points in $M$ is either $\{0\}$ or $M$. If $M^D = M$, then 
$$
t_D(M) = t_C^0(M)\,t^1_C(M;d), 
$$
for any $d\in D$ generating the image of $D$ in $\roman{Aut}\,M$. 
\endproclaim 
\demo{Proof} 
Since $M^D$ is a $C$-subspace of $M$, the first assertion is clear if $M$ is irreducible as linear $\Bbbk C$-module. Suppose, therefore, that $M = H(V)$, for some irreducible $\Bbbk C$-module $V$. Thus $M = V\oplus V^*$. Again, $V^D$ is either zero or $V$. In the first case, $(V^*)^D$ is zero, in the second it is $V^*$. The first assertion thus follows, and the second is now an instance of 3.5 Proposition 5 of \cite{13}. \qed 
\enddemo 
\subhead 
11.4 
\endsubhead 
We return to the situation and notation of 10.4, in order to prove the Symplectic Signs Lemma. In particular, we have chosen a prime element $\vp_F$ of $F$, and $\vp$ is a prime element of $E_0$ lying in $C_E(\vp_F)$. We are given an automorphism $\alpha\in \vU = \roman{Aut}(E|F)$. The group $C = C_E(\vp_F)/C_F(\vp_F)$ is finite abelian, of order relatively prime to $p$. The space $J_\theta^1/H_\theta^1$ provides a symplectic $\Bbbk C$-module. 
\par 
If $M$ is a symplectic $\Bbbk C$-module, we denote by $M_\alpha$ the space of $\vp^\alpha$-fixed points in $M$. In particular, $M_1$ is the space of $\vp$-fixed points. We have to prove: 
\proclaim{Proposition} 
If $M$ is a symplectic $\Bbbk C$-module, then 
$$
t_{\ags{\vp^\alpha}}(M)\,t^0_{\bs\mu_E}(M_\alpha)\,t^1_{\bs\mu_E}(M;\vp^\alpha/\vp) = t_{\ags\vp}(M)\,t^0_{\bs\mu_E}(M_1). 
\tag 11.4.1 
$$
\endproclaim 
\demo{Proof} 
The $t$-invariants are multiplicative relative to orthogonal sums (11.3.1). We henceforward assume, therefore, that the given symplectic $\Bbbk C$-module $M$ is {\it irreducible.} 
\par 
Consider first the case in which both $\vp$ and $\vp^\alpha$ act trivially on $M$. Each side of the relation (11.4.1) then reduces to $t^0_{\bs\mu_E}(M)$, and there is nothing to prove.  
\par
Suppose next that $\vp$ acts trivially, while $\vp^\alpha$ does not. The desired relation is then 
$$
t_{\ags{\vp^\alpha}}(M)\,t^1_{\bs\mu_E}(M;\vp^\alpha) = t^0_{\bs\mu_E}(M). 
\tag 11.4.2 
$$ 
The proposition of 11.3 shows that $\vp^\alpha$ acts on $M$ with only the trivial fixed point and, further, $t_{\ags{\vp^\alpha}}(M) = t^0_C(M)t^1_C(M;\vp^\alpha)$. However, since $C$ is generated by $\bs\mu_E$ and $\vp$, our hypothesis on $\vp$ implies $t^j_C(M) = t^j_{\bs\mu_E}(M)$, $j=0,1$. That is, $t_{\ags{\vp^\alpha}}(M) = t^0_{\bs\mu_E}(M)\,t^1_{\bs\mu_E}(M;\vp^\alpha)$, whence (11.4.2) follows. 
\par
The next case, where $\vp^\alpha$ acts trivially while $\vp$ does not, is exactly the same. 
\par
We are reduced to the case where neither $\vp$ nor $\vp^\alpha$ acts trivially. If $\bs\mu_E$ acts trivially, the result is immediate. We therefore assume it does not. The desired relation is 
$$
t_{\ags{\vp^\alpha}}(M)\,t^1_{\bs\mu_E}(M;\vp^\alpha/\vp) = t_{\ags\vp}(M). 
\tag 11.4.3 
$$ 
If the root of unity $\vp^\alpha/\vp$ acts trivially on $M$, this relation is immediate. We therefore assume it does not. Using 11.3 Proposition as before, we have 
$$ 
\align 
t_{\ags{\vp^\alpha}}(M) &= t^0_C(M)\,t^1_C(M;\vp^\alpha), \\ 
t_{\ags\vp}(M) &= t^0_C(M)\,t^1_C(M;\vp). 
\intertext{We therefore need to show} 
t^1_{\bs\mu_E}(M;\vp^\alpha/\vp) &= t^1_C(M;\vp^\alpha/\vp). 
\tag 11.4.4 
\endalign 
$$ 
If the symplectic $\Bbb F_pC$-module $M$ is hyperbolic, (11.4.4) follows directly from the definition of the character $t^1_C(M)$. We therefore assume $M$ is not hyperbolic over $C$ (one says that $M$ is {\it anisotropic\/}). Let $\zeta$ be some generator of $\bs\mu_E$. By definition, $t^0_C(M) = -1$, so  
$$ 
t_{\bs\mu_E}(M) = t^0_C(M)\,t^1_C(M;\zeta) = -t^1_C(M;\zeta).  
$$
On the other hand, 
$$
t_{\bs\mu_E}(M) = t^0_{\bs\mu_E}(M)\,t^1_{\bs\mu_E}(M;\zeta). 
$$ 
As in 11.2, $M \cong V_\chi = \Bbbk[\chi]$, where $\chi\in \Hom{}C{\bar\Bbbk^\times}$ satisfies $\chi^{-1} = \chi^\omega \neq \chi$, for some $\omega\in \vO$. In the present case, the $\Bbbk\bs\mu_E$-module $M$ has length 
$$ 
d = \big[\Bbbk[\chi]: \Bbbk[\chi|_{\bs\mu_E}]\big], 
$$ 
and is symplectic. Each irreducible linear component is equivalent to $M_0 = \Bbbk[\chi|_{\bs\mu_E}]$. If $d$ is odd, the space $M_0$ is therefore anisotropic over $\bs\mu_E$. We get $t^0_{\bs\mu_E}(M) = t^0_{\bs\mu_E}(M_0)^d = -1$ and $t^1_{\bs\mu_E}(M) = t^1_C(M)|_{\bs\mu_E}$. The desired relation (11.4.4) follows in this case. 
\par
In the final case, where $d$ is even, the symplectic $\Bbbk\bs\mu_E$-module $M$ is the direct sum of $d/2$ copies of the hyperbolic space $H(M_0)$. We assert that $\bs\mu_E$ acts on $M_0$ via a character of order $2$. Let $\Bbbk_0$ denote the subfield of $\Bbbk[\chi]$ for which $[\Bbbk[\chi]:\Bbbk_0] =2$. The group $\chi(C)$ is then contained in the kernel of the field norm $\N{\Bbbk[\chi]}{\Bbbk_0}$, and this kernel intersects $\Bbbk_0^\times$ in $\{\pm1\}$. As $\Bbbk_0$ contains $\Bbbk[\chi|_{\bs\mu_E}]$, the restriction of $\chi$ to $\bs\mu_E$ has order $\le2$. Indeed, it has order exactly $2$, since we have eliminated the case where $\bs\mu_E$ acts trivially. 
\par 
In particular, $M_0 \cong \Bbbk$. It follows that $t^1_{\bs\mu_E}(M)$ is non-trivial if and only if $p\equiv 3\pmod4$ and $d/2$ is odd (as follows from the definition of $t^1$ in \cite{13}). On the other hand, this same condition is equivalent to the character $t^1_C(M)$ being non-trivial on $\chi^{-1}(\pm1)$. 
\par 
This completes the proof of the proposition. \qed 
\enddemo 
We have also completed the proof of the Symplectic Signs Lemma. 
\head\Rm 
12. Main Theorem and examples  
\endhead 
We complete the proof of the Comparison Theorem of 7.3, and prove the Types Theorem of 7.6. 
\subhead 
12.1 
\endsubhead 
Let $\alpha\in \wP F$ and let $m\ge1$ be an integer. Set $E = Z_F(\alpha)$. Let $E_m/E$ be unramified of degree $m$ and set $\vD = \Gal {E_m}E$. 
\par 
Let $\vS\in \Gg mF{\scr O_F(\alpha)}$. By 1.3 Proposition, there exists $\rho\in \wW E$ such that $\rho|_{\scr P_F} \cong\alpha$. Let $\rho_m$ be the restriction of $\rho$ to $\scr W_{E_m}$. As in 1.6, there exists a $\vD$-regular character $\xi\in X_1(E_m)$ such that $\vS\cong \Ind_{E_m/F}\,\xi \otimes \rho_m$. 
\par
Let $K_m/F$ be the maximal unramified sub-extension of $E_m/F$ and put $\sigma = \Ind_{E_m/K_m}\, \xi\otimes \rho_m$. Thus $\sigma \in \Gg1{K_m}{\scr O_{K_m}(\alpha)}$. Since $E_m/K_m$ is totally ramified, we may apply 8.2 Corollary to get a {\it unique\/} character $\mu = \mu^{K_m}_{1,\alpha} \in X_1(E_m)$ such that 
$$
\upr L\sigma = \mu\odot_{\Phi_{E_m}(\alpha)} \upr N\sigma,  
\tag 12.1.1 
$$ 
for all $\sigma\in \Gg1{K_m}{\scr O_{K_m}(\alpha)}$. The field $K = K_m^\vD/F$ is the maximal unramified sub-extension of $E/F$. The group $\vD$ acts on both $\Gg1{K_m}{\scr O_{K_m}(\alpha)}$ and $\Aa1{K_m}{\Phi_{K_m}(\alpha)}$, via its natural action on $K_m$. Both the Langlands correspondence and the na\"\i ve correspondence preserve these actions (7.1 Proposition, 7.2 Proposition). The uniqueness property of $\mu^{K_m}_{1,\alpha}$ implies it is fixed by $\vD$, that is, 
$$
\mu^{K_m}_{1,\alpha} = \mu^K_{1,\alpha} \circ \N{E_m}E,  
\tag 12.1.2 
$$ 
for a character $\mu^K_{1,\alpha}\in X_1(E)$, uniquely determined modulo $X_0(E)_m$. 
\par 
We remark on one consequence of (12.1.1). Comparing central characters, as in 7.2, we get: 
\proclaim{Proposition} 
The character $\mu^{K_m}_{1,\alpha}$ agrees on $U_{K_m}$ with the discriminant character $d_{E_m/K_m}$ of the extension $E_m/K_m$. The character $\mu^K_{1,\alpha}$ agrees on $U_K$ with $d_{E/K}$. 
\endproclaim 
\subhead 
12.2 
\endsubhead 
For the next step, we have $\vS = \Ind_{K_m/F}\,\sigma$, so 
$$
\upr L\vS = \roman A_{K_m/F}\,\upr L\sigma, \quad \upr N\vS = \text{\it ind}_{K_m/F}\,\upr N\sigma. 
$$ 
In particular, 
$$
\upr L\vS = \mu^K_{1,\alpha}\odot_{\Phi_E(\alpha)} \roman A_{K_m/F}\,\upr N\sigma, 
\tag 12.2.1 
$$
by (9.1.3). The Unramified Induction Theorem of 9.1 (or 9.2 Theorem) gives a unique character $\mu^{K_m/F}_\alpha \in X_1(E_m)^\vD$ such that 
$$
\roman A_{K_m/F}\,\tau = \text{\it ind}_{K_m/F}\,\mu^{K_m/F}_\alpha \odot_{\Phi_{E_m}(\alpha)} \tau, 
$$
for every $\vD$-regular $\tau\in \Aa1{K_m}{\Phi_{K_m}(\alpha)}$. We write 
$$
\mu^{K_m/F}_\alpha = \mu^{K/F}_{m,\alpha}\circ \N{E_m}E, 
$$
for a character $\mu^{K/F}_{m,\alpha} \in X_1(E)$, uniquely determined modulo $X_0(E)_m$. Applying (9.1.3), we find 
$$
\upr L\vS = \mu^{K/F}_{m,\alpha}\,\mu^K_{1,\alpha}\odot_{\Phi_E(\alpha)} \upr N\vS, 
\tag 12.2.2 
$$
and this relation holds for all $\vS\in \Gg mF{\Phi_F(\alpha)}$. The character 
$$
\mu^F_{m,\alpha} = \mu^{K/F}_{m,\alpha}\,\mu^K_{1,\alpha} \in X_1(E) 
\tag 12.2.3 
$$
is uniquely determined modulo $X_0(E)_m$, and satisfies the requirements of the Comparison Theorem. \qed 
\subhead 
12.3 
\endsubhead 
We prove the Types Theorem of 7.6. We use the notation of the statement. In addition, we let $E_m/F$ be the maximal tamely ramified sub-extension of $P_m/F$ and $K = E\cap K_m$
\par 
The representation $\upr N\vS$ contains an extended maximal simple type $\vL_0\in \scr T(\theta)$ of the form $\lambda_\tau\ltimes \nu_0$, where $\nu_0\in \scr H(\theta)$ satisfies $\det\nu_0(\zeta) = 1$, for all $\zeta \in \bs\mu_{K_m}$. By the Comparison Theorem, the representation $\upr L\vS$ contains the representation $\mu^F_{m,\alpha} \odot \vL_0\in \scr T(\theta)$. 
\par 
We evaluate $\mu^F_{m,\alpha}$ on units by first using the decomposition (12.2.3). The factor $\mu^K_{1,\alpha}$ is $d_{E/K}$, which is unramified if $[E{:}K]$ is odd, of order $2$ on $U_K$ otherwise. The other factor of $\mu^F_{m,\alpha}$ is given by (10.7) Corollary, and the result then follows from 5.4 Lemma. \qed 
\baselineskip=12pt \parskip=2pt plus 1pt
 

\Refs 
\ref\no 1 
\by J. Arthur and L. Clozel 
\book Simple algebras, base change, and the advanced theory of the trace formula
\bookinfo Annals of Math. Studies {\bf 120} \publ Princeton University Press \yr 1989 
\endref 
\ref\no 2 
\by C.J. Bushnell and A. Fr\"ohlich 
\book Gauss sums and $p$-adic division algebras 
\bookinfo Lecture Notes in Math. \vol 987 \publ Springer \publaddr Berlin-Heidelberg New York  \yr 1983 
\endref 
\ref\no 3 
\by C.J. Bushnell and G. Henniart 
\paper Local tame lifting for $\roman{GL}(n)$ I: simple characters 
\jour Publ. Math. IHES \vol 83 \yr 1996 \pages 105--233 
\endref 
\ref\no4 
\bysame 
\paper Local tame lifting for $\roman{GL}(n)$ II: wildly ramified supercuspidals 
\jour Ast\'erisque \vol254 \yr 1999 
\endref 
\ref\no 5 
\bysame 
\paper Davenport-Hasse relations and an explicit Langlands correspondence 
\jour J. reine angew. Math. \vol 519 \yr 2000 \pages 171--199 
\endref 
\ref\no 6 
\bysame 
\paper Davenport-Hasse relations and an explicit Langlands correspondence II: twisting conjectures 
\jour J. Th. Nombres Bordeaux \vol 12  \yr 2000 \pages 309--347 
\endref 
\ref\no 7
\bysame 
\paper Local Rankin-Selberg convolution for $\roman{GL}(n)$: divisibility of the conductor 
\jour Math. Ann. \vol 321 \yr 2001 \pages 455--461 
\endref 
\ref\no 8 
\bysame 
\paper Local tame lifting for $\roman{GL}(n)$ IV: simple characters and base change 
\jour Proc. London Math. Soc. (3) \vol 87 \yr 2003 \pages 337--362 
\endref 
\ref\no 9 
\bysame 
\paper Local tame lifting for $\roman{GL}(n)$ III: explicit base change and Jacquet-Langlands correspondence 
\jour J. reine angew. Math. \vol 508 \yr 2005 \pages 39--100 
\endref 
\ref\no 10 
\bysame 
\paper The essentially tame local Langlands correspondence, I 
\jour J. Amer. Math. Soc. \vol 18 \yr 2005 \pages 685--710 
\endref 
\ref\no 11 
\bysame 
\paper The essentially tame local Langlands correspondence, II: totally ramified representations 
\jour Compositio Mathematica \vol 141 \yr 2005 \pages 979--1011 
\endref 
\ref\no 12 
\bysame 
\book The local Langlands Conjecture for $\roman{GL}(2)$ 
\bookinfo Grundlehren der mathematischen Wissenschaften {\bf 335} \publ Springer \yr 2006 
\endref 
\ref\no 13 
\bysame 
\paper The essentially tame local Langlands correspondence, III: the general case 
\jour Proc. London Math. Soc. (3) \vol 101 \yr 2010 \pages 497-553 
\endref 
\ref\no 14 
\bysame 
\paper The essentially tame local Jacquet-Langlands correspondence 
\jour Pure App. Math. Quarterly \vol 7 \yr 2011 \pages 469--538 
\endref 
\ref\no15 
\bysame 
\paper Intertwining of simple characters in $\roman{GL}(n)$ 
\jour  Internat. Math. Res. Notices \yr 2012 \pages{\tt doi:10.1093/imrn/rns162} 
\endref 
\ref\no16 
\by C.J. Bushnell, G. Henniart and P.C. Kutzko 
\paper Local Rankin-Selberg convolutions for $\roman{GL}_n$: Explicit conductor formula 
\jour J. Amer. Math. Soc. \vol 11 \yr 1998 \pages 703--730
\endref 
\ref\no 17 
\bysame 
\paper Correspondance de Langlands locale pour $\roman{GL}(n)$ et conducteurs de paires 
\jour Ann. Scient. \'Ecole Norm. Sup. (4) \vol 31 \yr 1998 \pages 537--560
\endref 
\ref\no 18 
\by C.J. Bushnell and P.C. Kutzko 
\book The admissible dual of $GL(N)$ via compact open subgroups 
\bookinfo Annals of Math. Studies {\bf 129} \publ Princeton University Press \yr 1993 
\endref 
\ref\no19 
\by A. Fr\"ohlich and A.M. McEvett 
\paper Forms over semisimple algebras with involution 
\jour J. Alg. \vol 12 \yr 1969 \pages 105--113 
\endref 
\ref\no 20 
\bysame  
\paper The representation of groups by automorphisms of forms 
\jour J. Alg. \vol 12 \yr 1969 \pages 114--133 
\endref 
\ref\no 21 
\by G. Glauberman 
\paper Correspondences of characters for relatively prime operator groups 
\jour Can\-ad. J. Math. \vol 20 \yr 1968 \pages 1465--1488 
\endref  
\ref\no 22 
\by M. Harris and R. Taylor 
\book On the geometry and cohomology of some simple Shimura varieties 
\bookinfo Annals of Math. Studies {\bf 151} \publ Princeton University Press \yr 2001 
\endref 
\ref\no23 
\by G. Henniart 
\paper La Conjecture de Langlands pour $\roman{GL}(3)$ 
\jour M\'em. Soc. Math. France \vol 12 \yr 1984 
\endref 
\ref\no24 
\bysame 
\paper Correspondance de Langlands-Kazhdan explicite dans le cas non-ramifi\'e 
\jour Math. Nachr. \vol 158 \yr 1992 \pages 7--26 
\endref 
\ref\no25 
\bysame 
\paper Caract\'erisation de la correspondance de Langlands par les facteurs $\varepsilon$ de paires 
\jour Invent. Math. \vol 113 \yr 1993 \pages 339--350 
\endref 
\ref\no26
\bysame 
\paper Une preuve simple des conjectures locales de Langlands pour $\roman{GL}_n$ sur un corps $p$-adique 
\jour Invent. Math. \vol 139 \yr 2000 \pages 439--455
\endref 
\ref\no27 
\bysame 
\paper Une caract\'erisation de la correspondance de Langlands locale pour $\roman{GL}(n)$ 
\jour Bull. Soc. Math. France \vol 130 \yr 2002 \pages 587--602 
\endref 
\ref\no28
\by G. Henniart and R. Herb 
\paper Automorphic induction for $GL(n)$ (over local non-archimed\-ean fields) 
\jour Duke Math. J. \vol 78 \yr 1995 \pages 131--192 
\endref 
\ref\no29
\by G. Henniart and B. Lemaire 
\paper Formules de caract\`eres pour l'induction automorphe 
\jour J. reine angew. Math. \vol 645 \yr 2010 \pages 41--84 
\endref 
\ref\no 30 
\bysame 
\paper Changement de base et induction automorphe pour $\roman{GL}_n$ en caract\'eristique non nulle 
\jour M\'em. Soc. Math. France \vol 108 \yr 2010 
\endref  
\ref\no31 
\by H. Jacquet, I. Piatetski-Shapiro and J. Shalika 
\paper Rankin-Selberg convolutions 
\jour Amer. J. Math. \vol 105 \yr 1983 \pages 367--483 
\endref  
\ref\no32 
\by P.C. Kutzko 
\paper The Langlands conjecture for $\roman{GL}_2$ of a local field 
\jour Ann. Math. \vol 112 \yr 1980 \pages 381-412 
\endref 
\ref\no33 
\by P.C. Kutzko and A. Moy 
\paper On the local Langlands conjecture in prime dimension 
\jour Ann. Math. \vol 121 \yr 1985 \pages 495--517 
\endref 
\ref\no34 
\by G. Laumon, M. Rapoport and U. Stuhler 
\paper $\Cal D$-elliptic sheaves and the Langlands correspondence 
\jour Invent. Math. \vol 113 \yr 1993 \pages 217--338 
\endref 
\ref\no35 
\by C. M\oe glin 
\paper Sur la correspondance de Langlands-Kazhdan 
\jour J. Math. Pures et Appl. (9) \vol 69 \yr 1990 \pages 175--226 
\endref 
\ref\no36 
\by F. Shahidi  
\paper Fourier transforms of intertwining operators and Plancherel measures for $\roman{GL}(n)$ 
\jour Amer. J. Math. \vol 106 \yr 1984 \pages 67--111
\endref  
\endRefs
\enddocument